\documentclass[12pt,a4paper]{amsart}

\usepackage{mathrsfs}
\usepackage{todonotes}

\usepackage{amssymb, tikz}

\usepackage{hyperref}

\usepackage{amsmath}
\usepackage[utf8]{inputenc}

\newtheorem{theorem}{Theorem}[section]

\newtheorem{lemma}[theorem]{Lemma}

\newtheorem{conclusion}[theorem]{Conclusion}

\theoremstyle{definition}
\newtheorem{definition}[theorem]{Definition}

\newtheorem{discussion}[theorem]{Discussion}

\theoremstyle{remark}
\newtheorem{remark}[theorem]{Remark}
\newtheorem{question}[theorem]{Question}

\newcommand{\dc}{{\rm dc}}

\newcommand{\cd}{{\rm cd}}

\newcommand{\nd}{{\rm nd}}

\newcommand{\Av}{{\rm Av}}
\newcommand{\CH}{{\rm CH}}
\newcommand{\Ax}{{\rm Ax}}
\newcommand{\otp}{{\rm otp}}
\newcommand{\stat}{{\rm stat}}
\newcommand{\club}{{\rm club}}
\newcommand{\gen}{{\rm gen}}

\newcommand{\cor}{{\rm cor}}

\newcommand{\Ord}{{\rm Ord}}

\newcommand{\Gen}{{\rm Gen}}
\newcommand{\Lim}{{\rm Lim}}

\newcommand{\id}{{\rm id}}

\newcommand{\Rang}{{\rm Rang}}

\newcommand{\rest}{{\restriction}}
\newcommand{\dom}{{\rm dom}}

\newcommand{\wilog}{{\rm without loss of generality}}

\newcommand{\then}{{\underline{then}}}

\newcommand{\mn}{{\medskip\noindent}}
\newcommand{\sn}{{\smallskip\noindent}}
\newcommand{\bn}{{\bigskip\noindent}}

\newcommand{\bbR}{{\mathbb R}}

\newcommand{\cC}{{\mathscr C}}

\newcommand{\bbD}{{\mathbb D}}
\newcommand{\cD}{{\mathscr D}}

\newcommand{\cH}{{\mathscr H}}

\newcommand{\cF}{{\mathscr F}}
\newcommand{\cE}{{\mathscr E}}
\newcommand{\cJ}{{\mathscr J}}
\newcommand{\cG}{{\mathscr G}}
\newcommand{\cI}{{\mathscr I}}

\newcommand{\cM}{{\mathscr M}}

\newcommand{\bbP}{{\mathbb P}}
\newcommand{\cP}{{\mathscr P}}
\newcommand{\gp}{{\mathfrak p}}
\newcommand{\gy}{{\mathfrak y}}
\newcommand{\gq}{{\mathfrak q}}
\newcommand{\gx}{{\mathfrak x}}

\newcommand{\bbQ}{{\mathbb Q}}

\newcommand{\cT}{{\mathscr T}}

\newcommand{\cf}{{\rm cf}}

\newcount\skewfactor
\def\mathunderaccent#1#2 {\let\theaccent#1\skewfactor#2
\mathpalette\putaccentunder}
\def\putaccentunder#1#2{\oalign{$#1#2$\crcr\hidewidth
\vbox to.2ex{\hbox{$#1\skew\skewfactor\theaccent{}$}\vss}\hidewidth}}
\def\name{\mathunderaccent\tilde-3 }

\newenvironment{PROOF}[2][\proofname.]
   {\begin{proof}[#1]\renewcommand{\qedsymbol}{\textsquare$_{\rm #2}$}}
   {\end{proof}}

\begin{document}

\title[NNR Revisited]{NNR Revisited}

\author {Mohammad Golshani}
\address{School of Mathematics\\
 Institute for Research in Fundamental Sciences (IPM)\\
  P.O. Box:
19395-5746\\
 Tehran-Iran.}
\email{golshani.m@gmail.com}
\urladdr{http://math.ipm.ac.ir/~golshani/}

\author {Saharon Shelah}
\address{Einstein Institute of Mathematics\\
Edmond J. Safra Campus, Givat Ram\\
The Hebrew University of Jerusalem\\
Jerusalem, 91904, Israel\\
 and \\
 Department of Mathematics\\
 Hill Center - Busch Campus \\
 Rutgers, The State University of New Jersey \\
 110 Frelinghuysen Road \\
 Piscataway, NJ 08854-8019 USA}
\email{shelah@math.huji.ac.il}
\urladdr{http://shelah.logic.at}
\thanks{The first author's research has been supported by a grant from IPM (No. 1400030417). The second
author's research has been partially supported by Israel Science Foundation (ISF) grant
no: 1838/19.  The authors thanks Alice Leonhardt for the beautiful typing the first version of this paper. The second author thanks Todd Eisworth for many corrections on an earlier version of this paper.
This is publication number 656 of the second author.}






\begin{abstract}
We show if we use countable support iteration of
forcing notions not adding reals that satisfy
additional conditions, then the limit
forcing does not add reals.  As a result we prove that we can amalgamate two
earlier methods and prove the consistency with ZFC + GCH of two
statements gotten separately earlier: Souslin hypothesis and non-club guessing.  We also answer a question of Justin Moore by
proving the consistency of one further case of ``strong failure of club
guessing" with GCH.
\end{abstract}

\maketitle
\setcounter{section}{-1}

\tableofcontents

\section{Introduction} \label{intro}
One of the major problems in the theory of iterated forcing is to prove some preservation theorems. One example
is the preservation of cardinals or cofinalities. By Solovay-Tennenbaum \cite{solovay-tennenbaum}, finite support iteration of c.c.c. forcing notions
is c.c.c. and hence it preserves all cardinals and cofinalities. In the case of countable support iterations, Shelah \cite{Sh:b} proved that the countable
support iteration of proper forcing notions is again proper, and hence it preserves $\aleph_1$.

In this paper we are interested in ``not adding new reals'', abbreviated NNR, in the case of countable support  iteration of proper forcing notions. There is a lot of work in this direction, see for example  \cite[Ch.V,\S7,~ Ch.VIII,\S4,~ Ch.XVIII,\S1,\S2]{Sh:f},
\cite[\S3]{Sh:666} and the references there.
In this paper we present a more general preservation theorem, and as an application we prove the consistency with ZFC + GCH of Souslin's hypothesis with non-club guessing. We also
prove the consistency with GCH of further cases of ``strong failure of club
guessing", in particular, 
    we answer a question of  Moore.

There are several limitations in preserving NNR at limit stages of countable support iterations that we will discuss some
of them in  Section \ref{Obstacles about NNR preservation}. By results of Shelah, if, e.g. $\bold V = \bold L$, then the preservation
of NNR fails. Indeed, assuming $\bold V = \bold L$,  Shelah has built a countable support iteration of length $\omega^2$
of NNR forcing notions of size $\aleph_1$ such that the limit adds a new real (see \cite[XVIII. Lemma 1.1]{Sh:f}).
Justin Moore \cite{moore2013}, solving a problem from
\cite[\S3]{Sh:666}, proved that  the continuum hypothesis implies the negation of the forcing axiom for the class of completely proper forcing notions. Indeed he proved in ZFC + CH that some proper forcing
notion not adding reals and satisfying a ``strong form of the medicine
against weak diamond" has no generic, in fact, is a tree with no
branches.


The paper is organized as follows.
In Section \ref{preliminaries}, we present some definitions and results which are needed for the rest of the paper.   In Section \ref{Obstacles about NNR preservation}, we discuss some obstacles  about preservation
of NNR in countable support iterations and suggest some ideas about how to overcome them.

In Section \ref{preservation},
we present sufficient conditions for countable support iteration of proper
forcing notions, not to add reals.  For this we define ``reasonable parameters
${\gp}$" and we have two main demands.  One (clause (c) of Definition
\ref{a24}) is a weakening of ``$\alpha$-proper for every $\alpha <
\omega_1$".  This time it has the form (on $\bbQ_i$), ${\gp}$-proper which
informally says that: if ${\gp} \in N,Y \subseteq \{M \in N:
M \text{ appropriate}\}$ is $\alpha$-large, then for some
$(N,\bbQ_i)$-generic condition $q \ge p,q$ forces that
$\{M \in Y:M[\name{\bold G}_{\bbQ_i}] \cap \bold V = M\}$ is $\alpha$-large
(the meaning of $\alpha$-large depends on ${\gp}$).
The other main demand (clause (d) of Definition \ref{a24}) is
a ``weak diamond preventive".

We then show that $\alpha$-properness for $\alpha < \omega_1$ is sufficient
for the first main demand (in Lemma \ref{a48}(3)).  The demand on the games for
${\gp}$ helps to prove the preservation of ${\gp}$-properness.

The preservation theorem in   Section \ref{preservation} does not, for
standard ${\gp}$, cover shooting a club $C \subseteq \omega_1$ running
away for $C_\delta \subseteq \delta = \sup(C_\delta),C_\delta$ small.  For this we will use, in Section \ref{delayed}, $({\gp},\alpha,\beta)$-proper for enough
pairs $\alpha \le \beta < \ell g({\gp})$ (so starting from $\beta$-large
we get $\alpha$-large; for many $\alpha$ we can choose $\beta = \alpha$,
but during the inductive proof we pass through cases of $\alpha <
\beta$).
Here we introduce various definitions and basic facts needed.

In Section \ref{examples},
we present the natural forcing showing $\kappa = 2$ is
interesting (not only $\kappa = \aleph_0$) (from \cite[Ch.VIII,\S4]{Sh:b}).
We show that the natural forcing (see above) for running away from $C_\delta
\subseteq \delta$, of small order type (see \cite[Ch.XVIII,\S2]{Sh:f}) falls
under our framework for delayed properness.
We give examples: running away from
$\langle C_{\delta,0},C_{\delta,1}:\delta < \omega_1 \text{ limit} \rangle,
C_{\delta,0},C_{\delta,1}$ are disjoint closed subsets of $\delta$
with no restrictions on their order type so we
ask for $C,C \cap C_{\delta,0}$ or $C \cap C_{\delta,1}$ to be bounded in
$\delta$ and more.

In Section
\ref{second},
we give a sufficient condition for the limit forcing
not to add reals.  We here are weakening the demand ``${\gp}$-proper", using
$({\gp},\alpha,f(\alpha))$-proper instead of
$({\gp},\alpha,\alpha)$-proper, what we called
delayed properness.  The price is that here ${\gp}$ has length of
large cofinality, so essentially we catch our tails on a club of
it.  Also the results here cover the examples.

In Section \ref{spelling}, we derive some forcing axioms from our preservation theorems and give several examples
that fit into our axioms.

Finally in Section \ref{moore}, we  answer a question of Justin Moore, which is related to the failure of weak club guessing
at $\omega_1$ in the presence of $\CH$.

The results and methods in this paper are all due to the second author. The first author's contribution was to fill in some details
and to write the paper. 

\section{Some preliminaries} \label{preliminaries}
In this section we present some preliminaries that are needed for the rest of the paper. We assume familiarity with the theory of iterated forcing
and countable support iterations. For a forcing notion $\bbP$ and conditions $p, q \in \bbP,$ we say $q$ is stronger than $p$ if $q \geq p$.
\begin{definition}
\label{y12}
\begin{enumerate}
\item $\bbQ$ is $\alpha$-proper
  if whenever $\chi$ is large enough regular, $\bar N = \langle N_i:i \le \alpha \rangle$ is an
increasing and continuous chain of countable elementary submodels of $({\cH}(\chi),\in)$ with $\alpha, \bbQ \in N_0$
and
$\bar N \restriction (i+1) \in N_{i+1}$,
if $p \in \bbQ \cap N_0$,
then there is $q,~p \le q \in \bbQ$ such that $q$ is
$(N_i,\bbQ)$-generic for each $i \le \alpha$.

\item We say $\bbQ$ is $(< \omega_1)$-proper if $\bbQ$ is
$\alpha$-proper for any $\alpha < \omega_1$.

\item We say $\bbQ$ is $(<^+ \omega_1)$-proper if it satisfies clause (1) for any $\alpha < \omega_1$ even omitting ``$\alpha \in
N_0$''.
\end{enumerate}
\end{definition}
\begin{definition}
\label{genPN}
Suppose $\bbP$ is a forcing notion, $p \in \bbP$ and  $N$ is a model with $\bbP \in N$. Then
\begin{enumerate}
\item $\Gen(N, \bbP) =\{ \bold G \subseteq \bbP \cap N: \bold G \text{~is a~} \bbP \cap N\text{-generic filter over~}N          \}.$

\item $\Gen^+(N, \bbP) =\{ \bold G \in \Gen(N, \bbP): G$ has an upper bound in $\bbP\}$.

\item $\Gen(N, \bbP, p) =\{ \bold G \in \Gen(N, \bbP): p \in \bold G\}$.
\end{enumerate}
\end{definition}
One important notion that is useful in proofs for showing that certain countable support iterations do not add reals is Shelah's notion
of completeness system.
\begin{definition}
\label{completeness system} (\cite[Ch. V, Definition 5.2]{Sh:f})
A \emph{completeness system} for a forcing notion $\bbP$ is a function $\mathbb{D}$ such that the following statements hold:
\begin{enumerate}
\item For a sufficiently large $\theta,$ the domain of $\mathbb{D}$ consists of pairs $(N, p),$ where
$N \prec (H(\theta), \in)$ is countable, $\bbP \in N$ and $p \in \bbP \cap N,$

\item For every $(N, p) \in \dom(\mathbb{D}),$ $\mathbb{D}(N, p)$ is a collection of subsets of $\Gen(N, \bbP, p)$.
\end{enumerate}
\end{definition}
\begin{definition}
\label{lambda-completeness system} (\cite[Ch. V, Definition 5.2]{Sh:f})
Suppose $\kappa$ is a cardinal. We say  $\mathbb{D}$ is a $\kappa$-completeness system for  $\bbP$, if it is a  completeness system
for $\bbP$ and for every $(N, p) \in \dom(\mathbb{D}),$ the intersection of fewer than $1+\kappa$ elements of $\mathbb{D}(N, p)$ is nonempty.
\end{definition}
\begin{definition}
\label{simple-completeness system} (\cite[Ch. V, Definition 5.4]{Sh:f})
A completeness system $\mathbb{D}$ for  $\bbP$ is \emph{simple}
if there is a second order formula $\Psi$ such that $\mathbb{D}(N, p)=\{\mathcal{G}_X: X \subseteq N          \}$, where
\[
\mathcal{G}_X=\{\bold G \in \Gen(N, \bbP, p): (N, \in, \bbP \cap N)\models \Psi(\bold G, X)    \}.
\]
\end{definition}
\begin{definition}
\label{D-complete forcing} (\cite[Ch. V, Definition 5.3]{Sh:f})
Suppose $\mathbb{D}$ is a simple completeness system for $\bbP$. Then $\bbP$ is said to be $\mathbb{D}$-complete, if for  every
 $(N, p) \in \dom(\mathbb{D}),$ $\Gen^+(N, \bbP, p)$ contains an element of $\mathbb{D}(N, p)$.
\end{definition}
The next theorem of Shelah gives a sufficient condition  for a countable support iteration of
forcing notions to not add new reals.
\begin{theorem}
\label{pre22}
(\cite[Ch. VIII, Theorem 4.5]{Sh:f})
A countable support iteration of forcing notions which are $<\omega_1$-proper and $\mathbb{D}$-complete
with respect to a simple 2-completeness system does not introduce reals.
\end{theorem}

\begin{definition}
\label{eccdef} (\cite[Ch. VII, Definition 1.2]{Sh:f})
The forcing notion $\bbP$ satisfies the $\kappa$-e.c.c. ($\kappa$-extra chain condition), if there
is a binary relation $R$ on $\bbP$ such that:
\begin{itemize}
\item For any sequence $\langle p_i:  i<\kappa\rangle$ of elements of $\bbP,$ there are pressing down functions $f_n: \kappa \to \kappa$ (i.e., for all $\alpha<\kappa, f_n(\alpha)< 1+\alpha$) for $n<\omega$
such that for all $0<i, j < \kappa$,
\[
\bigwedge\limits_{n<\omega} ( f_n(i)=f_n(j)) \implies p_i R p_j.
\]

\item If $\langle p_i: i \leq \omega \rangle$ and $\langle q_i: i \leq \omega \rangle$ are increasing sequences in $\bbP$ and for all $n<\omega, p_n R q_n$, then there is an $r$ such that for all $n<\omega, r \geq p_n, q_n$.
\end{itemize}
\end{definition}
\begin{definition}
\label{picdef} (\cite[Ch. VIII, Definition 2.1]{Sh:f})
The forcing notion $\bbP$  satisfies  the $\kappa$-pic ($\kappa$-properness isomorphism condition), if the following holds for any large enough regular cardinal $\lambda$: Suppose $i< j<\kappa,  N_i, N_j \prec (\cH(\lambda), \in, \lhd_\lambda)$  (where $\lhd_\lambda$ is a well-ordering of $\cH(\lambda)$)
are countable such that $\kappa, \bbP \in N_i \cap N_j$, $i\in N_i, j \in N_j, N_i \cap \kappa \subseteq j, N_i \cap i = N_j \cap j, p \in N_i \cap \bbP$ and $h: N_i \cong N_j$ is such that $h \restriction N_i \cap N_j$ is identity and $h(i)=j.$
Then there exists $q \in \bbP$ such that:
\begin{itemize}
\item $q \geq p, h(p)$ and for every maximal antichain $\cI \in N_i$ of $\bbP$, we have that $\cI \cap N_i$
is predense above $q$ and similarly for $\cI \in N_j$,

\item for every $r \in N_i \cap \bbP$ and $q'\geq q,$ there is $q'' \geq q'$ such that $$r \leq q'' \iff h(r) \leq q''.$$
\end{itemize}
\end{definition}
See \cite[Ch.VII, \S1]{Sh:f} and  \cite[Ch.VIII, \S2]{Sh:f}
for more information about the above two defined notions.
\begin{theorem}
\label{ecc and pic and chain condition}
Assume CH holds.
\begin{enumerate}
\item If $\bbP$ is a countable support iteration of length at most $\omega_2$ whose iterands are $<\omega_1$-proper, $\mathbb{D}$-complete
  for some $\aleph_1$-completeness system from $\bold V$ and satisfy the $\aleph_2$-e.c.c, then $\bbP$
satisfies the $\aleph_2$-c.c. The same result holds if we replace ``$\aleph_1$-completeness system'' by ``$\aleph_0$-completeness system''
or by ``2-completeness system''.

\item If $\bbP$ is a countable support iteration of length at most $\omega_2$ whose iterands satisfy the $\aleph_2$-pic, then $\bbP$
satisfies the $\aleph_2$-c.c.
\end{enumerate}
\end{theorem}
\begin{proof}
(1). For the case of $\aleph_1$-completeness system see \cite[Ch.VII, Lemmas 1.3]{Sh:f}. The case of $\aleph_0$-completeness system
follows from  \cite[Ch.VII, Lemmas 1.6]{Sh:f}  and the case of 2-completeness system follows from the above resuls combined with \cite[Ch.VIII, Theorem 4.5 and Lemma 4.13]{Sh:f}.

(2). See \cite[Ch.VIII, Lemma 2.4]{Sh:f}.
\end{proof}

\section{Obstacles for NNR preservation}
\label{Obstacles about NNR preservation}
In this section we give lengthy explanation of the problems and proofs for
NNR countable support iterations of proper forcing notions, and suggest some ideas about how to overcome them. These ideas will be made precise in the later sections of the paper.

\begin{definition}
\label{y3}
\begin{enumerate}
\item  Let $K_0$ be the family of countable support
iterations $\bar{\bbQ} = \langle \bbP_i,\name{\bbQ}_i:i < \alpha
\rangle$. We denote $\bbP_\alpha = \Lim(\bar{\bbQ})$.

\item We say $\bar{\bbQ} \in K_0$ is proper if for each $i<\alpha, ~\Vdash_{\bbP_i}\text{``}\name{\bbQ}_i$ is proper''. Note that it follows that
 $\bbP_j/\name{G}_{\bbP_i}$ is proper for $i < j \le \alpha$
(see \cite{Sh:b} or \cite{Sh:f}).

\item  We say $\bar{\bbQ} \in K_0$ is ${}^\omega \omega$-bounding if for
 each $i<\alpha, ~\Vdash_{\bbP_i}\text{``}\name{\bbQ}_i$ is ${}^\omega \omega$-bounding''
\footnote{The forcing notion $\bbQ$ is ${}^\omega \omega$-bounding  (in the
universe $\bold V$) if every $f \in ({}^\omega \omega)^{\bold V^{\bbQ}}$
is bounded by some $g \in ({}^\omega \omega)^{\bold V}$.}. It again follows that
 $\bbP_j/\dot{G}_{\bbP_i}$ is ${}^\omega \omega$-bounding for
$i < j \le \alpha$.

\item We say $\bar{\bbQ}$ is NNR if for each $i < \alpha,$ the forcing notion $\bbP_{i+1}$
adds no reals, or
equivalently if for $i < \alpha, \Vdash_{\bbP_i}\text{``} \name{\bbQ}_i$ adds no reals''  and for each $\beta <
\alpha, \bbP_\beta$ adds not reals.
\end{enumerate}
\end{definition}

It would be nice if also NNR is preserved in
limit stages of the iteration.  But this is wrong for at least two known  reasons, explained below:
\begin{enumerate}
\item[$\bigotimes_1$]   weak diamond
\sn
\item[$\bigotimes_2$]  existence of clubs.
\end{enumerate}
  Let us explain these obstacles in more details and the way to avoid them.

\underline{{\bf Weak diamond}}:

Let us first explain the obstacle arising from weak diamond.
Given a stationary set $S \subseteq \omega_1,$ recall that the weak diamond
$\Phi_S$ says: for each function $F:$$^{<\omega_1}2 \to 2$, there exists $g: \omega_1 \to 2$ such that for each
$f:\omega_1 \to 2$, the set
\[
\{ \delta \in S: g(\delta)=F(f \restriction \delta)        \}
\]
is stationary. By Devlin-Shelah
 \cite{DvSh:65} (see also \cite[Ch.XII,\S1]{Sh:b} or
\cite[AP,\S1]{Sh:f}), $\Phi_{\omega_1}$ is equivalent to $2^{\aleph_0} < 2^{\aleph_1}$.

Now let $\bar \eta = \langle \eta_\delta:\delta < \omega_1,\delta \text{ limit}
\rangle$ be a ladder system, where $\eta_\delta = \langle \eta_\delta(n):n < \omega \rangle$ is
 an increasing $\omega$-sequence of ordinals cofinal in $\delta$.  Let $D$ be a non-principal ultrafilter on $\omega$.
For $f \in {}^{\omega_1}2$ and a limit ordinal $\delta < \omega_1$  let
\[
\text{Av}_D(f,\eta_\delta)
= \ell \iff \{n:f(\eta_\delta(n)) = \ell\} \in D.
\]
Consider the following natural question:
\begin{question}
\label{y6}
(CH)  Given $\bar e = \langle e_\delta:\delta < \omega_1, \delta
\text{ limit} \rangle,e_\delta \in \{0,1\}$, is there
$f \in {}^{\omega_1}2$ such that for a club
of $\delta < \omega_1$ we have $e_\delta = \Av_D(f,\eta_\delta)$?
\end{question}
Naturally, trying to prove the consistency of this statement, we should use a countable support iteration
$\bar{\bbQ} =
\langle \bbP_i,\name{\bbQ}_i:i < \omega_2 \rangle$, where for each $i<\omega_2$, for some $\bar{e} \in V^{\bbP_i}$ as in Question \ref{y6}, $\bbQ_i$ is defined in
 $V^{\bbP_i}$ as $\bbQ_i = \bbQ_{\bar e}$, where

\[\begin{array}{r}
\bbQ_{\bar e} = \{f: \text{ for some }
\zeta < \omega_1,f \in {}^\zeta 2 \text{ and for every limit ordinal}\\
\delta \le \zeta
  \text{ we have } \Av_D(f,\eta_\delta) = e_\delta\}.
\end{array}\]
This is a very nice forcing notion,
 it is proper (even $< \omega_1$-proper,
see below) and NNR.
For example, let us show that   $\bbQ_{\bar e}$ is NNR.

Thus let $p \in \bbQ_{\bar e}$,  $\name\tau \in V^{\bbQ_{\bar e}}$ and $p\Vdash\text{``}\name\tau:\omega \rightarrow \Ord$
is a function''.
Let $\chi$ be a large enough regular cardinal and let $\bar N=\langle N_i:i \le \omega^2 \rangle$ be an increasing and continuous chain of countable
elementary submodels of $({\cH}(\chi), \in)$ with
${\bar e},\bbQ_{\bar e}, p, \name\tau \in N_0$ and $\bar N \rest (i+1)
\in N_{i+1}$. Let $\delta(i) = N_i \cap \omega_1$.  So $\langle
\delta(i):i \le \omega^2 \rangle$ is a strictly increasing and continuous sequence of countable ordinals.
Since $\eta_{\delta(\omega^2)}$  has
order type $\omega$,
$$W = \{i < \omega^2:\exists n(\delta_i \le
\eta_{\delta(\omega^2)}(n) < \delta_{i+1})\}$$
 has order type $\omega$ as well.
So, for each $n<\omega,$ we can find $\ell_n < \omega$ such that $$\bigwedge\limits_{n < \omega} \,
\bigwedge\limits_{m < \omega} \omega \cdot n + \ell_n +m \notin W.$$
 We choose, by induction on $n<\omega,~p_n \in
\bbQ_{\bar e}$ and $a_n \in \Ord$ such that:
\begin{itemize}
\item $p \le p_n,$
\item $p_{n-1} \le p_n,$
\item $p_n \in
N_{\omega \cdot n + \ell_n +1},$
\item $p_n$ forces some value for
$\name \tau(n)$, say $p_n \Vdash \name\tau(n)=\check{a}_n$
\item On $[\delta_{\omega \cdot m + \ell_m},
\delta_{\omega \cdot (m+1)+ \ell_{m+1}}) \cap \Rang(\eta_\delta) \backslash
\dom(p)$, $p_n$  agrees with $e_{\delta(\omega^2)}$.
\end{itemize}
This is easily seen to be possible by the choice of $\ell_n$'s. Then $q=\bigcup\limits_{n<\omega}p_n$ is a condition in
$\bbQ_{\bar e}$ and it forces $\name\tau= \langle \check{a}_n: n<\omega   \rangle \in \check{V}.$

Also note that for every $\alpha < \omega_1, {\cI}_\alpha =
\{f \in \bbQ_{\bar e}:\alpha \subseteq \dom(f)\}$ is a dense open
subset of $\bbQ_{\bar e}$, and hence if $\bold G$ is $\bbQ_{\bar e}$-generic over $V$, then $f=\bigcup\limits_{f \in \bold G}f: \omega_1 \to 2$
is as requested in Question \ref{y6}, for $\bar e$.
But clearly the weak diamond  tells us for this case that the answer is no,
that is:
$$\exists \bar e \forall f \in {}^{\omega_1}2 \exists^{\stat} \delta
(e_\delta \ne \Av_D(f,\eta_\delta)).$$

In fact this holds for any function $\Av':\bigcup\limits_{\delta < \omega_1}
{}^\delta 2 \rightarrow \{0,1\}$.
So
if $\bar{\bbQ}$ is going to preserve NNR, the desired demand on $\bbQ_i$'s
 should exclude the $\bbQ_{\bar e}$'s.
We now explain a way to overcome the above difficulty. Let us first give a definition.

\begin{definition}
\label{y9}
\begin{enumerate}
\item Let $K_1$ be the class of proper ${}^\omega \omega$-bounding
iterations $\bar{\bbQ} \in K_0$.
\item Let $K_2$ be the class of NNR
iterations $\bar{\bbQ} \in K_1$.

\item  Let $K_3$ be the class of
$\bar{\bbQ} \in K_2$ such that
if
\begin{enumerate}
\item  $ \chi$ is a large enough regular cardinal,

\item $ N \prec ({\cH}(\chi),\in)$ is countable,

\item $\bar{\bbQ} \in N$,

\item  $ i \in \ell g(\bar{\bbQ}) \cap N$,

\item $p \in \bbP_{i+1} \cap N$,

\item  $q_0,q_1 \in \bbP_i$ are $(N,\bbP_i)$-generic
(i.e. $q_\ell \Vdash ``N[\name{\bold G}_{\bbP_i}] \cap \bold V = N"$),

\item  $q_\ell \Vdash ``\name{\bold G}_{\bbP_i} \cap
N = \bold G^*"$,

\item $ p \restriction i \le q_\ell$,
\end{enumerate}
then
we can find $q'_0,q'_1$ and $\bold G^{**}$ such that for
 $\ell=1,2$ we have
 \begin{enumerate}

\item[(i)]  $q_\ell \le q'_\ell \in \bbP_{i+1}$,

\item[(j)]  $p  \le q'_\ell$,

\item[(k)]  $q'_\ell \Vdash ``\name{\bold G}_{\bbP_{i+1}} \cap N
= \bold G^{**}"$,

\item[(l)]  $q'_\ell$ is $(N,\bbP_{i+1})$-generic, so
$\bold G^{**} \subseteq \bbP_{i+1} \cap N$ is generic over $N$.
\end{enumerate}
\end{enumerate}
\end{definition}
\mn
Clause (3) of the above definition tries
to say the following.
We know $\name{\bold G}_{\bbP_i} \cap N$ (as being $\bold G^*$) and we
are looking at $N[\bold G^*]$ (formally, only its isomorphism type).
So we know $\name{\bbQ}^N_i[\bold G^*]$.
We would like to find $\bold G' \subseteq \name{\bbQ}^N_i[\bold G^*]$
generic over $N[\bold G^*]$, so that $\bold G^*,\bold G'$ will
determine $\bold G^{**}$.  But we need a guarantee that $\bold G'$
will have an upper bound in $\name{\bbQ}_i[\bold G_{\bbP_i}]$.
If we know $\bold G_{\bbP_i}$, this is  fine; but in a sense, we are given
2 candidates by $q_0,q_1$ and can increase them to
$q'_0 \rest i,q'_1 \restriction i$, and have to find $\bold G'$
``accepted" by both.

The weak diamond obstacle was overcame in \cite[Ch.V,\S7]{Sh:f} using $\aleph_1$-completeness systems
and in \cite[Ch.XVIII,\S4]{Sh:f} using 2-completeness systems. Here we show that
being in $K_3$ is sufficient to overcome this difficulty.
Indeed, we will assume something like the following.
Many times in some sense $q_0,q_1 \in \bbP_i$ are
$(N,\bbP_i)$-generic, $p \in \name{\bbQ}_i \cap N,q_\ell
\Vdash_{\bbP_i} ``\name{\bold G}_{\bbP_i} \cap N = \bold G^*"$
and for some $\bold G',q'_0 \ge q_0,q'_1 \ge q_1$ in $\bbP_{i+1}$ we have
$\bold G'\subseteq (\name{\bbQ}_i \cap N)[\bold G^*]$ and
$q'_\ell \Vdash_{\bbP_i} ``\name{\bold G}_{\bbQ_i} \cap N[\bold G^*]
= \bold G'"$ and $p \in \bold G'$.


Unfortunately, this is not sufficient to overcome with the other obstacle $\bigotimes_2$. There is an example where
 for some incomparable $q_0$ and $q_1$ in $\bbQ_i, \name E_i$ a $\bbQ_i$-name of
a club and for some $\alpha(q_0,q_1)$ we have:
\[
q_\ell \le q'_\ell
\]
\[
q'_\ell \Vdash ``\name E_i \cap \delta = E^\delta_i \Rightarrow
E^{\delta_0}_i \cap E^{\delta_1}_1 \backslash \alpha(q_0,q_1)
\text{ is finite}".
\]
This leads us to suggest another idea to overcome the obstacle which arises from the existence of clubs.

\underline{{\bf Existence of clubs}:}

Let us now explain the obstacle that arises from working with clubs.
This problem was already overcame  either by using $(< \omega_1)$-properness
or by a kind of ``finite powers are proper''.

As is shown in \cite{Sh:b}, \cite{Sh:f}, $(< \omega_1)$-properness is an antidote to such problems, i.e. against
$\bigotimes_2$. Though this is fine for many applications, like
 specializing an Aronszajn tree and many others,
but this requirement is too strong (see \cite{Sh:177}). For example consider the following
question.
\begin{question}
\label{y15}
Let $\bar C = \langle C_\delta:\delta < \omega_1,\delta \text{ limit}
\rangle$ where $C_\delta \subseteq \delta = \sup(C_\delta),\otp(C_\delta)
= \omega$ or at least $< \delta$. Is there a club $E$ of $\omega_1$
such that for each limit ordinal $\delta < \omega_1,~$ we have $\delta >
\sup(C_\delta \cap E)$ (i.e. is this consistent with CH)?
\end{question}
The natural forcing notion for adding such a club is given by

\[\begin{array}{r}
\bbQ^1_{\bar C} = \{f:\text{ for some non-limit } \alpha < \omega_1
\text{ we have } f \in {}^\alpha 2,f^{-1}(\{1\})\\
 \text{is closed and }
   \delta < \alpha \text{ limit } \Rightarrow
\sup(f^{-1}(\{1\}) \cap C_\delta) < \delta\}.
\end{array}\]
This forcing notion is
not even $\omega$-proper, for if $\langle N_i:i \le \omega \rangle$ satisfies
$C_{N_\omega \cap \omega_1} = \{N_i \cap \omega_1:i < \omega\}$,
then  no $f \in \bbQ^1_{\bar c}$ is $(N_i,\bbQ^1_{\bar C})$-generic, for infinitely
many $i$'s.

A solution to this problem was suggested in \cite[Ch.XVIII,\S2]{Sh:f}) by demanding that
each $\bbP_i \times \bbP_i$ is proper for $i < \ell g(\bar{\bbQ})$.
While this is fine for $\bbQ^1_{\bar C}$, this seems to exclude
specializing an Aronszajn tree without adding reals.

 In the proofs, we usually arrive to  a situation as follow:
\begin{itemize}
\item   $ \bar{\bbQ} \in N_0 \in N$,

\item  $ N_0 \prec ({\cH}(\chi),\in)$ and
$N \prec ({\cH}(\chi),\in)$ are countable,

\item  $q_\ell$ is $(N,\bbP_i)$-generic and
$(N_0,\bbP_i)$-generic  (for $\ell < 2$),

\item  $ q_\ell$ forces that $\name{\bold G}_{\bbP_i} \cap
N = \bold G_\ell$ (for $\ell < 2$),

\item  $ \bold G^* = \bold G_1 \cap N_0 =
\bold G_2 \cap N_0$,

\item  $i,j,p \in N_0[\bold G^*],i \le j \le
\ell g(\bar{\bbQ})$,

\item  $p \in P_j,p \restriction i \in \bold G^*$
(and possibly more).
\end{itemize}

We would like to find $\bold G' \subseteq \bbP^N_j/\bold G^*$ generic
over $N_0$ such that $q_0$ and $q_1$ both force that it has an upper
bound in $P_j/\name{\bold G}_{\bbP_i}$.  If $j=i+1$ this means
$\bold G' \subseteq \name{\bbQ}_i[\bold G^*]$ is generic over
$N_0$ such that $q_0,q_1$ both force that $\bold G'$ has an upper bound in
$\name{\bbQ}_i[\bold G_{\bbP_i}]$.

It is natural to demand $\bold G' \in N$, as otherwise the two possible generic
extensions (for $q_0$ and $q_1$) become unrelated.  For the case $j=i+1$,
the medicine against $\bigotimes_1$ should help us.  But we need it for every $j$.
Naturally we prove it by induction on $j$, and the successor case can be
reduced to the case $j=i+1$.

But to continue to a limit case, we need $\bold G' \in N$ and more:
for some intermediate $N_1$ with $N_0 \in N_1 \in N$, we also need
$\bigwedge\limits_{\ell} [q_\ell \Vdash N_1[\name{\bold G}_{\bbP_i}]
\cap \bold V = N_1]$.
So the clubs of elementary submodels which $q_0,q_1$ induce on
$\{M \prec N:M \in N\}$ should have non-trivial intersection.  This is a
major point and it has always appeared in some form.  Here the medicine
against $\bigotimes_2$ should help, in some way there will be many possible
$N_1$'s; but its help has a price, that is we have to carry it during the induction.
On the other hand the models playing the role of $N_1$ may change, we may
``consume it and discard it".

Note that the discussion is on two levels. Necessary limitations of universes
with $\CH$ on the one hand, and
how we try to carry the inductive proof
on appropriate iterations on the other hand; the connection though is
quite tight.

So we shall try for $j \in \ell g(\bar{\bbQ}) \cap N_0$ to extend the
situation with $i$ being replaced by $j$ while $\bold G^*$ is being
increased to $\bold G^{**}$.  We shall prove by induction some suitable facts,
with $\bold G^{**}$ the object we are really interested in.
We are given $q_1,q_2 \in \bbP_i$ and would like to find suitable
$q'_1,q'_2 \in \bbP_j$ such that $q'_\ell \restriction i = q_\ell$.
This last requirement helps us in limit steps to find an upper bound.

So the real action occurs for $j$ limit, hence we choose $\zeta_n \in N \cap
[i,j)$ such that $\zeta_0=i,\zeta_n < \zeta_{n+1}$ (sometimes better to have
$i$ and each $\zeta_n$ non-limit) and
$\bigcup\limits_{n < \omega} \zeta_n = \sup(j \cap N)$.
You can think of:
\begin{enumerate}
\item[]  in each case of limit $j$, proving the inductive statement,
we choose a

\item[] ``surrogate" for $N$ called $N_1$, during the induction it
serves like $N$, in

\item[] the limit dealing with $\zeta_0,\zeta_1,\ldots$ using
the induction hypothesis on $N_1$

\item[] we get $\bold G^{**}$ which may not
be in $N_1$ but is in $N$.
\end{enumerate}
So we try to choose by induction on $n,$ the conditions $q_{0,n},q_{1,n}$ and $\bold G^*_n$
such that:
\begin{itemize}
\item $q_{\ell,n} \in \bbP_{\zeta_n}$ is
$(N,\bbP_{\zeta_n})$-generic,
\item  $q_{\ell,0} =
q_\ell,$
\item $q_{\ell,n+1} \restriction \zeta_n = q_{\ell,n},$
\item $\bold G^*_n \in N_1,$
\item $\bold G^*_n \subseteq P_{\zeta_n} \cap N$ is generic over $N$, and
\item $q_{\ell,n} \Vdash ``\name{\bold G}_{\bbP_{\zeta_n}} \cap N = \bold G^*_n$.
\end{itemize}
The construction of the $\bold G^*_n$ should
use little information on the actual
$q_{\ell,n}$ so that the choices of the $\bold G^*_n$ can be carried say inside
$N_1$ so that $\langle \bold G^*_n:n < \omega \rangle \in N$.  In fact several
models will play a role like $N_1$.

By the proof of the preservation of ${}^\omega \omega$-bounding we can choose
some $N_1$ and demand ``$q_{\ell,n}$ gives to each $\bbP_{\zeta_n}$-name of an
ordinal $\name\tau_n \in N_1$, only finitely many possibilities".

Let us now explain how $(< \omega_1)$-properness or remaining proper under products can help in such arguments.
If the forcings are $(< \omega_1)$-proper , then we can assume in the beginning that
$\langle N_{1,\gamma}:\gamma \in A \rangle \in N$ is an increasing
and continuous chain of countable elementary submodels of some $(H(\chi), \in)$, $N_0 \prec N_{1,\gamma} \prec N,\langle N_{1,\gamma}:\gamma \le
\beta \rangle \in N_{\beta +1}$ with $A = (j+1) \cap N \backslash i$
and assume $q_\ell$ is $(N_{1,\gamma},\bbP_i)$-generic for
$\gamma \in A$ (similarly for $q'_0,q'_1,j$ in the conclusion) and
demand $q_{\ell,n}$ is $(N_{1,\gamma},\bbP_{\zeta_n})$-generic for
$n < \omega$ and $\gamma \in A \backslash \zeta_n$.


If components of the iteration remain proper under products, then
we demand things like ``$(q_0,q_1)$ is $(N_1,\bbP_i \times \bbP_i)$-generic"
so this gives many common $N_1$'s, but to preserve this we need more
complicated situations.
Instead of a ``tower" of models of countable length, we have a finite tower
of models  where on the bottom we are computing $\bold G^{**} \cap
\bbP_{\zeta_n}$ and as we go up, less and less is demanded.

In this paper, we will deal with a
condition which follows from both ``$(< \omega_1)$-properness'' and (essentially) ``the
square of the forcing notion is proper". We call this ${\gp}$-properness
where ``$\bbQ$ is ${\gp}$-proper" says that if $Y$ is a large
family of $M \prec N$ and if $p \in \bbQ \cap N$ and $\bbQ \in N$, then for
some $q$ we have $q \geq p$  is $(N,\bbQ)$-generic and
$q \Vdash ``\{M \in Y:M[\name{\bold G}_{\bbQ}] \cap \bold V = M\}$ is large".





\section {Preservation of not adding reals} \label{preservation}

In this section we define the notion of ${\gp}$-properness, for a reasonable parameter ${\gp}$, and prove
some preservation theorems.
\begin{definition}
\label{pseudo-filter}
\begin{enumerate}
\item A pseudo-filter on a set $N$ is a family  $D$ of subsets of $N$ which is closed under
supersets. If $D$ is a pseudo-filter on $N$, then we set
$D^-={\cP}(N)\setminus D$.

\item If $D$ is a filter on $N$, then set $D^+=\{X \subseteq N: N \setminus X \notin D   \}$.
\end{enumerate}
\end{definition}
\begin{definition}
\label{a3}
We say ${\gp} = (\bar \chi,\bar R,\bar{\cE},\bar D) = (\bar
\chi^{\gp},\bar R^{\gp},\bar{\cE}^{\gp},\bar D^{\gp})$ is a
reasonable parameter, when for some ordinal $\alpha^*$, denoted
$\ell g({\gp})$, we have:
\begin{enumerate}
\item[(a)]   $\bar \chi = \langle \chi_\alpha:\alpha < \alpha^* \rangle,$ where $
\chi_\alpha$ is a regular cardinal and ${\cH}((\bigcup\limits_{\beta < \alpha}
\chi_\beta)^+) \in {\cH}(\chi_\alpha)$.

\item[(b)]  $\bar R = \langle R_\alpha:\alpha < \alpha^* \rangle,$ where
$R_\alpha \in {\cH}(\chi_\alpha)$.

\item[(c)]   $\bar{\cE} = \langle {\cE}_\alpha:\alpha < \alpha^*
\rangle$, where ${\cE}_\alpha \subseteq [{\cH}(\chi_\alpha)]^{\le
\aleph_0}$ is stationary.

\item[(d)]   $\bar D = \langle D_\alpha:\alpha < \alpha^* \rangle,$
where
$D_\alpha$ is a function with domain ${\cE}_\alpha,$ and for $a \in {\cE}_\alpha,
D_\alpha(a)$ is a pseudo-filter on $a$.

\item[(e)]  for $\alpha < \alpha^*$ set ${\gp}^{[\alpha]} =:
\langle \bar \chi \restriction \alpha,\bar R \restriction (\alpha +1),
\bar{\cE} \restriction \alpha,\bar D \restriction \alpha \rangle$,
so it belongs to ${\cH}(\chi_\alpha)$.

\item[(f)]  if $a \in {\cE}_\alpha$, then for some countable
$N \prec ({\cH}(\chi_\alpha),\in)$,
$a$ is the universe of $N$, so we may write
$D_\alpha(N)$ instead of $D_\alpha(a)$ and $N \in {\cE}_\alpha$ instead of
$|N| \in {\cE}_\alpha$.

\item[(g)]  if $\alpha < \alpha^*$ and $N \in {\cE}_\alpha$, then
${\gp}^{[\alpha]} \in N$, so $\alpha \in N$.

\item[(h)]   for $N \in {\cE}_\alpha$ and $X \subseteq N$ we have:
\begin{center}
$X \in D_\alpha(N) \iff (\bigcup\limits_{\beta < \alpha}
{\cE}_\beta) \cap X  \in D_\alpha(N).$
\end{center}
\item[(i)]   if $N \in {\cE}_\alpha,X \in D_\alpha(N),\beta \in
\alpha \cap N$ and $y \in N \cap {\cH}(\chi_\beta)$,
then  for some $M \in {\cE}_\beta \cap X$ we have $X \cap M \in
D_\beta(M)$ and $y \in M$
\end{enumerate}
\end{definition}
Let us explain a little about the intended meaning of the above definition.
The requirement (a)
is just technical. About $R_\alpha$, we could require $R_\alpha$
is a relation on ${\cH}(\chi_\alpha)$, in a sense it codes a club of $[{\cH}(\chi_\alpha)]^{\le \aleph_0}$.
In clause (e), we considered $\bar R \restriction (\alpha +1)$ and not   $\bar R \restriction \alpha$. This makes it an easy demand on
${\cE}_\alpha,$ i.e., if $N \in {\cE}_\alpha$, then $R_\alpha \in N$. Clause (h) says that each $D_\alpha(N)$ has concentrated on  $\bigcup\limits_{\beta < \alpha} {\cE}_\beta$, and the last clause (i) is some kind of density, as it implies that $N \cap {\cH}(\chi_\beta) \subseteq \bigcup\{M: M \in   {\cE}_\beta \cap X         \}$.
\begin{remark}
\label{a6}
\begin{enumerate}
\item Note that $\langle {\cE}_\alpha:\alpha < \ell g({\gp}) \rangle$ are
pairwise disjoint by items (g) and  (e), so $D(N)$ can be well
defined as $D_\alpha(N)$ for the unique $\alpha$ such that $N \in
{\cE}_\alpha$.

\item Clearly, by clause (h), only $D_\alpha(N) \cap {\cP} (\bigcup\limits_{\beta < \alpha}
{\cE}_\beta)$ matters.
\end{enumerate}
\end{remark}
\begin{remark}
Some natural choices for $D(N)$ are as follows:
\begin{enumerate}
\item[$(a)$]  $D(N)$ is a filter on $N$.

\item[$(b)$]  $D(N)=\{X \subseteq N:X \ne
\emptyset \mod F\}$ for some filter $F$ on $N$.

\item[$(c)$]  $D(N) = F^+$ for a filter $F$ on $N$.
\end{enumerate}

\end{remark}
\begin{definition}
\label{a9}
Suppose ${\gp}=(\bar \chi,\bar R,\bar{\cE},\bar D)$ is a reasonable parameter as in Definition \ref{a3}.
\begin{enumerate}
\item We say $\bar D$ is standard, if for every $\alpha < \alpha^*
(= \ell g(\gp))$ and $N \in {\cE}_\alpha$ we have

\begin{equation*}
\begin{array}{clcr}
D_\alpha(N) = \{X \subseteq N:&\text{ for every } \gamma \in N \cap
\alpha \text{ and} \\
  &y \in N \cap \bigcup\{{\cH}(\chi_\beta):\beta \in N \cap \alpha\}, \\
  &\text{for some } \beta \in N \cap (\alpha \setminus \gamma) \text{ and} \\
  &M \in X \cap {\cE}_\beta, \text{ we have } y \in M \\
  &\text{and } X \cap M \in D_\beta(M)\}.
\end{array}
\end{equation*}

\item We say ${\gp}$ is standard if  $\bar D$ is standard.

\item We define the partial order $\le_{\gp}$ on $\alpha^* = \ell g(\gp)$ as
follows:  $\alpha \le_{\gp} \beta$ iff

\begin{enumerate}
\item[(a)]  $\alpha \le \beta,$

\item[(b)]  $N \in {\cE}_\beta \wedge \alpha \in N \Rightarrow N \cap
{\cH}(\chi_\alpha) \in {\cE}_\alpha$,

\item[(c)]  $N \in {\cE}_\beta \wedge \alpha \in N \wedge Y \in
D_\beta(N) \Rightarrow Y \cap \bigcup\limits_{\gamma < \alpha}
{\cE}_\gamma \in D_\alpha(M)$, where $M=N \cap
{\cH}(\chi_\alpha)$.
\end{enumerate}

\item We say ${\gp}$ is simple if $\alpha \le \beta < \alpha^* \Rightarrow
\alpha \le_{\gp} \beta$.

\item If $N \prec ({\cH}(\chi),\in)$ and $N \cap {\cH}(\chi_\alpha)
 \in {\cE}_\alpha$, (hence $\alpha, {\gp} \restriction
\alpha, R_\alpha \in N$), then  we let
$D_\alpha(N) = D^{\gp}_\alpha (N)$ to be $D_\alpha(N \cap
{\cH}(\chi_\alpha))$.
\end{enumerate}
\end{definition}
When ${\gp}$ is standard, we may drop $\bar D^{\gp}$ and just write ${\gp} = (\bar
\chi^{\gp},\bar R^{\gp},\bar{\cE}^{\gp})$.  Also if ${\gp}$ is
clear from the context, we may remove the superscript ${\gp}$.
We now define several games related to a reasonable parameter $\gp$.
\begin{definition}
\label{a15}
Suppose ${\gp}$ is a reasonable parameter.
\begin{enumerate}
\item For  $0< \alpha < \ell g({\gp})$ and $N \in
{\cE}^{\gp}_\alpha$, the game $\Game_\alpha(N,{\gp})$ is defined as follows. The play lasts $\omega$ moves, in the $n$-th move:
\begin{enumerate}

\item  the \underline{challenger}  chooses
$X_n \in D_\alpha(N)$ such that $m < n \Rightarrow X_n \subseteq X_m$

\item  the \underline{chooser}  chooses
$M_n \in X_n$ and $Y_n \subseteq M_n \cap X_n$ satisfying
$Y_n \in D(M_n) \cap N$

\item  the \underline{challenger} chooses $Z_n \subseteq Y_n$
such that $Z_n \in D(M_n)$.
\end{enumerate}
At the end, the chooser wins if $\bigcup\{\{M_n\} \cup Z_n:n < \omega\} \in
D_\alpha(N)$.

\item Assume $N \in N' \prec ({\cH}(\chi),\in)$, ${\gp} \restriction
\alpha \in N'$ and $N \prec N'$ are countable.  The game
$\Game'_\alpha(N,N',{\gp})$ is defined similar to $\Game_\alpha(N,\gp)$,
 but during the $n$-th
move, we demand that all the chosen objects belong to $N'$ (this
means only then $X_n \in N'$), and at the end of the $n$-th move,
the chooser also chooses $X'_n \subseteq X_n,X'_n \in D_\alpha(N) \cap N'$
and the challenger in the next move has to satisfy
$X_{n+1} \subseteq X'_n$.

\item   Omitting $N'$, i.e., writing $\Game'_\alpha(N,\gp)$ we mean:
for any  $N'$ as in (2), the demand $\Game'_\alpha(N,N',{\gp})$ holds.

\item We say that ${\gp}$ is a winner or a $\Game$-winner (resp. $\Game'$-winner), if
 for every $0< \alpha < \ell g({\gp})$ and $N \in
{\cE}^{\gp}_\alpha$, the chooser has a winning strategy in the game $
\Game_\alpha(N,{\gp})$ (resp. $\Game'_\alpha(N,{\gp})$).

\item We say that $\gp$ is a non-$\Game$-loser (resp. a non-$\Game'$-loser) if for $0<\alpha < \ell g({\gp})$ and $N \in {\cE}_\alpha$ the challenger has no winning strategy
in $\Game_\alpha(N,{\gp})$ (resp. $\Game'_\alpha(N,{\gp})$).


\end{enumerate}
\end{definition}

\begin{lemma}
\label{a18}
\begin{enumerate}

\item If ${\gp}$ is a reasonable parameter with the
standard $\bar D^{\gp}$, then  ${\gp}$ is a winner.

\item If ${\gp}$ is a $\Game_\alpha$-winner, then  ${\gp}$ is a
$\Game'_\alpha$-winner. If ${\gp}$ is a $\Game$-winner, then
${\gp}$ is a $\Game'$-winner. Similarly for a non-loser.
\end{enumerate}
\end{lemma}
\begin{proof}
(1). Suppose ${\gp}$ is a standard reasonable parameter. Let $0<\alpha < \ell g({\gp})$ and $N \in {\cE}^{\gp}_\alpha.$ Let $\langle y_n: n<\omega  \rangle$ be an enumeration of $N \cap \bigcup\{{\cH}(\chi_\beta): \beta \in \alpha \cap N             \}$ such that each $y \in N \cap \bigcup\{{\cH}(\chi_\beta): \beta \in \alpha \cap N             \}$ appears infinitely often in the enumeration and let $\langle \gamma_n: n<\omega  \rangle$ be an increasing sequence of ordinals in $N \cap \alpha$ with $\sup\limits_{n<\omega}\gamma_n = \sup(N \cap \alpha).$ We define the following winning strategy for chooser in the game $\Game_\alpha(N, {\gp})$: in the $n$-th move, the challenger chooses some $X_n \in D^{\gp}_\alpha(N)$. In particular, we can find $\beta_n \in N \cap \alpha \setminus \gamma_n$ and $M_n \in X_n \cap {\cE}^{\gp}_{\beta_n}$ such that
$ y_n \in M_n$ and $M_n \cap X_n \in D^{\gp}_{\beta_n}(M_n)$. Set also $Y_n=M_n \cap X_n.$ Then the challenger choose some
$Z_n \subseteq Y_n.$

We show that
$$X=\bigcup\{\{M_n\} \cup Z_n:n < \omega\} \in
D^{\gp}_\alpha(N).$$
Thus let $\gamma \in N \cap \alpha$ and $y \in N \cap \bigcup\{{\cH}(\chi_\beta): \beta \in \alpha \cap N             \}$. Pick $n<\omega$
such that $\gamma_n > \gamma$ and $y=y_n$. Then $\beta_n$ and $M_n$ are such that $\beta_n \in N \cap \alpha\setminus \gamma$,
$M_n \in X  \cap {\cE}^{\gp}_{\beta_n}$ and we have $y_n \in M_n$ and $X \cap M_n \in D_{\beta_n}^{\gp}(M_n)$.
 Thus $X \in D^{\gp}_\alpha(N)$. Hence the above process defines a winning strategy for chooser, as required.

(2). Suppose  ${\gp}$ is a $\Game_\alpha$-winner. Let  $N \in {\cE}_\alpha$ and assume that $N'$ is such that $N \in N' \prec ({\cH}(\chi),\in)$, ${\gp} \restriction
\alpha \in N'$ and $N'$ is countable. We define a winning strategy for chooser in the game  $\Game'_\alpha(N,N',{\gp})$.

  Let  $\sigma$ be a winning strategy for chooser in the game $\Game_\alpha(N, {\gp}).$ By elementarity, we may assume that $\sigma$ is in $N'$.
  We define the strategy $\sigma'$ for chooser in the game $\Game'_\alpha(N,N',{\gp})$ as follows. At the $n$-th move, the challenger chooses some
  $X_n \in N'$. Then chooser picks the sets $M_n$ and $Y_n$ via the strategy $\sigma$ and he also takes $X'_n$ to be $X_n.$ As $\sigma$
   is in $N'$, all these objects are also in $N'$. Then challenger chooses
 some $Z_n \in N'.$ It is evident that
  $\sigma'$ is a winning strategy for chooser in the game $\Game'_\alpha(N,N',{\gp})$,
  as required. The other cases of the lemma can be proved in a similar way.
\end{proof}

\begin{definition}
\label{a21}
Assume ${\gp}$ is a reasonable parameter, $\alpha < \ell g({\gp})$, $N \in
{\cE}^{\gp}_\alpha,~y \in N$ and $\bbP \in N$ is a forcing notion.  Set
$${\cM}_{\bbP}[\name{\bold G}_{\bbP},N,y] = \{M \in N:\bbP,y
\in M \text{ and } \name{\bold G}_{\bbP} \cap M
\text{ is } (\bbP \cap M)\text{-generic over } M\},$$
where $\name{\bold G}_{\bbP}$ is the canonical $\bbP$-name for the generic filter.
\end{definition}
We consider ${\cM}_{\bbP}[\name{\bold G}_{\bbP},N,y]$ as a $\bbP$-name and then $\cM_{\bbP}[\bold G,N,y]$ is well defined
for any $\bbP$-generic filter $\bold G$.
If $\bbP$ is clear from the context, we may omit it. Note that ${\cM}_{\bbP}
[\bold G,N,y] = {\cM}_{\bbP}[\bold G \cap N,N,y]$, so we may write
$\bold G \cap N$ instead of $\bold G$.  If $y = \emptyset$ we may omit it.

\begin{definition}
\label{a24}
We say $\bar{\bbQ} \in K_0$ is a ${\gp}-$NNR$^0_{\aleph_0}$ iteration if the following conditions are satisfied:
\mn
\begin{enumerate}
\item[(a)]   $\bar{\bbQ} = \langle \bbP_i,\name{\bbQ}_i:i < j(*)
\rangle$ is a countable support iteration of proper forcing notions such that
 $\bar{\bbQ}, {\cP}(\Lim(\bar{\bbQ})) \in {\cH}
(\chi^{\gp}_0)$.

\item[(b)]   forcing with $\bbP_{j(*)} = \Lim(\bar{\bbQ})$ does not add
reals.

\item[(c)] (long properness) suppose that:

\begin{enumerate}
\item[$(*)_1$]  $(\alpha) \quad i \le j \le j(*),\alpha < \ell
g({\gp})$,

\item[${{}}$]  $(\beta) \quad N \in {\cE}^{\gp}_\alpha,\{i,j,\bar{\bbQ}\}
\in N$,

\item[${{}}$]  $(\gamma) \quad$ the condition $q \in \bbP_i$
is $(N,\bbP_i)$-generic,

\item[${{}}$]  $(\delta) \quad q \Vdash$``$ \name{\bold G}_{\bbP_i} \cap
N = \bold G$'',

\item[${{}}$]  $(\varepsilon) \quad p \in \bbP_j \cap N$ and $p \rest i \in
\bold G$,

\item[${{}}$]  $(\zeta) \quad Y \subseteq {\cM}_{\bbP_i}[\bold G,N,y]$
 where $y = \langle \bar{\bbQ},i,j \rangle$ and $Y \in D^{\gp}_\alpha(N)$
\footnote{Note that
the ordinal $i$ is
reconstructible from $\bar{\bbQ}$ and $\bold G$.},
\end{enumerate}
then there are $\bold G',q'$ such that:
\begin{enumerate}

\item[$(*)_2$]  $(\eta) \quad q' \in \bbP_j, ~p \le q'$ and
$q \le q' \restriction i$,

\item[${{}}$]  $(\theta) \quad q'$ is $(N,\bbP_j)$-generic,

\item[${{}}$]  $(\iota) \quad q'\Vdash$``$\name{\bold G}_{\bbP_j} \cap
N = \bold G'$'',

\item[${{}}$]   $(\kappa) \quad Y \cap {\cM}_{\bbP_j}[\bold G',N,y]
\in D^{\gp}_\alpha(N)$.
\end{enumerate}

\item[(d)]  (anti weak diamond or anti-w.d.)  suppose that:
\begin{enumerate}

\item[$(\bullet)_1$] $(\alpha) \quad i \le j \le j(*)$ and $\alpha < \ell
g({\gp}),$

\item[${{}}$] $(\beta) \quad N_0 \in N_1 \in {\cE}^{\gp}_\alpha,~N_0 \in
\bigcup\limits_{\beta < \alpha} {\cE}^{\gp}_\beta,$

\item[${{}}$] $(\gamma) \quad \otp(N_0 \cap [i,j)) < \alpha$
\footnote{so naturally $\ell g({\gp}) = \omega_1$. We use
the parallel of ``$\aleph_0$-completeness system" rather than
``2-completeness system" of \cite{Sh:f} as things are complicated
enough anyhow; see Definition \ref{a39} and Theorem \ref{a42}.},

\item[${{}}$] $(\delta) \quad n < \omega$ and for $\ell < n$
we have $q_\ell \in \bbP_i$ is $(N_1,\bbP_i)$-generic,

\item[${{}}$] $(\epsilon) \quad q_\ell \Vdash$``
$\name{\bold G}_{\bbP_i} \cap N_1 = \bold G^\ell$'',

\item[${{}}$] $(\zeta) \quad \bigwedge\limits_{\ell <n}
[\bold G^\ell \cap N_0 = \bold G^*]$ where
$\bold G^* \subseteq \bbP_i \cap N_0$
 is generic over $N_0$,

 \item[${{}}$] $(\eta) \quad  Y =: \bigcap\limits_{\ell <n} {\cM}_{\bbP_i}[\bold G^\ell,N_1] \in
D^{\gp}_\alpha(N_1)$,

\item[${{}}$] $(\theta) \quad p \in \bbP_j \cap N_0$ is such that $p \restriction i \in \bold
G^*$.
\end{enumerate}
Then:
\begin{enumerate}
\item[$(\bullet)_2$]   for some $\bold G^{**} \subseteq
\bbP_j \cap N_0$ generic over $N_0$ we have $p \in \bold G^{**} \in N_1$ and
 $\bigwedge\limits_{\ell < n}
\,\,\bigvee\limits_{q \in G^\ell}[q \Vdash
``\bold G^{**}$ has an upper bound in $\bbP_j/\name{\bold G}_{\bbP_i}"]$.
\end{enumerate}
\end{enumerate}
\end{definition}

\begin{remark}
\label{a25}
We may like to phrase clause (c) as a condition on each
$\name{\bbQ}_i$, for this see Definitions \ref{b18},
\ref{b24} and \ref{b37}.
\end{remark}
We now state and prove the main result of this section.
\begin{theorem}
\label{a27}
Assume ${\gp}$ is a reasonable parameter
of length $\omega_1$, $\bar{\bbQ}$ is a countable support iteration such that $\bar{\bbQ}, {\cP}(\Lim \bar{\bbQ}) \in
{\cH}(\chi^{\gp}_0)$, $ \delta = \ell g(\bar{\bbQ})$ is a limit ordinal and
for every $\alpha < \delta,\bar{\bbQ} \restriction \alpha$ is
a ${\gp}-\text{NNR}^0_{\aleph_0}$ iteration and ${\gp}$ is a $\Game$-winner.
Then $\bar{\bbQ}$ is a ${\gp}-\text{NNR}^0_{\aleph_0}$ iteration.
\end{theorem}

\begin{proof}
We show that  items (a)-(d)  of Definition \ref{a24} are  satisfied by $\bar{\bbQ}$.

\noindent
\underline{Proof of clause (a)}:
This holds trivially by our assumptions.

\noindent
\underline{Proof of clause (b)}:
This follows from clause (d) of Definition \ref{a24} proved below. To see this, let $p \in \bbP_{\delta}=\Lim \bar{\bbQ},$
$\name{r}$ be a $\bbP_\delta$-name  and suppose that $p \Vdash$``$\name{r}$ is a real''.
In Definition \ref{a24}(d) set $i=0, j=\delta, n=1$ and pick $N_0 \in N_1$ so that the hypotheses in $(\bullet)_1$ are satisfied and $\name{r} \in N_0$. Thus, by $(\bullet)_2,$ we can find $\bold G^{**} \subseteq
\bbP_\delta \cap N_0$ which is generic over $N_0$, $p \in \bold G^{**}$  and
 $\bold G^{**}$ has an upper bound $q \geq p$ in $\bbP_\delta$.
By genericity of $\bold G^{**}$, $q$ decides $\name{r}$ to be a real in $V$,
and we are done.

We first prove clause (d) and then return to clause (c).

\noindent
\underline{Proof of clause (d)}:
Let $i, j, \alpha, N_0, N_1, n, q_0,\dotsc, q_{n-1}, G^0, \dotsc, G^{n-1}, G^*$ and $p$ be as in the
assumptions of clause (d).  Let $\alpha' = \otp(N_0 \cap [i,j))$. Then
$\alpha' < \omega_1$, so $\alpha' \in N_1$.
If $j < \delta$, then by our assumption
$\bar{\bbQ} \restriction j$ is a ${\gp}-$NNR$^0_{\aleph_0}$
iteration and the result follows. Thus  assume that $j = \delta$. If $i=j$, the conclusion is trivial, so
assume $i<j$.

Let $i_m \in N_0 \cap j$ be such that $i_0 = i$
and $\langle i_m:m < \omega \rangle \in N_1$ is an increasing sequence with
$\bigcup\limits_{m < \omega} i_m = \sup(N_0 \cap j)$.
Choose $\langle M_k: k < 5 \rangle$ such that:
\begin{itemize}
\item $y^* := \{i,j,\alpha,\alpha',
\bar{\bbQ},N_0,\langle i_m:m < \omega \rangle\} \in M_k$,

\item $M_k \in
{\cE}^{\gp}_{\alpha'} \cap N_1 \cap \bigcap\limits_{\ell < n}
{\cM}_{\bbP_i}[G^\ell,N_1]$,

\item  $M_0 \in M_1 \in M_2 \in M_3 \in M_4$,

\item  $\bigcap\limits_{\ell < n} {\cM}_{\bbP_i}
[\bold G^\ell,M_0,y^*] \in D^{\gp}_\alpha(M_0)$.

\end{itemize}
Note that for $\ell < n$ and $k<5$,  $N_0 \prec M_k
\prec N_1$ and  $\bold G^\ell \cap
M_k$ is a generic subset of $\bbP_i \cap M_k$.
Now for $\ell < n$ we can choose $q'_\ell \in \bold G^\ell \cap M_4$ such that:
\begin{itemize}
\item $q'_\ell$ forces (for $\bbP_{i_0} = \bbP_i$) a value for
$\name{\bold G}_{\bbP_{i_0}} \cap M_3$, which necessarily is
$\bold G^\ell \cap M_3$,

\item  $q_\ell \le
q'_\ell$,

\item $q'_\ell$ is $(M_k,\bbP_{i_0})$-generic, forcing
$\name{\bold G}_{\bbP_{i_0}} \cap M_k = \bold G^\ell \cap M_k$ for $k =
0,1,2,3$,

\item $q'_\ell$ forces $\name{\bold G}_{\bbP_{i_0}} \cap N_0 = \bold G^*$.
\end{itemize}
Let $\langle {\cI}^*_m:m < \omega \rangle \in M_0$ list
the maximal antichains of $\bbP_j$ that belongs to $N_0$.  We choose,
by induction on $m < \omega$, the objects $r_m, p_m, n_m, \bold G^*_m$,
 $\langle \bold G^\ell_m: \ell < n_m \rangle$ and $Y_m$ such that:

\begin{enumerate}
\item[$(*)_1$]  $(a) \quad r_m \in \bbP_{i_m} \cap M_4$,

\item[${{}}$]  $(b) \quad \dom(r_m) \subseteq [i,i_m)$,

\item[${{}}$]  $(c) \quad r_{m+1} \restriction i_m = r_m$,

\item[${{}}$]  $(d) \quad q'_\ell \cup r_m \in \bbP_{i_m}$ and is
$(M_k,\bbP_{i_m})$-generic for $k=0,1,2,3$ and is
$(N_0,\bbP_{i_m})$-
 generic,\footnote{note that $q'_\ell$ and $r_m$ have disjoint domains.}

\item[${{}}$]  $(e) \quad$ for every predense subset ${\cI}$
of $\bbP_{i_m}$ which belongs to $M_2$, for some
finite ${\cJ}
\subseteq {\cI} \cap M_2$, the set ${\cJ}$ is predense above
$q'_\ell \cup r_m,$
for each $\ell < n$,

\item[${{}}$]  $(f) \quad n_m < \omega$ and for $\ell < n_m$
$\bold G^\ell_m \in M_1$ is a subset of $\bbP_{i_m} \cap M_0$ generic
over $M_0,$

\item[${{}}$]  $(g) \quad$ if $\ell < n_{m+1}$, then
$\bold G^\ell_{m+1} \cap \bbP_{i_m} \in \{\bold G^k_m:k < n_m\}$,

\item[${{}}$]  $(h) \quad n_0 = n$ and $\bold G^\ell_0 = \bold G^\ell \cap M_0$,

\item[${{}}$]  $(i) \quad q'_\ell \cup r_m \Vdash_{\bbP_{i_m}}
``\name{\bold G}_{\bbP_{i_m}} \cap M_0 \in \{\bold G^\ell_m:\ell < n_m\}$'',

\item[${{}}$] $(j) \quad \bold G^*_m$ is a subset of $\bbP_{i_m}
\cap N_0$ generic over $N_0$,

\item[${{}}$]  $(k) \quad \bold G^*_m \subseteq \bold G^\ell_m$
for $\ell < n_m$,

\item[${{}}$]  $(l) \quad p_m$ is such that:

\item[${{}}$]  $\quad\quad (l$-$1) \quad p_m\in \bbP_j \cap N_0$,

\item[${{}}$]  $\quad\quad (l$-$2) \quad p_m \rest i_m
\in \bold G^*_m,$

\item[${{}}$]  $\quad\quad (l$-$3) \quad p_{m+1} \in {\cI}^*_m,$

\item[${{}}$]  $\quad\quad (l$-$4) \quad p_0 = p,$

\item[${{}}$]  $\quad\quad (l$-$5) \quad p_m \le p_{m +1}$.

\item[${{}}$]  $(m) \quad Y_m = \bigcap\limits_{\ell < n_m} {\cM_{\bbP_{i_m}}}
[G^\ell_m,M_0,y^*] \in D_{\alpha'}(M_0)$ where $y^* = \{N_0,\langle i_m:m < \omega \rangle,\bar{\bbQ},i,j\}$.
\end{enumerate}
The construction is clear for $m=0$.  So suppose that we have it for $m$ and we shall
choose for $m+1$.
We do this is several steps.

\noindent
\underline{Stage A}:  Choose $p_{m+1} \in N_0 \cap {\cI}^*_m$ such that
$p_m \le p_{m+1}$ and $p_{m+1} \restriction i_m \in \bold G^*_m$.

\noindent
\underline{Stage B}:  Choose $\bold G^*_{m+1} \subseteq
\bbP_{i_{m+1}} \cap N_0$ generic over $N_0$ such that $\bold G^*_m \subseteq
\bold G^*_{m+1} \in M_0, p_{m+1} \restriction i_{m+1} \in
\bold G^*_{m+1}$ and
\begin{center}
$\bigwedge\limits_{\ell < n_m} \,
\bigvee\limits_{r \in \bold G^\ell_m} [r \Vdash_{\bbP_{i_m}}
``\bold G^*_{m+1}$ has an upper bound in $\bbP_{i_{m+1}}/
\name{\bold G}_{\bbP_{i_m}}]$.
\end{center}

This is easy by applying clause (d) of the Definition \ref{a24} for
$i_m, i_{m+1}, \alpha', p_{m+1} \restriction i_m, \bold G^*_m, \langle
\bold G^\ell_m:\ell < n_m \rangle, N_0, M_0$,  for the forcing notion $\bar{\bbQ} \restriction i_{m+1}$, which is, by induction
hypothesis, a
${\gp}-NNR^0_{\aleph_0}$ iteration, .
We  also use the fact that $\otp(N_0 \cap [i_m,i_{m+1})) < \otp(N_0 \cap [i_m,j))
= \alpha'$.

\noindent
\underline{Stage C}:
As $\bbP_{i_m}$ is proper and adds no new  reals, $\name{\bold G}_{\bbP_{i_m}} \cap M_1$ is a $\bbP_{i_m}$-name of
an object from $\bold V$,
so
$${\cI} = \{p \in \bbP_{i_m}:p \text{ forces a value for }
\name{\bold G}_{\bbP_{i_m}} \cap M_1 \text{ in } \bold V\}$$
 is a dense
open subset of $\bbP_{i_m}$ and ${\cI} \in M_2$.  By clause $(*)_1(e)$ of the
induction hypothesis, there is a finite ${\cJ} \subseteq {\cI}
\cap M_2$ such that: $\ell < n \Rightarrow {\cJ}$ is predense above
$q'_\ell \cup r_m$.  Without loss of generality ${\cJ}$ is minimal.  Let
$n_{m+1} = |{\cJ}|$.

Let ${\cJ} = \{p^\ell_m:\ell < n_{m+1}\}$  and for each $\ell < n_{m+1}$ choose  $H^\ell_m \in
M_2$ such that $p^\ell_m \Vdash
``\name{\bold G}_{\bbP_{i_m}} \cap M_1 = H^\ell_m$''. As ${\cJ}$ is minimal, $H^\ell_m \cap M_0 \in
\{\bold G^\ell_m:\ell < n_m\}$ so for some  $h:n_{m+1} \rightarrow n_m$ and every $\ell < n_{m+1}$ we have $H^\ell_m \cap M_0 =
\bold G^{h(\ell)}_m$.

Let $Y = \bigcap\limits_{\ell < m_n} {\cM_{\bbP_{i_m}}}[\bold G^\ell_m,M_0,y^*] \in
D^{\gp}_{\alpha'}(M_0)$.
Now we choose by induction on $\ell \le n_{m+1}$ a condition
$r^\ell_m \in M_1$ such that:
\begin{enumerate}
\item[$(*)_2$]   $(\alpha) \quad r^\ell_m \in \bbP_{i_{m+1}} \cap M_1$
  and $r^\ell_m \restriction i_m \in H^\ell_m$,

\item[${{}}$]   $(\beta) \quad r^\ell_m$ is $(M_0,\bbP_{i_m})$-generic
and forces a value for $\name{\bold G}_{\bbP_{i_m}} \cap M_0$, say
 $\bold G^\ell_{m+1}$,

\item[${{}}$]  $(\gamma) \quad r^\ell_m$ is above $\bold G^*_{m+1}$,
and moreover above $p^{h(\ell)}_m$,

\item[${{}}$]   $(\delta) \quad Y \cap \bigcap\limits_{k < \ell} {\cM_{\bbP_{i_m}}}
[\bold G^k_{m+1},M_0,y^*] \in D^{\gp}_{\alpha'}(M_0)$.
\end{enumerate}
The construction can be easily done by applying clause (c) of Definition \ref{a24}
to $i_m, i_{m+1}, \alpha', M_0$, large
enough member of $H^\ell_m, p^{h(\ell)}_m$
and
$Y^\ell_m = Y \cap \bigcap\limits_{k < \ell}{\cM}[\bold G^k_{m+1},M_0,y^*] \in M_1,$
for the ${\gp}-NNR^0_{\aleph_0}$-iteration $\bar{\bbQ} \restriction i_{m+1}$.

\noindent
\underline{Stage D}:
By  \cite[Ch.XVIII, Claim 2.6]{Sh:f}, we can choose $r_{m+1}$ as required such that $\{r^\ell_m:\ell < n_{m+1}\}$
is predense over it.

This completes the inductive construction. Let us now show that this is sufficient to get clause (d).
Let $\lhd^* \in N_1$ be a well-ordering of $M_4$. During the construction above we chose inductively members of $M_4$ and
all the parameters used are from $M_4$, so if we  always choose the $\lhd^*$-first object, the
construction is determined and is in $N_1$. By clause $(*)_1(c)$, we have
\begin{enumerate}
\item[$(\alpha)$] Let   $r = \bigcup\limits_{m} r_m$ be  the unique
$r \in \bbP_j$ satisfying $m < \omega \Rightarrow r \restriction i_m
= r_m$. Then $r \in \bbP_j \cap N_1$.
\end{enumerate}
Also, by the choice of $\langle {\cI}^*_m:m < \omega \rangle$ and
clause $(*)_1(k)$, we have
\begin{enumerate}
\item[$(\beta)$]  $\bold G^{**} = \{p' \in \bbP_j \cap N_0:
\bigvee\limits_{m < \omega}[p' \le p_m]\}$ belongs to $M_4$ and is
a subset of $\bbP_j \cap N_0$ generic over $N_0$.
\end{enumerate}
It is also clear from $(*)_1(\ell)$ that
\begin{enumerate}
\item[$(\gamma)$]   $q'_\ell \cup r$ is above $\bold G^{**}$
(in $\bbP_j$).
\end{enumerate}
So we have finished proving clause (d).

\noindent
\underline{Proof of clause (c)}:
We prove this by induction on $\alpha$.
Let $i, j, \alpha, N, p, q$ and  $Y$ be as there.  If $j < \delta$ we can apply ``$\bar{\bbQ}
\restriction j$ is a ${\gp}-NNR^0_{\aleph_0}$ iteration'', so \wilog \,
$j =\delta$.  If $i=j$ the statement is trivial, so assume $i<j$.

Choose
$i_n \in N \cap j$, for $n < \omega$, such that $i_0 = i, i_n <
i_{n+1}$ and $\bigcup\limits_{n < \omega} i_n = \sup(j \cap N)$.
Let $\langle (y_n,\beta_n):n < \omega \rangle$ list the pairs $(y,\beta) \in
N \times (\alpha \cap N)$ such that $y \in {\cH}(\chi^{\gp}_\beta)$.
Let $\sigma$ be a winning strategy for the chooser in the game
$\Game_\alpha(N, \gp)$
and let $\langle {\cI}_n:n < \omega \rangle$ list the dense open subsets of
$\bbP_j$ which belong to $N$.

We choose by induction on $n < \omega$, the objects $q_n, \name p_n, \name M_n$
and $ \name Y_n$ such that:

\begin{enumerate}
\item[$(*)_3$]  $(a) \quad q_n \in \bbP_{i_n}$ with $q_0 = q$,

\item[${{}}$]  $(b) \quad q_n$ is $(N,\bbP_{i_n})$-generic,

\item[${{}}$]  $(c) \quad q_{n+1} \restriction i_n = q_n$,

\item[${{}}$]  $(d) \quad \name p_n$ is a $\bbP_{i_n}$-name of a member
of $(\bbP_j/\name{\bold G}_{\bbP_{i_n}}) \cap N$

\item[${{}}$]  $(e) \quad \name p_n$ is forced to belong to ${\cI}_n$,

\item[${{}}$]  $(f) \quad \name M_n$ is a $\bbP_{i_n}$-name of a
member of ${\cE}^{\gp}_{\beta_n} \cap N$,

\item[${{}}$]  $(g) \quad$if $\bold G_j \subseteq \bbP_j$ is generic
over $\bold V$ such that $q_n \in \bold G_j, \name p_n[\bold G_j \cap \bbP_{i_n}]
\in \bold G_j$
and if  $M = \name M_n[\bold G_j]$, then

\begin{enumerate}
\item[${{}}$]   $(\alpha) \quad \bold G_j \cap M$ is a subset of
$\bbP_j \cap M$ generic over $M$,

\item[${{}}$]  $(\beta) \quad {\cM_{\bbP_j}}[\bold G_j \cap M,M,y^*] \cap Y
\in D^{\gp}_{\beta_n} [M]$,

\item[${{}}$]  $(\gamma) \quad \name p_n[\bold G_j]$ belongs to $M$,
\end{enumerate}

\item[${{}}$]  $(h) \quad \langle Y_m \cap {\cM_{\bbP_{i_m}}}[\name{\bold G}_{i_m},
\name M_m,y^*], \name Y_m, \bbP_{i_m}, \name M_m: m \le n \rangle$ is
forced by $q_n$ to be an initial segment of a play of the game
$\Game_\alpha(N)$ in which the chooser uses the fixed winning
strategy $\sigma$.
\end{enumerate}

The proof is straight by the induction hypothesis on $\beta$, the fact that
$\name M_n, \name Y_n$ are $\bbP_{i_n}$-names
of objects from $\bold V$ and
$\bar{\bbQ} \restriction i_n$ is a ${\gp}-NNR^0_{\aleph_0}$-iteration.
Now let $q'= \bigcup\limits_{n<\omega}q_n$
and $\bold G' =\{p' \in \bbP_j \cap N: \bigvee\limits_{n<\omega} p' \leq   q_n        \}$. Then
$\bold G'$ and $q'$ are as required;
see also the proof of clause $(c)'$ of Theorem \ref{d9} for more details, where a more general
result is proved.

The theorem follows.
\end{proof}
\begin{remark}
\label{a30}
\begin{enumerate}
\item It is possible to use ``adding no reals+
clause (d)'' in the proof of clause (c) in order to weaken ``winner'' to ``not
loser''. Also we can use $\Game'_\alpha(N,N',\bbP)$, see Section \ref{second}.
The assumption ``$N_0 \in \bigcup\limits_{\beta < \alpha} {\cE}_\beta$'' can be
replaced by $``N_0 \in {\cE}'_0$ with ${\cE}'_0 \subseteq [{\cH}
(\chi^{\gp}_0)]^{\aleph_0}$ stationary''. We can also put extra
restrictions on $\bold G^*$ (and $\bold G^{**}$), for example we can require
${\cM}[\bold G^*,N_0,y^*]$ is large.

\item  The use of $\langle \chi_\alpha:\alpha < \ell g({\gp})
\rangle$ is not really necessary
as  all the properties depend just on $N_\ell \cap
{\cP}(\bbP_{\ell g(\bar{\bbQ})})$.

\item The proof of clause (c) being preserved can be applied to any $\bar{\bbQ}$
satisfying (a) + (c) of Definition \ref{a24} (so possibly adding reals),
but then we have to replace $D^{\gp}_\alpha(N)$ by a definition of such pseudo
filters with the winning strategy being absolute enough, e.g. for standard
$\bar D$.
\end{enumerate}
\end{remark}

\begin{definition}
\label{a33}
Let $\bar{\bbQ}$ be a countable support iteration of forcing notions.  It will be called ${\gp}$-proper if
it satisfies items (a) and (c) of Definition \ref{a24}.
\end{definition}

We  like to consider the parallel of having 2-completeness
systems and also to demand only non-losing rather than winning
in the assumption of Theorem \ref{a27}.

\begin{definition}
\label{a39}
Let $\kappa \in [2,\omega]$. We say that  $\bar{\bbQ}$
is a ${\gp}-NNR^0_\kappa$-iteration if items (a)-(c) of Definition \ref{a24} are satisfied
and clause (d) is replaced by $\kappa$-anti w.d., which is the same as
in (d) there, but with
$n < 1 + \kappa$ and $N_0 \in {\cE}^{\gp}_0$.
\end{definition}

The next theorem is a natural generalization of Theorem \ref{a27}.
\begin{theorem}
\label{a42}
Assume $\gp$ is a  reasonable parameter of length
$\omega_1$ which is a non-$\Game'$-loser, $2 \le \kappa < \aleph_0,
\bar{\bbQ}$ is a countable support iteration with ${\cP}(\Lim(\bar{\bbQ})) \subseteq
{\cH}(\chi^{\gp}_0),$  $\delta = \ell g(\bar{\bbQ})$
is a limit ordinal and for $i < \delta, \bar{\bbQ} \restriction i$ is a
${\gp}-NNR^0_\kappa$-iteration.  Then  $\bar{\bbQ}$ is a
${\gp}-NNR^0_\kappa$-iteration.
\end{theorem}
\begin{proof}
The proof is similar to the proof of Theorem \ref{a27}, with some changes, as
in \cite[Ch.VIII, Claim 4.10]{Sh:b}  and
\cite[Ch.XVIII, Proof 2.10C]{Sh:f}, so we do not give the details. The only main change is that
during the proof of clause (d), we add the following extra conditions to items (a)-(m):
\begin{enumerate}
\item[$(n)$]   $n_m$ is a power of 2, say $2^{n^*_m}$ and so we can rename
$\{\bold G^\ell_m:\ell < n_m\}$ as $\{\bold G^\eta_m:\eta \in {}^{n^*_m}2\}$,

\item[$(o)$]  The following conditions are satisfied:

\item[${{}}$] $(\alpha) \quad$ for
$\eta \in {}^{(n^*_m \ge)}2$, $ M_\eta \in M_1 \cap {\cE}^{\gp}_{j_\eta}$,  where $j_0 = \otp([i,j) \cap N_0)$ and if
 $\eta = \nu ^{\frown} \langle i \rangle,$ then $j_\eta = \otp([i,j) \cap M_\nu)$,

\item[${{}}$]   $(\beta) \quad M_{\langle\rangle} = N_0$,

\item[${{}}$] $(\gamma) \quad M_\eta \in M_{\eta^{\frown} \langle 0
\rangle} \cap M_{\eta^{\frown} \langle 1 \rangle}$,

\item[${{}}$]  $(\delta) \quad \eta \triangleleft \nu_1 \in {}^{n^*_m}2 \wedge
\eta \triangleleft \nu_2 \in {}^{n^*_m}2 \Rightarrow \bold G^{\nu_1}_m
\cap M_\eta = \bold G^{\nu_2}_m \cap M_\eta$ so we
 call it $K^\eta_m$

\item[${{}}$] $(\varepsilon) \quad$ for  $\eta \in {}^{n^*_m >}2$, $M_{\eta^{\frown} \langle 0 \rangle} =
M_{\eta^{\frown} \langle 1 \rangle}$,  call it
$N_\eta$,

\item[${{}}$]  $(\zeta) \quad$ for
$\eta \in {}^{(m^*_m >)}2$ and  $\ell < 2$, $N_\eta \in \cE^{\gp}_{j_{\eta^{\frown} \langle \ell \rangle}}$,

\item[${{}}$] $(\eta) \quad Y^\eta_m = {\cM}[K^{\eta^{\frown} \langle 0
\rangle}_m,N_\eta] \cap {\cM}[K^{\eta^{\frown}\langle 1 \rangle}_m,
N_\eta] \in D_{j_{\eta^{\frown} \langle 0 \rangle}}(N_\eta)$.
\end{enumerate}
The rest of the argument is essentially the same as before.
\end{proof}
The following is an immediate consequence of Theorems \ref{a27} \and \ref{a42}.
\begin{conclusion}
\label{a45}
Suppose ${\gp}$ is
a non-$\Game'$-loser reasonable parameter with $\ell g({\gp}) = \omega_1$,
 $\bar{\bbQ}$ is a countable support iteration and  $2 \le n(*) \le
\aleph_0$.
Then, $\bar{\bbQ}$ is a ${\gp}-NNR^0_{n(*)}$-iteration iff for each $i <
\ell g(\bar{\bbQ})$

\begin{enumerate}
\item[$(*)_i$]   $\name{\bbQ}_i$ is a proper forcing and
$\bbP_i,\name{\bbQ}_i$ satisfy clauses (d) + (c) of the
Definition ``${\gp}-NNR^0_{n(*)}$-iteration" with $i,i+1$ here standing
for $i,j$ there.
\end{enumerate}
\end{conclusion}

\begin{proof}
By induction on $j = \ell g(\bar{\bbQ})$.
For $j=0$ there is nothing to prove and for $j$ a successor ordinal, this follows easily from the definitions. For $j$ a limit ordinal,
 the result follows from Theorem \ref{a27} (for $n(*)=\aleph_0$) or Theorem \ref{a42} (for $2 \leq n(*) < \aleph_0$).
\end{proof}
We point out here that Clause (c) of Definitions \ref{a24} and \ref{a39}
really follows from earlier properties which play parallel roles.
\begin{lemma}
\label{a48}
\begin{enumerate}
\item Assume that ${\gp}$ is a standard reasonable parameter, $\alpha < \ell g({\gp}),N \in {\cE}^{\gp}_\alpha,
Y \in D^{\gp}_\alpha(N)$ and $\delta \le \omega_1 \cap N$ is a limit ordinal.
Then we can find sequences $\bar N = \langle N_i:i < \delta \rangle$
and $\bar\gamma=
\langle \gamma_i:i < \delta \rangle$ such that:

\begin{enumerate}
\item[$(a)$]   $N_i \in N$  is countable,
$N \cap \alpha \subseteq N_i$ and $N_i \in Y$ (for $i < \delta$),

\item[$(b)$]   $N_i \subseteq \bigcup\limits_{\beta \in \alpha \cap N}
({\cH}(\chi^{\gp}_\beta),\in)$, and $\beta \in \alpha \cap N_i
\Rightarrow N_i \restriction {\cH}(\chi^{\gp}_\beta) \prec
({\cH}(\chi^{\gp}_\beta),\in)$,

\item[$(c)$]   $i < j \Rightarrow N_i \subseteq N_j$,

\item[$(d)$]   if $i$ is a limit  ordinal, then $N_i = \bigcup\limits_{j<i} N_j$ and
$N \cap \bigcup\limits\{{\cH}(\chi^{\gp}_\beta):\beta \in \alpha \cap N\}
= \bigcup\limits_{j < \delta} N_j$, so we can stipulate $N_\delta = N$,

\item[$(e)$]  $\beta \in \alpha \cap N \Rightarrow \langle {\cH}
(\chi^{\gp}_\beta) \cap N_j:j \le i \rangle \in N_{i+1}$,

\item[$(f)$]   $\gamma_i \in N_i \cap \alpha, N_i \cap {\cH}
(\chi^{\gp}_{\gamma_i}) \in {\cE}^{\gp}_{\gamma_i} \cap Y$ and $\bar
\gamma \restriction (i+1) \in N_{i+1}$,

\item[$(g)$]   if $i \le \delta$ is a limit ordinal and either $(i = \delta$ and
 $\beta \in \alpha \cap N_i)$ or  $(i < \delta$ and $ \beta
\in \gamma_i \cap N_i)$, then for some
$j < i, \gamma_j = \beta$ and $y \in N_j$.

\item[$(h)$]   if $i \le \delta$ is a limit ordinal, then
$\{N_j \cap {\cH}(\chi^{\gp}_{\gamma_j}):j < i \text{ and }
\gamma_j < \gamma_i\} \in D^{\gp}_{\gamma_i}(N) = D^{\gp}_{\gamma_i}(N \cap {\cH}
(\chi^{\gp}_{\gamma_i})).$

\item[$(i)$]   if $\delta < N \cap \omega_1$ then $\delta \in N_0$, if
$\delta = N \cap \omega_1$ then $i < \delta \Rightarrow i \in N_0$.
\end{enumerate}

\item If ${\gp}$ is a standard reasonable parameter, $\bar{\bbQ}$ is a
countable support iteration, $\ell g(\bar{\bbQ}) = \beta +1,\bar{\bbQ}
\restriction \beta$ is ${\gp}-NNR^0_{k(*)}$-iteration and
$\Vdash_{\bbP_\beta} ``\name{\bbQ}_\beta$ is proper and
$(<^+ \omega_1)$-proper", then $\bar{\bbQ}$ is $\gp$-proper (see Definition \ref{a33}).

\item If $\ell g({\gp}) = \omega_1$, then  in part $(2)$, it suffices  to assume $\Vdash_{\bbP_i} ``\name{\bbQ}_i$ is
$(< \omega_1)$-proper''.
\end{enumerate}
\end{lemma}

\begin{proof}
(1). By induction on $i<\delta$ we prove that there are sequences
$\langle N_j:j < i \rangle \in N$ and $\langle \gamma_j: j<i \rangle$  satisfying the
relevant requirements, so that for some sequence
$\langle N'_j:j < i \rangle \in N$ with $N'_j \prec ({\cH}(\chi^{\gp}_\alpha), \in ),$  $N_j= N'_j \cap \bigcup\limits_{\beta \in \alpha \cap N}
({\cH}(\chi^{\gp}_\beta),\in)$.

For $i=0$, let $N'_0 \prec ({\cH}(\chi^{\gp}_\alpha), \in )$ be a countable model such that $N'_0 \in N \cap Y$ and such that $N'_0$ contains all relevant information, in particular, $N \cap \alpha \subseteq N'_0, \delta \subseteq N'_0$. Let $N_0= N'_0 \cap \bigcup\limits_{\beta \in \alpha \cap N}
({\cH}(\chi^{\gp}_\beta),\in).$
Pick also $\gamma_0$ such that  $\gamma_0 \in N_0 \cap \alpha$ and $N_0 \cap {\cH}
(\chi^{\gp}_{\gamma_0}) \in {\cE}^{\gp}_{\gamma_0} \cap Y$. Such  $\gamma_0$ exists as $\gp$ is standard.

If $i=j+1$ is a successor ordinal, let $N'_i  \prec ({\cH}(\chi^{\gp}_\alpha), \in )$ be a countable model such that:
\begin{itemize}
\item $N'_i \in N \cap Y$,
\item $N'_j \subseteq N_i,$
\item $\beta \in \alpha \cap N \Rightarrow \langle {\cH}
(\chi^{\gp}_\beta) \cap N_k:k \le j \rangle \in N_{i}$,
\item $\langle \gamma_k: k \leq j \rangle \in N_{i}$.
\end{itemize}
Then, using $\gp$ is standard, take  $\gamma_i \in N_i \cap \alpha$ such that $N_i \cap {\cH}
(\chi^{\gp}_{\gamma_i}) \in {\cE}^{\gp}_{\gamma_i} \cap Y$.

For limit $i$ set $N'_i=\bigcup\limits_{j<i}N_j$ and  $N_i= N'_i \cap \bigcup\limits_{\beta \in \alpha \cap N}
({\cH}(\chi^{\gp}_\beta),\in)$. Let also $\gamma_i$ be such that $\{N_j \cap {\cH}(\chi^{\gp}_{\gamma_j}):j < i \text{ and }
\gamma_j < \gamma_i\} \in D^{\gp}_{\gamma_i}(N) = D^{\gp}_{\gamma_i}(N \cap {\cH}
(\chi^{\gp}_{\gamma_i})).$

As $N$ is countable, we can choose $N'_j$'s such that clause (d) holds as well. This completes our inductive construction.

(2) We  show that $\bar{\bbQ}$ satisfies  clause (c) of Definition \ref{a24}. Thus let $i, j, \alpha, N, q, p, \bold G$
and $Y$ be as there. Without loss of generality, $i=\beta$ and $j=\beta+1$.

By the definition of $(<^+ \omega_1)$-proper, if
$\bold G_\beta \subseteq \bbP_\beta$ is generic over
$N,$ then $ {\cM}_{\bbP_\beta}[\bold{G}_\beta, N, y^*] \in D^{\gp}_\beta(N)$. Let $\delta = N \cap \omega_1$ and let
$\langle N_i:i < \delta \rangle$ be as in $(1)$.
Without loss of generality
$p(\beta) \in N_0 \cap \name{\bbQ}_\beta[\bold G_\beta]$.
Let $q' \ge p$ be $(N_i,\name{\bbQ}_\beta[\bold G_\beta])$-generic
for every $i < \delta$ \footnote{formally we only need to look at $\bar N' = \langle N'_i:i < \delta
\rangle,N'_i = N_i \cap {\cH}(\chi^{\gp}_0)$ and apply the $(<^+ \omega_1)$-properness to it.
} so that $q'\restriction \beta \geq q.$
Let also $\bold G(\beta) \subseteq \name{\bbQ}_\beta[\bold G_\beta]$  be generic over $N$ such that $q'\restriction \beta$ forces
$\bold G(\beta)$ has an upper bound in  $\name{\bbQ}_\beta[\bold G_\beta]$. Set $\bold{G}'=\bold G_\beta \ast \bold G(\beta).$
Then $q', \bold G'$ are as required.

(3) Follows from (2), as $\alpha \in N \cap \omega_1 \Rightarrow \delta = \omega \alpha \in N \cap
\omega_1$.
\end{proof}

\begin{remark}
\label{a51}
\begin{enumerate}
\item The results of this section include as special cases \cite[Ch.V, \S5, \S7]{Sh:b}.
There is no direct comparison with
\cite[Ch.VIII, \S4]{Sh:b}, \cite[Ch.VIII, \S4]{Sh:f}, but we can
make the notion  somewhat more complicated, to include the theorems
there in our context, but  this is not really needed for the examples discussed there (see Section
\ref{examples}).
The condition in \cite[Ch.VIII, \S4]{Sh:b} and  \cite[Ch.VIII, \S4]{Sh:f} involves having many sequences $\langle N_i:
i \le \delta \rangle$ such that if $p_0,p_1 \in \bbP_\alpha, p_\ell$ is
$(N_i,\bbP)$-generic for $i, p_\ell
\Vdash ``\name{\bold G}_{\bbP_\alpha} \cap N_0 = \bold G^*$'', then there
is $\bold G' \subseteq \Gen(N_0[\bold G^*], \name{\bbQ}_\alpha[\bold G^*]),
\bold G' \in N_0,$ such that $ \nVdash_{\bbP_{\alpha}} ``\name{\bold G'}$ has no bound in
$\name{\bbQ}_\alpha$''.  This speaks on a family of sequences
from $[{\cH}(\chi)]^{\aleph_0}$ rather than members of ${\cH}(\chi)$.

\item For \cite[Ch.XVIII, \S2]{Sh:f}, the comparison is not so easy.
Our problem is to ``carry" good $(N,\bbP_i,\langle \bold G_\ell:
\ell < n \rangle),\bold G_\ell \in \Gen(N,\bbP_i)$ with a bound, such
that we can ``increase $i$" and we can find $N',y \in N' \in N,
N' \prec N$ such that $(N',P_i,\langle \bold G_i \cap
N:\ell < n \rangle)$ is good enough.  In \cite[Ch.XVIII]{Sh:f} we are
carrying genericity in some $\bbP_{\bar \alpha},$ where $\bar \alpha \in
\text{ trind}(i)$\footnote{see \cite[Ch.XVIII, Definiton 2.1]{Sh:f}.}, but here we have much less.  But what we need is the
implication ``if $(N,\bbP_i,\bar G)$ is good we can extend it to good
$(N,\bbP_{i+1},\bar G')$", so making good weaker
generates an incomparable notion and clearly there are other variants.

\item The iteration theorems proved in this section can be used to give alternative proofs of the consistency results in
\cite[Ch.XVIII, \S1]{Sh:f}
(see Section \ref{examples}).
\end{enumerate}
\end{remark}

\section {Delayed properness} \label{delayed}

In this section we introduce several notions that will be used in sections \ref{examples}
and \ref{second}.
  We concentrate  on simple reasonable parameters and we present two versions for it.
The simpler one is version 2 for
which simplicity is a very natural demand.
The proof of the next lemma is  straightforward in which  a general way to create simple reasonable parameters is introduced.
\begin{lemma}
\label{b3}
\begin{enumerate}
\item Assume

\begin{enumerate}
\item[$(a)$]   $\bar \chi = \langle \chi_\alpha:\alpha < \alpha^* \rangle$
increases fast enough, so that $
{\cH}((\bigcup\limits_{\beta < \alpha}\chi_\beta)^+) \in {\cH}(\chi_\alpha)$,

\item[$(b)$]   ${\cE}_\alpha \subseteq \{N: N$ is a countable elementary
submodel of $({\cH}(\chi_\alpha),\in)\}$ is stationary,

\item[$(c)$]   $R_\alpha \in {\cH}(\chi_\alpha)$ and $N \in
{\cE}_\alpha$ implies $\langle \chi_\beta:\beta < \alpha \rangle \in N,
\langle R_\beta:\beta \le \alpha \rangle \in N$ and $\langle
{\cE}_\beta:\beta < \alpha \rangle \in N$.
\end{enumerate}
Then there is a standard reasonable parameter ${\gp}$ with
$\ell g({\gp}) = \alpha^*, \chi_\alpha^{\gp}=\chi_\alpha,~{\cE}^{\gp}_\alpha = {\cE}_\alpha$
and $R^{\gp}_\alpha = R_\alpha$.

\item If in addition clause (d) below holds, then  ${\gp}$ is a simple
standard reasonable parameter (recall Definition \ref{a9}(4)), where
\begin{enumerate}
\item[$(d)$]   $\beta \in N \in {\cE}_\alpha,\beta < \alpha
\Rightarrow N \cap {\cH}(\chi_\beta) \in {\cE}_\beta$.
\end{enumerate}

\item If $\chi_\alpha = (\beth_{2 \alpha +1})^+$ for $\alpha < \alpha^*,
R_\alpha \in {\cH}(\chi_\alpha)$, then $\chi_\alpha$ increases fast
enough.  If $\langle \chi_\alpha:\alpha < \alpha^* \rangle,\langle R_\alpha:
\alpha < \alpha^* \rangle$ are as in part (1), $\chi \le \chi_0,{\cE}
\subseteq [{\cH}(\chi)]^{\le \aleph_0}$ stationary and we let ${\cE}
_\alpha = \{N: N$ is a countable elementary submodel of $({\cH}(\chi_\alpha),
\in)$ and $\langle \chi_\beta:\beta < \alpha \rangle,\langle R_\beta:\beta
\le \alpha \rangle,{\cE}$ belong to $N$ and $N \cap {\cH}(\chi) \in
{\cE}\}$, then  the assumptions of parts (1) and
(2) above hold.
\end{enumerate}
\end{lemma}
\begin{proof}
(1). Let ${\gp}=\langle \bar\chi, \bar{R}, \bar{\cE}, \bar D \rangle$, where $\bar\chi, \bar{R}$ and $\bar{\cE}$
are given as above and $\bar D$ is defined as in  Definition \ref{a9}(1). Then $\gp$ is easily seen to be a  standard reasonable parameter
as required.

Items (2) and (3) are clear.
\end{proof}
We now define an extension of the games $\Game_\alpha(N, {\gp})$ and $\Game_\alpha(N, N', {\gp})$ given in Definition \ref{a15}.
\begin{definition}
\label{b6}
Let ${\gp}$ be a reasonable parameter and $\alpha \le \beta
< \ell g({\gp})$.
\begin{enumerate}
\item For $N \in {\cE}^{\gp}_\beta$ such that $\alpha \in N$, we define a game
$\Game_{\alpha,\beta}(N) = \Game_{\alpha,\beta}(N,{\gp})$ of length $\omega $ as follows.
 In the $n$-th move:
\begin{enumerate}
\item[$(a)$]  the \underline{challenger} chooses $X_n \in D^{\gp}_\beta(N)$ such that
$m < n \Rightarrow X_n \subseteq X_m$,

\item[$(b)$]   the \underline{chooser} chooses $\alpha_n \in \alpha \cap N$,

\item[$(c)$]   the \underline{challenger} chooses $\beta'_n \in \beta \cap
N$ and $y'_n \in N \cap {\cH}(\chi^{\gp}_{\alpha_n})$,

\item[$(d)$]  the \underline{chooser} chooses $\beta_n \in \beta \cap N
\backslash \beta'_n$ together with $M_n \in X_n \cap {\cE}^{\gp}_{\beta_n}$,
 $y_n \in M_n \cap {\cH}(\chi^{\gp}_{\alpha_n})$ and $Y_n \in D^{\gp}_{\beta_n}
(M_n)$ satisfying:
\begin{itemize}
\item $\alpha_n \le \beta_n,$
\item $y'_n \in M_n,$
\item $\alpha_n \in M_n$,

\item  $Y_n \subseteq X_n$ and $Y_n \in N$
\end{itemize}
\item[$(e)$]  the \underline{challenger} chooses $M'_n \in Y_n \cap
{\cE}^{\gp}_{\alpha_n} \cap \left(M_n \cup \{M_n \cap {\cH}
(\chi^{\gp}_{\alpha_n})\}\right)$ satisfying $y_n,y'_n \in M'_n$ and
chooses $Z_n \in D^{\gp}_{\alpha_n}(M'_n) = D^{\gp}_{\alpha_n}(M'_n \cap {\cH}(\chi
^{\gp}_{\alpha_n}))$ such that $Z_n \subseteq Y_n$.
\end{enumerate}
At the end, the chooser wins the play if
\[
\bigcup\{Z_n \cup \{M'_n\}:n < \omega\} \in D_\alpha^+(N)
= D_\alpha^+(N \cap {\cH}(\chi^{\gp}_\alpha)),
\]
where $D_\alpha=D^{\gp}_\alpha$.

\item We call $\Game_{\alpha,\beta}(N) = \Game_{\alpha,\beta}(N,{\gp})$,
defined in clause (1), version 1 of the game. Version 2 of the game is defined similarly, where
\begin{itemize}
\item  in clause (d) we
require $\alpha_n \le_{\gp} \beta_n$.

\item in clause (e), we add the requirement
$M'_n = M_n \cap {\cH}(\chi^{\gp}_{\alpha_n})$,
\end{itemize}
 If we do not mention the version,
it means that it holds for both versions.

\item Assume $N \in N' \prec ({\cH}(\chi),\in)$.
We define the game $\Game'_{\alpha,\beta}(N,N',{\gp})$ similarly, where
items (a) - (e) are replaced by:
\begin{enumerate}
\item[$(a)'$]   the \underline{challenger} chooses $X_n \in D^{\gp}_\beta(N) \cap N'$ such
that $m < n \Rightarrow X_n \subseteq X'_m$,

\item[$(b)'$]  the \underline{chooser} chooses $\alpha_n \in \alpha \cap
N$ and $X'_n \subseteq X_n$ such that $X'_n \in D^{\gp}_\beta(N)\cap N'$,

\item[$(c)'$]  like (c) above,

\item[$(d)'$]   like (d) above, but we replace ``$Y_n \in N$'' by ``$Y_n \in N'$'',

\item[$(e)'$]  like (e) above, but add $Z_n \in N'$.
\end{enumerate}
Note that every proper initial segment of a play belongs to $N'$.
\end{enumerate}
\end{definition}

\begin{definition}
\label{b9}
Let ${\gp}$ be a reasonable parameter.
\begin{enumerate}
\item For $\alpha \le \beta < \ell g({\gp})$, we say ${\gp}$ is a
$\Game_{\alpha,\beta}$-winner (resp. non-$\Game_{\alpha,\beta}$-loser), if
 for some $x \in {\cH}(\chi^{\gp}_\beta)$ we have:

\begin{enumerate}
\item[${{}}$]   if $\{x,\alpha\} \in N \in {\cE}^{\gp}_\beta$, then  the
chooser wins the game $\Game_{\alpha,\beta}(N,{\gp})$

\item[${{}}$]  (resp. the challenger does not win in the game $\Game_{\alpha,\beta}
(N,{\gp})$).
\end{enumerate}
\item  Similarly we can define when $\gp$ is a  $\Game'_{\alpha,\beta}$-winner (resp. non-$\Game'_{\alpha,\beta}$-loser).

\item  For any function $f:\ell g({\gp}) \rightarrow {\cP}(\ell g
({\gp}))$ we can replace $\alpha,\beta$ by $f$, so that $\gp$ is a $\Game_f$-winner (resp. non-$\Game_f$-loser)  if for every
$\alpha < \ell g({\gp})$ and $\beta \in f(\alpha)$,
${\gp}$ is a $\Game_{\alpha,\beta}$-winner (resp. non-$\Game_{\alpha,\beta}$-loser).
\end{enumerate}
\end{definition}
Given  a reasonable parameter $\gp$, we define some families of functions as follows.

\begin{definition}
\label{b999}
Let ${\gp}$ be a reasonable parameter.
\begin{enumerate}
\item Let ${\cF}^{\gp}$ be the family of functions $f$ from
$\ell g({\gp})$ to ${\cP}(\ell g({\gp}))$ such that for each
$\alpha < \ell g({\gp}),f(\alpha)$ is a nonempty subset of
$[\alpha,\ell g({\gp}))$.

\item Let ${\cF}^{\gp}_{\club}$ be
the set of $f \in {\cF}^{\gp}$ such that for each $\alpha < \ell g
({\gp}),f(\alpha)$ is a club of $\ell g({\gp})$.

\item Let
$\cF^{\gp}_{\nd}$ be the set of $f \in \cF^{\gp}$ such that for each
$\alpha < \ell g({\gp}),f(\alpha)$ is an end segment
of $\ell g({\gp})$, we then may identify $f(\alpha)$ with
$\min(f(\alpha))$.

\item Call $f \in {\cF}^{\gp}$ decreasing continuous if
\begin{enumerate}
\item[$(a)$]  $\alpha < \beta < \ell g(\gp) \Rightarrow f(\alpha) \supseteq
  f(\beta)$,

\item[$(b)$]  for limit $\delta < \ell g(\gp)$ we have $f(\delta) =
  \bigcap\{f(\alpha):\alpha < \delta\}$.
\end{enumerate}
  Let also $f \le g$ mean that $(\forall \alpha < \ell g({\gp}))
(g(\alpha) \subseteq f(\alpha))$.

\item Let ${\cF}^{\gp}_{\dc}$ be the set of decreasing continuous functions $f \in
   \cF^{\gp}_{\club}$.
\end{enumerate}
\end{definition}
We are  interested in dealing with winning strategies for the games $\Game_f$, where $f$ is a function coming from one of the
  above family of functions. The next lemma gives some obvious monotonicity properties for these classes of functions.
\begin{lemma}
\label{b12}
Assume ${\gp}$ is a reasonable parameter.
\begin{enumerate}
\item If $\alpha \le_{\gp} \alpha' \le \beta = \beta'
< \ell g({\gp})$, and ${\gp}$ is a $\Game_{\alpha,\beta}$-winner,
then  it is $\Game_{\alpha',\beta'}$-winner.  Similarly for
$\Game'$-winner, non-$\Game$-loser and non-$\Game'$-loser.

\item  If ${\gp}$ is a $\Game_{\alpha,\beta}$-winner, then
${\gp}$ is a $\Game'_{\alpha,\beta}$-winner and a
non-$\Game_{\alpha,\beta}$-loser.
If ${\gp}$ is a $\Game'_{\alpha,\beta}$-winner or
non-$\Game_{\alpha,\beta}$-loser, then  ${\gp}$ is
non-$\Game'_{\alpha,\beta}$-loser.

\item Assume $f,g \in {\cF}^{\gp}$ and $f \le g$.
If ${\gp}$ is a $\Game_f$-winner (or $\Game'_f$-winner)
(or non-$\Game_f$-loser) (or non-$\Game'_f$-loser), then
${\gp}$ is a $\Game_g$-winner (or $\Game'_g$-winner)
(or non-$\Game_g$-loser) (or non-$\Game'_g$-loser).
\end{enumerate}
\end{lemma}

\begin{proof}
(1) Suppose $\sigma$ is a winning strategy for chooser in the game $\Game_{\alpha,\beta}$. We show that it is a winning strategy for chooser in the game $\Game_{\alpha',\beta'}$ as well. Suppose not, so, following the notation of Definition \ref{b6},
\[
\bigcup\{Z_n \cup \{M'_n\}:n < \omega\} \notin D^+_{\alpha'}(N) \,
(= D^+_{\alpha'}(N \cap {\cH}(\chi^{\gp}_{\alpha'}))).
\]
It then follows that
\[
N \cap {\cH}(\chi^{\gp}_{\alpha'}) \setminus \bigcup\{Z_n \cup \{M'_n\}:n < \omega\} \in D_{\alpha'}(N).
\]
But then, as $\alpha \leq_{\gp} \alpha',$
\[
N \cap {\cH}(\chi^{\gp}_{\alpha}) \setminus \bigcup\{Z_n \cup \{M'_n\}:n < \omega\} \in D_\alpha(N).
\]
But, as $\sigma$ is a winning strategy for chooser in the game $\Game_{\alpha,\beta}$, we have
\[
\bigcup\{Z_n \cup \{M'_n\}:n < \omega\} \in D^+_{\alpha}(N),
\]
a contradiction.

The proof of clause (2) is similar to the proof of Lemma \ref{a18}(2)
 and the proof of clause (3) is straightforward.
\end{proof}

\begin{lemma}
\label{b15}
\begin{enumerate}
\item Assume ${\gp}$ is a standard reasonable parameter.  Then
${\gp}$ is a winner.

\item  If ${\gp}$ is a reasonable parameter and it is a winner, then
${\gp}$ is a $\Game_{\alpha,\alpha}$-winner.

\item  If ${\gp}$ is a reasonable parameter and it is a winner and
$\alpha \le \beta < \ell g({\gp})$, then  ${\gp}$
is a $\Game_{\alpha,\beta}$-winner (hence $\Game_f$-winner for any $f:\ell g({\gp}) \rightarrow {\cP}(\ell g
({\gp}))$).

\item Similarly with $\Game'$-winner, $\Game'_{\alpha,\beta}$ winner and/or
with the ``non-loser" cases.
\end{enumerate}
\end{lemma}

\begin{proof}
Clause (1) follows from Lemma \ref{a18}. The proof of items (2)-(4) is easy and follows from the Definitions \ref{a15}
and \ref{b6}.
\end{proof}
We now define an interpretation of reasonable parameters in the forcing extensions and show that these interpretations are
reasonable parameters in the corresponding extension.
\begin{definition}
\label{b18}
Let ${\gp}$ be a reasonable parameter and let $\bbP$ be a proper
forcing notions which adds no new reals.   Suppose that ${\cP}(\bbP)
\in {\cH}(\chi^{\gp}_0)$ and $\bold G_{\bbP} \subseteq \bbP$ is generic over
$\bold V$.
We interpret ${\gp}$ in $\bold V^{\bbP}$ as ${\gp}' =
{\gp}^{\bold V[\bold G_{\bbP}]}$,  defined
as follows:
\begin{enumerate}
\item[$(a)$]   $\chi^{{\gp}'}_\alpha = \chi^{\gp}_\alpha$,

\item[$(b)$]   $R^{{\gp}'}_\alpha = \langle R^{\gp}_\alpha, \bbP,\bold G_P \rangle$,

\item[$(c)$]   ${\cE}^{{\gp}'}_\alpha = \{N[\bold G_P]:N \in {\cE}^{\gp}
_\alpha$, $\bbP \in N$ and $N[\bold G_{\bbP}] \cap V = N\}$,

\item[$(d)$]  $D^{{\gp}'}_\alpha(N[\bold G_{\bbP}]) = \{\{M[\bold G_{\bbP}] \in
{\cE}^{{\gp}'}_\alpha:M \in Y \cap \bigcup\limits_{\beta < \alpha}
{\cE}^{\gp}_\beta\}:Y \in D^{\gp}_\alpha(N)\}$.
\end{enumerate}
We also use ${\gp}^{\bbP}$ for ${\gp}' =
{\gp}^{\bold V[\bold G_{\bbP}]}$.
\end{definition}
The proof of the next lemma is straightforward.
\begin{lemma}
\label{b21}
Let ${\gp},\bbP$ and $\bold G_{\bbP}$ be as in Definition \ref{b18}.
\begin{enumerate}
\item  ${\gp}^{\bold V[\bold G_{\bbP}]}$ is
a reasonable parameter in $\bold V[\bold G_{\bbP}]$.

\item If ${\gp}$ is, in $\bold V$, a $\Game$-winner (or
non-$\Game$-loser or $\Game'$-winner or non-$\Game'$-loser), then
${\gp}^{\bold V[\bold G_{\bbP}]}$ is so in $\bold V[\bold G_{\bbP}]$.

\item If ${\gp}$ is, in $\bold V$, a $\Game_{\alpha,\beta}$-winner (or non-
$\Game_{\alpha,\beta}$-loser or $\Game'_{\alpha,\beta}$-winner or
non-$\Game'_{\alpha,\beta}$-loser), then
${\gp}^{\bold V[\bold G_{\bbP}]}$ is so in $\bold V[\bold G_p]$.
\end{enumerate}
\end{lemma}
\begin{proof}
(1). It is easily seen, using the fact that $\bbP$ is an NNR proper forcing notion, that ${\gp}^{\bold V[\bold G_{\bbP}]}$ satisfies items (a)-(i)  of Definition \ref{a3} in $\bold V[\bold G_{\bbP}]$.

(2). Suppose ${\gp}$ is a winner  in $\bold V$. We show that ${\gp}'={\gp}^{\bold V[\bold G_{\bbP}]}$ is a winner in $\bold V[\bold G_{\bbP}]$. Suppose $\alpha < \ell g(\gp)$ and $N' \in {\cE}^{{\gp}'}_\alpha.$ Then for some $N \in {\cE}^{{\gp}}_\alpha, N'=N[\bold G_{\bbP}]$. Let $\sigma$ be a winning strategy for $\gp$ for the game $\Game_\alpha(N, \gp)$. We define the winning strategy $\sigma'$ for the game $\Game_{\alpha}(N', {\gp}')$ in $\bold V[\bold G_{\bbP}]$ as follows. Following the notation of Definition \ref{a15},
in the $n$-th move, challenger chooses some $X'_n \in D^{{\gp}'}_\alpha(N')$. So for some
 $X_n \in D^{\gp}_\alpha(N)$, we have $X'_n= \{M[\bold G_{\bbP}] \in
{\cE}^{{\gp}'}_\alpha:M \in X_n \cap \bigcup\limits_{\beta < \alpha}
{\cE}^{\gp}_\beta\}$. Let $M_n \in X_n$ and $Y_n \subseteq M_n \cap X_n$ with $Y_n \in D^{\gp}(M_n) \cap N$
be the play that chooser does using the strategy $\sigma$. Let $M'_n=M_n[\bold G_{\bbP}]$
and $Y'_n=\{M[\bold G_{\bbP}]: M \in Y_n              \}$ be what chooser replies via $\sigma'.$
Then challenger chooses some $Z'_n \subseteq Y'_n$ with $Z'_n \in D^{{\gp}'}(M_n')$. Thus  for some $Z_n \subseteq Y_n$ with $Z_n \in D^{{\gp}}(M_n) \cap  \bold V,$
we have $Z'_n = \{M[\bold G_{\bbP}]: M \in Z_n         \}.$
Since $\sigma$ is a winning strategy,
$\bigcup \{\{M_n\} \cup Z_n: n<\omega   \} \in D^{\gp}(N).$
 It then follows that
 $$\bigcup \{\{M'_n\} \cup Z'_n: n<\omega   \} \in D^{{\gp}'}(N').$$
Thus $\sigma'$ is a winning strategy for chooser for the game $\Game_{\alpha}(N', {\gp}')$.
The proof for
non-$\Game$-loser or $\Game'$-winner or non-$\Game'$-loser is the same.

 Clause (3) can be proved similarly.
\end{proof}


\begin{definition}
\label{b24}
Let ${\gp}$ be a reasonable parameter
and let $\bbQ$ be a forcing notion.
\begin{enumerate}
\item For $\alpha \le \beta < \ell g({\gp})$, we say $\bbQ$
is $({\gp},\alpha,\beta)$-proper if
${\cP}(\bbQ) \in {\cH}(\chi^{\gp}_0)$ and:

\begin{enumerate}
\item[$(*)$]   for some $x \in {\cH}(\chi^{\gp}_\beta)$, if  $N \in
  \cE^{\gp}_\beta, \{x,\bbQ,\alpha\} \in N, p \in N \cap \bbQ$ and
$Y \in D^{\gp}_\alpha(N)$, then for some $q$ we have:

\begin{enumerate}
\item[$(a)$]  $p \le q \in \bbQ$,

\item[$(b)$]  $q$ is $(N,\bbQ)$-generic,

\item[$(c)$] for some $N' \in ({\cE}^{\gp}_\alpha \cap N \cap
Y) \cup \{N \cap {\cH}(\chi^{\gp}_\alpha)\}$ satisfying $\alpha =
\beta \Rightarrow N' = N$ we have
$q \Vdash_{\bbQ} ``{\cM}_{\bbQ}[\name{\bold G}_{\bbQ},N',y^*] \cap Y \in
D^{\gp}_\alpha(N')"$ where $y^* = \langle x,p,\bbQ \rangle$. Note that this implies $q$
is $(N',\bbQ)$-generic.
\end{enumerate}
We call the above, version 1 of $({\gp},\alpha,\beta)$-properness.
Version 2 of $({\gp},\alpha,\beta)$-properness is defined similarly, but we demand $N' = N \cap {\cH}(\chi^{\gp}
_\alpha)$ and  $\alpha \le_{\gp} \beta$.
\end{enumerate}

\item We say $\bbQ$ is $({\gp},f)$-proper, where $f \in {\cF}^{\gp}$ (see Definition \ref{b999}), when
for every $\alpha < \ell g({\gp})$ and $\beta \in f(\alpha)$, the
forcing notion $\bbQ$ is $({\gp},\alpha,\beta)$-proper.

\item We say $\bbQ$ is ${\gp}$-proper if $\bbQ$ is $({\gp},
\alpha,\alpha)$-proper for $\alpha < \ell g({\gp})$.  We say $\bbQ$
is almost ${\gp}$-proper if $\bbQ$ is $({\gp},f)$-proper for some
$f \in {\cF}^{\gp}$.
\end{enumerate}
\end{definition}
The next lemma shows some relations between the above defined notions.
\begin{lemma}
\label{b27}
Assume ${\gp}$ is a simple reasonable parameter.
\begin{enumerate}
\item If $\alpha' \le \alpha \le \beta \le \beta' < \ell g({\gp})$
(for version 2 we demand $\alpha' \le_{\gp} \beta'$ and $\alpha
\le_{\gp} \beta$) and $\bbQ$ is a $({\gp},\alpha,\beta)$-proper
forcing notion, then  $\bbQ$ is a $({\gp},\alpha',\beta')$-proper
forcing notion.

\item Assume $f,f'$ are in ${\cF}^{\gp}$ and $f \le f'$.
If $\bbQ$ is a $({\gp},f)$-proper forcing notion, then  $\bbQ$ is a
$({\gp},f')$-proper forcing notion.
\end{enumerate}
\end{lemma}

\begin{proof}
(1) Suppose $\bbQ$ is  $({\gp},\alpha,\beta)$-proper as witnessed by $x \in \cH(\chi^{\gp}_{\beta}) \subseteq \cH(\chi^{\gp}_{\beta'})$. We show that
$x$ witnesses that $\bbQ$ is  $({\gp},\alpha',\beta')$-proper as well. Thus let  $N \in
  \cE^{\gp}_{\beta'}$ with $\{x,\bbQ,\alpha\} \in N, p \in N \cap \bbQ$ and
$Y \in D^{\gp}_{\alpha'}(N)$. Then $N \cap \cH(\chi^{\gp}_{\beta}) \in  \cE^{\gp}_{\beta}$. Let $q \geq p$ be as in Definition \ref{b24}
and it
witnesses that $\bbQ$ is  $({\gp},\alpha,\beta)$-proper with respect to $N \cap \cH(\chi^{\gp}_{\beta}), p$ and $Y$. Then $q$ witnesses $\bbQ$ is
$({\gp},\alpha',\beta')$-proper with respect to $N, p$ and $Y$.

(2) is clear, as for every $\alpha < \ell g({\gp})$, $f'(\alpha) \subseteq f(\alpha)$.
\end{proof}
It follows from the above lemma that if ${\gp}$ is a simple reasonable parameter and $f \in{\cF}^{\gp}$, then
\begin{center}
${\gp}$-proper $~\implies ({\gp}, f)$-proper $~\implies~$ almost ${\gp}$-proper.
\end{center}
We may like to consider $({\gp},f)$-properness for iterations which may
add reals.  Then we have to replace $D^{\gp}_\alpha(N)$ by a definition
which is absolute enough (and the non-loser versions have to be
absolute enough as well). In such situations, it
is natural to restrict ourselves to those sets $Y$ which are ${\gp}$-closed, see below.

\begin{definition}
\label{b33}
Let ${\gp}$ be a simple reasonable parameter, $\alpha < \ell g({\gp})$ and $N \in {\cE}^{\gp}_\alpha$.
A subset $Y \subseteq N$ is called ${\gp}$-closed if:

\begin{enumerate}
\item[$(a)$]   $Y \subseteq N \cap \bigcup\limits_{\beta < \alpha}
{\cE}^{\gp}_\beta$,

\item[$(b)$]   if $M \in N \cap {\cE}^{\gp}_\beta$,  $\beta < \alpha$
(hence $\beta \in \alpha \cap M \subseteq \alpha \cap N$),
$\gamma \in M \cap \beta$ and $M \cap {\cH}(\chi^{\gp}_\gamma) \in
{\cE}^{\gp}_\gamma$,\footnote{note that this requirement is redundant if ${\gp}$ is simple or
just $\gamma \le_{\gp} \beta$.} then
$$M \cap
{\cH}(\chi^{\gp}_\gamma) \in Y \iff M \in Y,$$

\item[$(c)$]   if $\beta < \alpha, M_\ell \in N \cap {\cE}^{\gp}_\beta$
(hence $\beta \in \alpha \cap N$), $M_\ell \subseteq M_{\ell +1}$ for
$\ell < \omega$ and $M = \bigcup\limits_{\ell < \omega} M_\ell \in
N \cap {\cE}^{\gp}_\beta$, and even $\langle M_\ell:\ell < \omega \rangle
\in N$, then  $$\left(\bigwedge\limits_{\ell<\omega} M_\ell \in Y\right) \Rightarrow M \in Y.$$
\end{enumerate}
\end{definition}

\begin{definition}
\label{b37}
Assume $\gp$ is a reasonable parameter, $\bbP$ is a forcing notion and
$\name{\bbQ}$ is a $\bbP$-name of a forcing notion.
\begin{enumerate}
\item For $\alpha \le \beta < \ell g(\gp)$, we say  $\name{\bbQ}$
has $(\kappa,\alpha,\beta)$-anti w.d. above $\bbP$ (or
$(\bbP,\name{\bbQ})$ has
$(\kappa,\alpha,\beta)$-anti-w.d.),
if
clause (A) implies clause (B), where:

\begin{enumerate}
\item[$(A)$]  $(a) \quad N_0 \in \cE^{\gp}_\alpha$ and $N_1 \in \cE^{\gp}_\beta$,

\item[${{}}$]  $(b) \quad N_0 \in N_1$ and $\{\bbP,\name{\bbQ}\} \in N_0$,

\item[${{}}$]  $(c) \quad n< 1 + \kappa$,

\item[${{}}$]  $(d) \quad p_\ell \in \bbP$ is $(N_\iota,\bbP)$-generic
  for $\ell < n$
  and $\iota = 0,1$,

\item[${{}}$]  $(e) \quad p_\ell \Vdash ``\name{\bold G}_{\bbP} \cap N_1 =
  \bold G^\ell"$ for $\ell < n$,

\item[${{}}$]  $(f) \quad \bold G^\ell \cap N_0 = \bold G^*$ for $\ell
  < n$,

\item[${{}}$]  $(g) \quad Y = \bigcap\limits_{\ell < n} \cM_{\bbP}[\bold
  G^\ell,N_1, y]$ belongs to $D^{\gp}_\beta(N_1)$, where $y= \langle N_0, \bbP, \name{\bbQ} \rangle,$

\item[${{}}$]  $(h) \quad \name q$ is a $\bbP$-name of a member of
  $\name{\bbQ}$,

\item[${{}}$]  $(i) \quad \name q \in N_0$,
\end{enumerate}

\begin{enumerate}
\item[$(B)$]  there is a triple $( \langle p'_\ell: \ell < n \rangle,\name q',\bold G^{**})$ such that:
\sn
\item[${{}}$]   $(a) \quad \name q'$ is a $\bbP$-name of a member of
  $\name{\bbQ}$,

\item[${{}}$]  $(b) \quad p '_\ell \Vdash ``\name q \le \name q'$"
  for $\ell < n$,

\item[${{}}$]  $(c) \quad \bold G^{**} \in \Gen(N_0, \bbP * \name{\bbQ}),$

\item[${{}}$]  $(d) \quad (p'_\ell,\name q') \Vdash ``\bold G_{\bbP *
\name{\bbQ}} \cap N_0 = \bold G^{**}"$ for $\ell < n$.
\end{enumerate}
\item Given a function  $f \in \cF^{\gp}$, we say  $\name{\bbQ}$
has $(\kappa, f)$-anti w.d. above $\bbP$ (or
$(\bbP,\name{\bbQ})$ has
$(\kappa, f)$-anti-w.d.), if  for every $\alpha < \ell g({\gp})$, $\name{\bbQ}$
has $(\kappa,\alpha,f(\alpha))$-anti w.d. above $\bbP$.
\end{enumerate}
\end{definition}
\begin{remark}
\label{remarkb37}
We may assume that the conditions $p_\ell, \ell < n,$ in clause (1)(A)(d) of the above definition are pairwise
incompatible. If $p_\ell$ and $p_\iota$ are compatible, then by clause (1)(A)(e), $C^\ell=C^\iota,$
so we may replace $p_\ell, p_\iota$ by a common extension $p_{\ell, \iota}$ of them and
take $C^{\ell,\iota}=C^\ell.$
\end{remark}
As the sets ${\cH}(\chi_\alpha)$ may change with forcing, we may prefer to
use ${\cE}_\alpha \subseteq [\chi_\alpha]^{\le \aleph_0}$. For this reason, we define the notion of
ordinal based parameter, and we will show that any ordinal based parameter gives naturally a reasonable parameter.

\begin{definition}
\label{b42}
\begin{enumerate}
\item We call ${\gp}$ an o.b. (ordinal based) parameter if
$${\gp} = (\bar \chi^{\gp},\bar R^{\gp},
\bar{\cE}^{\gp},\bar D^{\gp}) (= (\bar \chi,\bar R,
\bar{\cE},\bar D)),$$
 where for some ordinal $\alpha^*$, called
$\ell g({\gp})$, we have:

\begin{enumerate}
\item[$(a)$]   $\bar \chi = \langle \chi_\alpha:\alpha < \alpha^* \rangle,$
where $\chi_\alpha$ is a regular cardinal and ${\cH}((\bigcup\limits_{\beta <
  \alpha} \chi_\beta)^+) \in {\cH}(\chi_\alpha)$,

\item[$(b)$]  $\bar R = \langle R_\alpha:\alpha < \alpha^* \rangle,$ where
$R_\alpha$ is an $n(R_\alpha)$-place relation on some bounded subset of
$\chi_\alpha$, \footnote{we could have asked ``on $\chi_\alpha$'', as there is no real difference.}

\item[$(c)$]   $\bar{\cE} = \langle {\cE}_\alpha:\alpha < \alpha^*
\rangle$, where ${\cE}_\alpha \subseteq [\chi_\alpha]^{\le \aleph_0}$ is
stationary,

\item[$(d)$]   $\bar D = \langle D_\alpha:\alpha < \alpha^* \rangle$ and
$D_\alpha$ is a function with domain ${\cE}_\alpha$ and for each $a \in
{\cE}_\alpha,D_\alpha(a)$ is a pseudo-filter on $a$,

\item[$(e)$]   let ${\gp}^{[\alpha]} = \langle \bar \chi \restriction
\alpha,\bar R \restriction (\alpha +1),\bar{\cE} \restriction \alpha,
\bar D \restriction \alpha \rangle$,

\item[$(f)$]   if $a \in {\cE}_\alpha$ and $X \subseteq a$,
  then
  $$X \in D_\alpha(a) \iff X \cap
\bigcup\limits_{\beta < \alpha} {\cE}_\beta \in D_\alpha(a).$$
\end{enumerate}

\item An o.b. parameter ${\gp}$ is simple if in addition, it satisfies:

\begin{enumerate}
\item[$(g)$]   if $a \in {\cE}_\alpha$, $X \in D_\alpha(a)$ and
$\beta \in \alpha \cap a$, then  $a \cap \chi_\beta \in {\cE}_\beta$.
\end{enumerate}

\item For ${\gp}$ as above let ${\gq} ={\gp}^{\bold V}$ be defined by
\begin{itemize}
\item
$\ell g({\gq}) = \ell g({\gp}),$
\end{itemize}
and we define by induction on $\alpha < \ell g({\gp})$
\begin{itemize}

\item $\chi^{\gq}_\alpha = \chi^{\gp}_\alpha,$

\item $R^{\gq}_\alpha = R^{\gp}_\alpha,$

\item ${\cE}^{\gq}_\alpha = \{N \prec ({\cH}(\chi^{\gq}_\alpha),\in):
N \text{ is countable}, N \cap \chi^{\gp}_\alpha \in
{\cE}^{\gp}_\alpha$ and ${\gq}^{[\alpha]} \in N\}$.
Note that ${\gq}^{[\alpha]}$ is well defined by the induction hypothesis.

\item $D^{\gq}_\alpha(N)$ is defined as

\[\begin{array}{ll}
\qquad \qquad \quad D^{\gq}_\alpha(N) = \{Y':&Y' \subseteq \bigcup\limits_{\beta < \alpha}
{\cE}^{\gq}_\beta \text{ and for some } y \in N \cap \bigcup\limits_{\beta <
\alpha} {\cH}(\chi^{\gp}_\beta)  \\
  & \text{ and  } Y \in D^{\gp}_\alpha(N \cap \chi^{\gp}_\alpha)
\text{ we have } Y' \supseteq \{M: M \\
  & \in N \cap \bigcup\limits_{\beta < \alpha}
{\cE}^{\gq}_\beta \text{ and } y \in M \text{ and } M \cap
\chi^{\gp}_\alpha \in Y\}\}.
\end{array}\]
\end{itemize}
\item For an o.b. parameter ${\gp}$, we say $\bar{\bbQ}$ is an $NNR^0_\kappa$-iteration for
${\gp}$ if it is an $NNR^0_\kappa$-iteration for ${\gp}^{\bold V}$.
We say ${\gp}$ is simple if ${\gp}^{\bold V}$ is.
Similarly for $\Game$-winner, non-$\Game$-loser, etc.
\end{enumerate}
\end{definition}
The next lemma shows that each o.b. parameter leads to a canonical reasonable parameter, namely ${\gp}^{\bold V}$,  and
also shows the advantage of working with o.b. parameters  than reasonable parameters.
\begin{lemma}
\label{b45}
Assume ${\gp}$ is an o.b. (simple) parameter in the universe $\bold V$.
\begin{enumerate}
\item ${\gp}^{\bold V}$ is a (simple) reasonable parameter.

\item If $\bbP \in {\cH}(\chi^{\gp}_0)$ is a proper forcing notion (or at
least it preserves the stationarity of ${\cE}^{\gp}_\alpha$, for each
$\alpha < \ell g({\gp})$), then $\Vdash_{\bbP} ``{\gp}$ is an o.b. (simple)
parameter''.

\item If forcing with $\bbP$ adds no reals, then also $\Game$-winner,
non-$\Game$-loser, etc., are preserved.
\end{enumerate}
\end{lemma}
\begin{proof}
Straightforward.
\end{proof}
\begin{definition}
\label{b48}
Let ${\gp}$ be an o.b. parameter.
\begin{enumerate}
\item We say $\bbQ \in \bold V$ is an $NNR^0_\kappa$-forcing for ${\gp}$ or is a
${\gp}-NNR^0_\kappa$-forcing notion, when the following holds:
\begin{itemize}
\item[($\ast$)] if for some transitive class
$\bold V_0$ with ${\gp} \in \bold V_0$ and some $NNR^0_\kappa$-iteration
$\bar{\bbQ} = \langle \bbP_i,\name{\bbQ}_i:i  <
\alpha\rangle \in \bold V_0$ we have $\bold V^{\Lim(\bar{\bbQ})}_0 = \bold V$, then
we can let $\bbP_\alpha= \Lim(\bar{\bbQ}),$ $\name{\bbQ}_\alpha = \bbQ$ and get an $NNR^0_\kappa$-iteration $\bar{\bbQ}'=\langle \bbP_i,\name{\bbQ}_i:i  <
\alpha+1\rangle$.
 (i.e. $\bold V = \bold V_0[\bold G_\alpha],
\bold G_\alpha \subseteq \bbP_\alpha$ is generic over $\bold V_0$
and there is a $\bbP_\alpha$-name $\name{\bbQ}_\alpha$
such that $\langle \bbP_i,\name{\bbQ}_i:i \le \alpha\rangle$ is
an $NNR^0_\kappa$-iteration and
$\bbQ = \name{\bbQ}_\alpha[\bold G_\alpha]$).
In particular $\bbQ$ is proper and does not add reals.
\end{itemize}
\item  If we omit ``for ${\gp}$'' we mean for any ${\gp}$ which makes
sense.  Alternatively, we can put a family of ${\gp}$'s.

\item  We add ``over $x$" if this holds whenever $x \in \bold V_0$.
We can use the same definition for other versions of $\text{NNR}$.
\end{enumerate}
\end{definition}
In the next section we will present several examples of forcing notions that fit into the above definition, in the sense that they are
${\gp}-NNR^0_\kappa$ for some $2 \leq \kappa \leq \aleph_0$ and some reasonable parameter $\gp.$

\section {Examples: shooting  thin clubs} \label{examples}

In this section we present some examples that fit into our framework.  We already
know that $(< \omega_1)$-proper forcing notions are ${\gp}$-proper for standard reasonable parameter $\gp$ of length $\omega_1$ (by
Lemma \ref{a48}).  We first deal with a forcing notion which is $\text{NNR}^0_\kappa$-proper, for some  $\kappa< \aleph_0$.
Second, we deal with shooting clubs of $\omega_1$ running away
from some $C_\delta \subseteq \delta = \sup(C_\delta)$ which are small (see
\cite[Ch.XVIII, \S1]{Sh:f}).  These are the most natural non-$\omega$-proper
forcing notions which do not add reals.

\begin{definition}
\label{c3}
\begin{enumerate}
\item Let $\bar C = \langle (C_\delta,n_\delta):\delta < \omega_1,~\delta
\text{ limit} \rangle$, where $C_\delta$ is an unbounded subset of
$\delta$ of order type $\omega$ and $1 \le n_\delta < \omega$.  Let $\bar u = \langle u_\delta:\delta <
\omega_1, ~\delta \text{ limit} \rangle$, where $u_\delta \in [2n_\delta +1]
^{n_\delta}$.
Then we define $\bbQ = \bbQ_{\bar C,\bar u}$ by

\begin{equation*}
\begin{array}{clcr}
\quad \quad \bbQ_{\bar C,\bar u} = \{f:&\text{for some } \alpha < \omega_1,f
\text{ is a function from } \alpha \text{ to }\\
  & \omega \text{ such that for every limit ordinal } \delta \le
\alpha, \\
  &\text{ for some }k < 2n_\delta +1,
k \notin u_\delta \text{ and for every}\\
  &\, i \in C_\delta \text{ large enough we have } f(i) = k\}.
\end{array}
\end{equation*}
$\bbQ_{\bar C, \bar u}$ is ordered by inclusion.

\item Assume $\bar C = \langle C_\delta:\delta < \omega_1, ~\delta \text{ limit} \rangle,
C_\delta$ a closed subset of $\delta$ of order type less than $ \omega \cdot \delta$
and for $\delta_1 < \delta_2$ limit ordinals, $\sup(C_{\delta_1} \cap C_{\delta_2}) <
\delta_1$, and for limit $\delta^*$ we have $\{C_\delta \cap \delta^*:\delta
< \omega_1\}$ is countable.  Assume further that $\bar \kappa = \langle
\kappa_\delta:\delta < \omega_1, \delta \text{ limit} \rangle$ with $\kappa_\delta \in
\{2,3,\dotsc,\aleph_0\}$ and $\bar D = \langle D_\delta:\delta < \omega_1 \rangle,
D_\delta$ is a family of subsets of $\dom(D_\delta)$, such that the
intersection of any $< \kappa_\delta$ of them is non-empty.

Let $\bar f = \langle f_{\delta, x}:\delta < \omega_1, x \in
\dom(D_\delta) \rangle$ satisfy $f_{\delta,x}:C_\delta \rightarrow
\omega$ and $\bar A = \langle A_\delta:\delta < \omega_1 \rangle$ satisfy
$A_\delta \in D_\delta$.  Then  we define
$\bbQ = \bbQ_{\bar C,\bar D,\bar \kappa,\bar f,\bar A}$ as

\begin{equation*}
\begin{array}{clcr}
\quad \quad \bbQ_{\bar C,\bar D,\bar \kappa,\bar f,\bar A}=\{f:&\text{ for some } \alpha < \omega_1,f \text{ is a function from }
\alpha \text{ to } \omega \\
  &\text{ such that for every limit } \delta \le \alpha \text{~and for some} \\
  &\, x \in A_\delta \text{ we have } f_{\delta,x} \subseteq^* f
\text{ i.e. for every large } \\
  &\text{ enough }  i \in C_\delta \text{ we have } f_{\delta,x}(i) = f(i)\}.
\end{array}
\end{equation*}
$ \bbQ_{\bar C,\bar D,\bar \kappa,\bar f,\bar A}$ is ordered by inclusion.
\end{enumerate}
\end{definition}

\begin{lemma}
\label{c6}
\begin{enumerate}
\item The forcing notion $\bbQ=\bbQ_{\bar C,\bar u}$ from Definition \ref{c3}(1) is proper, does not add
reals, and is
$(< \omega_1)$-proper. It is also $\bbD$-complete for some simple
2-completeness system, hence

\begin{enumerate}
\item[$(*)$]    if $\bar{\bbQ}$ is a countable support iteration,
$\ell g(\bar{\bbQ}) = \alpha +1,\bar{\bbQ} \restriction \alpha$ is
NNR$^0_2$-iteration and $\Vdash_{\bbP_\alpha} ``\name{\bbQ}_\alpha
= \bbQ_{\bar C, \bar u}", \,\, (\bar C, \bar u)$ as above $\in \bold V$, then
 $\bar{\bbQ}$ is an NNR$^0_2$-iteration.
\end{enumerate}
\item Similarly for the forcing notion $\bbQ= \bbQ_{\bar C,\bar D,\bar \kappa,\bar f,\bar A}$ from Definition \ref{c3}(2).
\end{enumerate}
\end{lemma}

\begin{proof}
We only prove (1), as (2) can be proved in a similar way. Set $\bbQ=\bbQ_{\bar C,\bar u}$.

Let us start by showing that $\bbQ$ does not add reals. Thus suppose that $f \in \bbQ$, $\name{t}$ is a $\bbQ$-name and $f \Vdash$``$\name{t}:\omega \to 2$ is a  real''.
Let $\chi$ be a large enough regular cardinal and let $\langle N_n: n <\omega \rangle$ be a chain of countable elementary submodels of $(H(\chi), \in)$  such that $\bbQ, f, \name{t} \in N_0$  and for each $n<\omega, N_n\in N_{n+1}$. Set $N=\bigcup_{n <\omega}N_n$ and
$\delta=N \cap \omega_1$. Also for each $n<\omega$ set $\delta_n=N_n \cap \omega_1.$ Pick some $k < 2n_\delta+1, k \notin u_\delta.$ We define by induction on $n<\omega$ a sequence $\langle f_n: n<\omega \rangle$ of conditions in $\bbQ$, such that:
\begin{enumerate}
\item $f_n \in N_n,$
\item $f_0=f,$
\item for some $\delta_{n-1} < \alpha_n < \delta_n, f_n: \alpha_n \to \omega$, where $\delta_{-1}=0,$
\item $f_n \leq f_{n+1}$,
\item $\alpha_n \notin C_\delta,$
\item for all $i \in C_\delta \cap \alpha_n \setminus (\alpha_0+1), f_n(i)=k$.
\item $f_n$ decides $\name{t} \restriction n$.
\end{enumerate}
The construction can be done quite easily, noting that for each $n, C_\delta \cap \delta_n$ is bounded in $\delta_n$.
Let $f'=\bigcup_{n<\omega}f_n.$ Then $f'\in \bbQ$ extends $f$ and it decides $\name{t}$.

To show that $\bbQ$ is $<\omega_1$-proper, fix $\alpha < \omega_1$ a limit ordinal and let $\bar{N}=\langle N_\xi: \xi < \alpha \rangle$
be an increasing and continuous chain of countable elementary submodels of some $(H(\chi), \in)$ such that $\alpha, \bbQ \in N_0$ and
$\bar{N} \restriction \xi+1 \in N_{\xi+1}$. For each $\xi < \alpha$ set $\delta_\xi=N_\xi \cap \omega_1.$

Let $f \in \bbQ\cap N_0.$ By essentially the same argument as above, we can find a sequence $\langle f_\xi: \xi < \alpha  \rangle$ of conditions in $\bbQ$, with $f_\xi: \beta_\xi \to \omega$ such that:
\begin{enumerate}
\item $f_0=f$,
\item $f_\xi \leq f_{\xi+1}$,
\item if $\xi$ is a limit ordinal, then $f_\xi=\bigcup_{\zeta<\xi} f_\zeta,$
\item $\delta_\xi < \beta_{\xi+1} < \delta_{\xi+1}$,
\item $f_{\xi+1} \in N_{\xi+1}$ is $(\bbQ, N_\xi)$-generic,
\item $f'=\bigcup_{\xi< \alpha}f_\xi$ is a condition.
\end{enumerate}
Then $f'\geq f$ is $(\bbQ, N_\xi)$-generic for every $\xi < \alpha.$

Let us show that $\bbQ$ is $\bbD$-complete for some simple
2-completeness system. Let $\theta$ be large enough regular and let $\bbD$ be a function whose domain consists of those pairs
$(N, f)$ where $N \prec (H(\theta), \in)$ is countable, $\bbQ \in N$ and $f \in \bbQ \cap N.$
Now suppose that $(N, f) \in \dom(\bbD)$. Set
\[
\bbD(N, f) = \{A_x: x \text{~is a finitary relation on~}N     \},
\]
where for any such $x$,
\[
A_x=\{  \bold G \in \text{Gen}(N, \bbQ, f): N \cup \cP(N) \models \Psi(x, \bold G, N, \bbQ, f)           \}
\]
and the formula $\Psi$ says that ``if $x=(y, k),$ where $y$ is an $\omega$-sequence cofinal in $\delta=N \cap \delta$ and
 $k< 2n_\delta+1, k \notin u_\delta$, then $\bigcup \bold G \restriction y$ is eventually  equal to $k$.

Let us show that
 $\bbD$ is a simple
2-completeness system. Thus suppose that $x_0, x_1$ are given, and we have to show that $A_{x_0} \cap A_{x_1}$ is non-empty. The only non-trivial case is when $x_0, x_1$ satisfy the hypotheses of the formula $\Psi.$ Thus suppose that $x_0=(y_0, k_0)$ and $x_1=(y_1, k_1)$, where
$y_0, y_1$ are  $\omega$-sequences cofinal in $\delta$ and $k_0, k_1 < 2n_\delta+1, k_0, k_1 \notin u_\delta$.
If $y_0 \cap y_1$ is cofinal in $\delta$, then we must have $k_0=k_1$ and so by the previous arguments we can find an increasing sequence
$\langle  f_n: n<\omega      \rangle$ of extensions of $f$ in $N$ such that $f^*=\bigcup\limits_{n<\omega}f_n$ gives rise to a filter
$\bold G \in \text{Gen}(N, \bbQ, f)$ such that $f^* \restriction (y_0 \cup y_1)$ is eventually equal to $k_0$ and we are done.
Otherwise, $y_0 \cap y_1$ is bounded in $\delta,$ so for some $\eta < \delta, y_0 \cap y_1 \subseteq \eta$. By enlarging $\eta$
we may also assume that $\dom(f)  \subseteq \eta.$ Again, by similar arguments as above, it is not difficult to build an increasing
sequence $\langle  f_n: n<\omega      \rangle$ of extensions of $f$ in $N$ such that $f^*=\bigcup\limits_{n<\omega}f_n$ gives rise to a filter
$\bold G \in \text{Gen}(N, \bbQ, f)$ such that
$f^* \restriction (y_0 \setminus \eta)$ is eventually equal to $k_0$ and $f^* \restriction (y_1 \setminus \eta)$ is eventually equal to $k_1$.
Then $\bold G \in A_{x_0} \cap A_{x_1}$ and we are done.

It is now clear that $\bbQ$ is $\bbD$-complete, which completes the proof.
\end{proof}


\begin{definition}
\label{c9}
Assume $\bar C = \langle C_\delta:\delta < \omega_1,~\delta \text{ limit } \rangle$,
where $
C_\delta$ is an unbounded subset of the limit ordinal
$\delta$ (think of the case $C_\delta$ of order type $< \delta$ but not
necessarily).  Let

\begin{equation*}
\begin{array}{clcr}
\bbQ_{\bar C} = \{c:&\text{ for some } \alpha < \omega_1,c
\text{ is a closed subset of } \alpha \\
  &\text{ and for every limit ordinal } \delta \le \alpha
\text{ we have} \\
  & ~\delta = \sup(c \cap \delta)
\Rightarrow c \cap C_\delta \text{ is bounded in } \delta\}.
\end{array}
\end{equation*}
Order $\bbQ_{\bar C}$ by
\[
c_1 < c_2 \iff c_1 \text{ is an initial segment of } c_2.
\]
\end{definition}

\begin{remark}
\label{c10}
\begin{enumerate}
\item For more information about the above forcing notion see \cite[Ch.XVIII, \S1, 1.9]{Sh:f}. Note that $\bbQ_{\bar C}$ may be
non-$\omega$-proper.

\item Note that $c \in \bbQ_{\bar C}$ is a closed subset of $\alpha$, but not
necessarily a closed subset of $\omega_1$.
\end{enumerate}
\end{remark}
In general the forcing notion $\bbQ_{\bar C}$ might be trivial, say for example when for every limit ordinal $\delta,~C_\delta=\delta.$ We are interested in the cases that this forcing notion is non-trivial, and we first deal with the simple case of $\otp(C_\delta)
= \omega$.
\begin{lemma}
\label{c12}
Assume ${\gp}$ is a simple reasonable parameter, $\bar C$ is as in
Definition \ref{c9} and  $\bigwedge\limits_{\delta} \otp(C_\delta) = \omega$. Let
$f \in {\cF}^{\gp}$ be defined as $f(0) = 0$ and  $f(\beta) = 1+ \beta$ for $\beta>0$. Then
$\bbQ_{\bar C}$ is $({\gp},f)$-proper.
\end{lemma}
In Lemma \ref{c18} we prove a stronger result, which includes the above lemma as a very special case.
 Note that if $\bigwedge\limits_{\delta} \otp(C_\delta) <
\delta$, then we can split the analysis by restricting ourselves to
$\{N:N \cap \omega_1 \in S_\gamma\}$, where
$\gamma < \omega_1$ is such that $S_\gamma = \{\delta:\otp(C_\delta)
= \gamma\}$ is stationary.
\begin{lemma}
\label{c18}
\begin{enumerate}
\item Assume

\begin{enumerate}
\item[$(a)$]   ${\gp}$ is a simple reasonable parameter,

\item[$(b)$]   $f \in {\cF}^{\gp}$ and ${\gp}$ is non-$\Game_f$-loser,

\item[$(c)$]   $\gamma(*) < \omega_1$,

\item[$(d)$]  $\bar C = \langle C_\delta:\delta < \omega_1 \rangle,
C_\delta \subseteq \delta = \sup(C_\delta)$.

\item[$(e)$]  $\otp(C_\delta) \le \omega^{\gamma(*)}$ for every $\delta$,

\end{enumerate}
Define $g \in {\cF}^{\gp}$ by recursion as
\begin{itemize}
\item $g(0) = 0,$

\item $g(1) = f(1) + \gamma(*),$

\item $g(\alpha +1) = f(g(\alpha)) + \gamma(*) +1$, for $\alpha > 0,$

\item for limit ordinals
$\alpha, g(\alpha) = \underset{\beta < \alpha}{\sup} g(\beta)$.
\end{itemize}
Then  for every $\alpha < \ell g({\gp})$, the forcing notion
$\bbQ = \bbQ_{\bar C}$ is $({\gp},\alpha,g(\alpha))$-proper (version 1).

\item In part (1), we can get ``version 2'' of $({\gp},\alpha,g(\alpha))$-properness when the following is satisfied: if $g(\delta) = \delta,N \in
{\cE}^{\gp}_\alpha,C \subseteq \omega_1 \cap N = \sup(C),\otp(C) \le
\omega^{\gamma(*)}$
and $Y \in D^{\gp}_\alpha(N)$, then
$Y' = \{M \in Y:M \cap \omega_1 \notin C\} \in D^{\gp}_\alpha(N)$.

\item If we weaken clause $(b)$ to $(b)_{f,g}$, where
\begin{enumerate}
\item[] $(b)_{f,g}:  f \in {\cF}^{\gp}$, $f(f(\alpha)) = f(\alpha)$ for
every $\alpha < \ell g({\gp})$ and ${\gp}$ is a

\item[] \quad\quad\quad non-$\Game'_f$-loser,
\end{enumerate}
 then  for
$\alpha < \ell g({\gp})$, the forcing notion $\bbQ_{\bar C}$ is
$({\gp},\alpha,f(\gamma(*) + \alpha))$-proper.
\end{enumerate}
\end{lemma}
\begin{proof}
We only prove clause (2). Note that version 2 is harder to prove,
and using the extra freedom, we can avoid the need for the extra
assumption from (2) to prove (1).
First observe that
\begin{enumerate}
\item[$(*):$]   for each $\alpha < \omega_1$ the set ${\cI}^*_\alpha =
\{p \in \bbQ_{\bar C}:\text{ there is } \beta \in p
\text{ which is } \ge \alpha\}$ is an open dense  subset of $\bbQ_{\bar C}$.
\end{enumerate}
 We prove (2) by induction on $\alpha$. Let $\beta = g(\alpha)$.
Let $N \in {\cE}^{\gp}_\beta$ be countable with $\bbQ_{\bar C},\alpha,\beta, f, g \in N$, $p \in N \cap \bbQ_{\bar C}$,
 and
suppose $Y \in D^{\gp}_\beta(N)$ is given.  Let $\delta = \delta_N = N \cap \omega_1$.

\underline{Case 1:  $\alpha = 0$.}
In this case, we are reduced to show that
$\bbQ_{\bar C}$ is proper.  Let $\langle {\cI}_n:n < \omega
\rangle$ list the dense open subsets of $\bbQ_{\bar C}$ that belong
to $N$. We shall choose by induction on $n<\omega$, a condition $p_n$ such that:
\begin{enumerate}
\item[$(i)$]   $p_0 = p,$

\item[$(ii)$] for each $n$, $p_n \in N$,

\item[$(iii)$]   $p_n \le p_{n+1} \in {\cI}_n$,

\item[$(iv)$]   the set $p_{n+1} \cup \{\sup p_{n+1}\} \backslash
(p_n \cup \{\sup p_n\})$ is disjoint from $C_\delta$.
\end{enumerate}

Set $p_0=p$.  Now assume $p_n$ has
been chosen and we shall choose $p_{n+1}$ as requested.  Let
 $F_n \in N$ be a function with domain $\bbQ_{\bar C}$ such that for all
$q \in \bbQ_{\bar C},~ q \le F_n(q) \in {\cI}_n$.

For $\alpha < \omega_1$ let $q^{[\alpha]} = q \cup \{\sup(q)\} \cup
\{\sup(q) + 1 + \alpha\}$, so clearly $q \le q^{[\alpha]} \in \bbQ_{\bar C}$
and the function $(q,\alpha) \mapsto q^{[\alpha]}$ belongs to $N$.  Define
a function $H:\omega_1 \rightarrow \omega_1$ by $H(\alpha) =
\sup(F_n(p^{[\alpha]}_n))$. Clearly it is well defined and belongs
to $N$. Let
$$C = \{\beta < \omega_1:
\beta \text{ a limit ordinal, } \omega \beta = \beta,
(\forall \alpha < \beta)(H(\alpha) < \beta) \text{ and } \sup(p_n) < \beta\}.$$
It is easily seen that $C$ is a club of $\omega_1$ which belongs to $N$ and $\gamma(*) \in N$,
hence we can find $\beta^* \in C$ such that $\otp(\beta^* \cap C)$ is
divisible by $\omega^{\gamma(*)}$. But $\otp(C_\delta \cap \beta^*)
< \omega^{\gamma(*)}$, hence for some $\beta \in C$
we have $\sup(C_\delta \cap \beta) < \beta$.  Let $p_{n+1} = F_n
(p^{[\sup(C_\delta \cap \beta)+1]})$.

Set
$q=\bigcup\limits_{n < \omega} p_n$. Then for each
$\alpha < N \cap \omega_1,{\cI}^*_\alpha \in \{{\cI}_n:n < \omega\}$, hence
$$\exists \beta \left( \beta \in q
 \text{ and } \alpha \le \beta < \omega_1\right).$$
It follows that $q$ is $(N, \bbQ_{\bar C})$-generic.

\underline{Case 2:  $\alpha =1$.}
Set $Y' = \{M \in Y:M \cap \omega_1 \notin C_\delta\}$. It is clear that $g(\delta)=\delta$ (as $g \in N$ and $\delta=N \cap \omega_1$), hence by the hypotheses in clause (2) we have $Y'\in
D^{\gp}_{g(1)}(N)$.
Let $\langle {\cI}_n:n < \omega \rangle$
list the dense open subsets of $\bbQ_{\bar C}$ which belong to $N$
and let $\delta = \lim\limits_{n < \omega} \alpha_n$ where
$\langle \alpha_n: n < \omega \rangle$ is increasing.  We now simulate a strategy
for the challenger in the game $\Game_{\alpha,\beta}(N)$, where in the $n$-th move, we let the challenger to choose $Z_n=\emptyset$ (so the chooser has to use $Y_n = \emptyset$ as well)  and at the end of
the $n$-th move,   the challenger also chooses
$p_{n+1} \in \bbQ_{\bar C} \cap N$ such that:
\begin{itemize}
\item $p_0 = p,$
\item $p_n \le p_{n+1} \in
{\cI}_n$,
\item $p_{n+1}$ is
$(M_n,\bbQ_{\bar C})$-generic and $\sup(p_{n+1}) > \alpha_n$,
\item the set
$(p_{n+1} \cup \{\sup p_{n+1}\}) \backslash (p_n \cup \{\sup p_n\})$
is disjoint to $C_\delta$.
\end{itemize}
  This is possible by Case 1 and its proof,
because $M_n \cap \omega_1 \notin C_\delta$ which holds as $M_n \in
Y'$.  As this is a legal strategy for the challenger,  it cannot
be a winning strategy, hence for some such play the chooser wins,
hence $\{M_n:n < \omega\} \in D^{\gp}_1(N)$.  Now
$q = \bigcup\limits_{n < \omega} p_n$ is well defined, and
$\sup(q) = \delta$ and $q \cap C_\delta \subseteq p \cup
\{\sup(p)\}$ and $q \Vdash_{\bbQ_{\bar C}}
``\{M_n:n < \omega\} \subseteq
{\cM}[\name{\bold G}_{\bbQ_{\bar C}},N]"$,  so $q$ is as required as the
chooser has won the play.

\underline{Case 3:  $\alpha > 1,\alpha$ successor.}
The proof is similar to the Case 2, only we use the induction hypothesis instead of using
Case 1.

\underline{Case 4:  $\alpha$ a limit ordinal.}
The proof is again similar to the Case 2.

This  completes the induction hypothesis and hence the proof of clause (2) of the lemma.
\end{proof}

\begin{definition}
\label{c21}
Suppose $S \subseteq \omega_1$ is stationary, ${\cD}_{\omega_1}$ is the club filter
 on $\omega_1$ and  $f \in {}^{\omega_1}\omega_1$.
\begin{enumerate}
\item  We say $f$ is a $({\cD}_{\omega_1} + S,
\gamma)$ function, when $S \Vdash_{({\cD}^+_{\omega_1},\supseteq)}$ ``in
$\bold V[\name{\bold G}],\{x \in \bold V^{\omega_1}/\name{\bold G}:
\bold V^{\omega_1}/\name{\bold{G}} \models$`` $x$ is an ordinal $< f/ \name{\bold G}\}$
has order type $\gamma$''.

\item Assume $\bar C = \langle C_\delta:\delta < \omega_1 \rangle$, where
$C_\delta$ is an unbounded subset of $\delta$.
We define, by induction on $\gamma$, when ``$\bar C$ obeys $f$ on $S$'' for $f \in {}^{\omega_1}\omega_1$
which is a $({\cD}_{\omega_1} + S,\gamma)$ function:
\begin{itemize}
\item if
$\gamma < \omega_1$, this means
$$\{\delta \in S:\otp(C_\delta) \le
\omega^{1+f(\delta)}\} = S \mod {\cD}_{\omega_1}.$$
\item
if $\gamma \geq \omega_1$, it means that for some $g:\omega_1 \rightarrow \omega_1$
and pressing down function $h$ on $S$, for
every $\zeta < \omega_1$ for which $h^{-1}\{\zeta\}$ is stationary,
for some $\beta < \gamma$ and $f_\beta$, a
$({\cD} + h^{-1}\{\zeta\},\beta)$ function, we have
$\langle C_{g(\delta)} \cap \delta:\delta \in h^{-1}\{\zeta\}\rangle$
 obeys $f_\beta$.
 \end{itemize}
\end{enumerate}
\end{definition}
The next lemma can be proved as in Lemma \ref{c18}.
\begin{lemma}
\label{c24}
Assume
\begin{enumerate}
\item[$(a)$]   ${\gp}$ is a simple reasonable parameter such that
$\ell g({\gp})$ is of uncountable cofinality,
\item[$(b)$]   $S \in {\cD}^+_{\omega_1}$ and $N \in \bigcup\limits_{\alpha}
{\cE}^{\gp}_\alpha \Rightarrow N \cap \omega_1 \in S$,

\item[$(c)$]   ${\gp}$ is a non-$\Game_{\alpha,\alpha}$-loser
(or just non $\Game'_{\alpha,\alpha}$-loser) for all
$\alpha \in C^*,$ where $C^*$ is a club of $\ell g({\gp})$
with $0 < \min(C^*)$,

\item[$(d)$]   $\bar C$ obeys $f$ on $S$ which is a
$({\cD}_{\omega_1} +S,\gamma)$-function,

\item[$(e)$]   for all $\alpha, g(\alpha) = \min(C^* \backslash \alpha)$.
\end{enumerate}
Then  $\bbQ_{\bar C}$ is $({\frak p},g)$-proper.
\end{lemma}
Let us now give
another example which fits into our general framework (see also \cite[Ch.XVIII]{Sh:f}). Recall that a filter $\cD$ on a countable set is called a $P$-filter if
it contains all co-finite sets and if $A_n \in \cD$
for $n < \omega$, then for some $A \in \cD$ and all $n<\omega$ we have
$|A \backslash A_n| < \aleph_0$
\begin{definition}
\label{c27}
\begin{enumerate}
\item We say $\bar{\cD} = \langle {\cD}_\delta:\delta < \omega_1,~\delta
\text{ limit} \rangle$ is an $\omega_1$-filter-sequence if:

\begin{enumerate}
\item[$(a)$]   ${\cD}_\delta$ is a filter on $\delta$, containing the
co-bounded subsets of $\delta$,

\item[$(b)$]   ${\cD}_\delta$ is a $P$-filter and some $C_\delta
\in \cD_\delta$ has order type $\omega$,

\item[$(c)$]   for every club $C \subseteq \omega_1$ and $\alpha <
\omega_1$, the set $A^\alpha_C[\bar{\cD}]$ is stationary, where,
by induction on $\alpha$,  we define $A^\alpha_C[\bar{\cD}]$
by
$A^\alpha_C[\bar{\cD}] = \{\delta < \omega_1:\delta
\text{ is a limit ordinal, } \delta \in C
\text{ and for every } \beta < \alpha \text{ we have}$
$\delta = \sup(\delta \cap A^\beta_C[\bar{\cD}])$, moreover
$\delta \cap A^\beta_C[\bar{\cD}] \in {\cD}_\delta\}$.
\end{enumerate}
\item A reasonable parameter ${\gp}$ obeys $\bar{\cD}$, if
for each $\alpha < \ell g({\gp})$ and $N \in {\cE}^{\gp}_\alpha$,
$\bar{\cD} \in N$ and we have

\begin{equation*}
\begin{array}{clcr}
{D}^{\gp}_\alpha(N) = \{Y:&Y \subseteq N \cap \bigcup\limits_{\beta
  < \alpha} {\cE}^{\gp}_\beta \text{ is } \gp\text{-closed  and if } \alpha > 0, \text{ then there are}\\
  &\text{}\bar \beta = \langle \beta_n:n < \omega \rangle \text{ and }\bar M =
\langle M_n:n < \omega \rangle \text{ satisfying: } \\
  &(a) \quad \beta_n \in N \cap \alpha,  \\
  &(b) \quad \text{ either for all}~ n, ~\alpha = \beta_n +1 \\
  &\qquad  \text{ or } \beta_n < \beta_{n+1},
          \sup_{n < \omega} \beta_n = \sup(\alpha \cap N), \\
  &(c) \quad M_n \in Y \cap {\cE}^{\gp}_{\beta_n},M_n \in M_{n+1}, \\
  &(d) \quad \bigcup\limits_{n < \omega} M_n = N \cap
\bigcup\limits_{\beta \in \alpha \cap N} {\cH}(\chi^{\gp}_\beta), \\
  &(e) \quad \{M_n \cap \omega_1:n < \omega\} \in {\cD}_{N \cap \omega_1}\}.
\end{array}
\end{equation*}

\item A forcing notion $\bbQ$ is a $\bar{\cD}-NNR^0_\kappa$-forcing
if for every reasonable parameter ${\gp}$ which obeys $\bar{\cD}$, $\bbQ$ is an $NNR^0_\kappa$-forcing over ${\gp}$
(see Definition \ref{b48}).

\item For a $P$-filter ${\cD}$ on $\omega$, we say a reasonable parameter
${\gp}$ obeys ${\cD}$ if for every $N \in {\cE}^{\gp}_\alpha$

\begin{equation*}
\begin{array}{clcr}
D^{\gp}_\alpha(N) = \{Y:&Y \subseteq N \cap \bigcup\limits_{\beta < \alpha}
{\cE}^{\gp}_\beta \text{ is } \gp\text{-closed and if } \alpha > 0,
\text{ then there are} \\
  & \bar \beta,\bar M \text{ satisfying items }
(a), (b) \text{ and } (d) \text{ of clause } (2) \text{ and} \\
  &(e) \quad \{n:M_n \in Y\} \in {\cD}\}. \\
\end{array}
\end{equation*}

\item In parts (1) - (4) above, we may replace the word ``filter'' by ``ultrafilter''
if the ${\cD}_\alpha$'s are ultrafilter.
\end{enumerate}
\end{definition}

\begin{lemma}
\label{c30}
\begin{enumerate}
\item If $\diamondsuit_{\aleph_1}$ holds, then there is
an $\omega_1$-ultrafilter sequence.

\item If $\bar{\cD}$ is an $\omega_1$-filter sequence and
$\langle (\chi_\alpha,{\cE}_\alpha):\alpha < \omega_1 \rangle$ is as  in
Definition \ref{a3}, then there is a
reasonable parameter ${\gp}$ of length $\omega_1$ obeying $\bar{\cD}$
which is a non-$\Game$-loser. Furthermore, for all $\alpha<\omega_1$, $\chi^{\gp}_\alpha = \chi_\alpha$ and
${\cE}^{\gp}_\alpha = {\cE}_\alpha$.

\item If $\diamondsuit_{\aleph_1}$ holds, $\langle (\chi_\alpha,{\cE}_\alpha):
\alpha < \omega_1 \rangle$ is as above and ${\cD}$ is
a $P$-filter on $\omega$, then  some reasonable parameter ${\gp}$
of length $\omega_1$ is $P$-filter like, non-$\Game_{\id}$-loser with $\chi^{\gp}_\alpha
= \chi_\alpha,{\cE}^{\gp}_\alpha = {\cE}_\alpha$.
Similarly for ultrafilters.

\item Instead of $\diamondsuit_{\aleph_1}$ it is enough to assume CH and
that for some $\langle C_\delta:\delta < \omega_1 \text{ limit} \rangle$
and some normal filter $D$ on $\omega_1$, and for every club $C$ of $\omega_1,
\{\delta:\delta > \sup (C_\delta \backslash C)\} \in D$.
\end{enumerate}
\end{lemma}

\begin{proof}
We only prove (1) and (2), the other parts can be proved similarly.

(1). Let
$\langle S_\delta: \delta < \omega_1 \rangle$ be a $\diamondsuit_{\aleph_1}$-sequence, where each
$S_\delta \subseteq \delta$ and let
$\langle  E_\delta:\delta < \omega_1  \text{ limit}   \rangle$
be a ladder system, where each $E_\delta \subseteq \delta$ has order type $\omega$ and $E_\delta=S_\delta$ if $S_\delta$ is an $\omega$-sequence
cofinal in $\delta.$
Let $\bar{\cD}=\langle \cD_\delta: \delta < \omega_1 \text{ limit}       \rangle$,
where $\cD_\delta$ is any $P$-ultrafiler on $\delta$  extending $\cD_{^\delta}^{\text{c}} \cup \{ E_\delta\}$,
where $\cD_{\delta}^{\text{c}}$ is the filter of co-bounded subsets of $\delta$.
We show that $\bar{\cD}$ is as required. Items (a) and (b) of Definition \ref{c27}(1) are clearly satisfied.
Let us prove clause (3). Thus suppose that $C\subseteq \omega_1$ is a club, $\alpha < \omega_1$, and suppose by induction that for all
$\beta<\alpha$, the set $A^{\beta}_C(\bar{\cD})$ is stationary. We show that $A^{\alpha}_C(\bar{\cD})$
is stationary as well. We may assume that $C$ only contains limit ordinals. Then
\[
A^{\alpha}_C(\bar{\cD}) = C \cap \bigcap\limits_{\beta<\alpha}\{ \delta:\delta=\sup(\delta \cap A^{\beta}_C(\bar{\cD}))     \} \cap \{\delta: \forall \beta < \alpha,~\delta \cap A^{\beta}_C(\bar{\cD}) \in \cD_\delta    \}.
\]
It follows from the induction hypothesis that the set $C \cap \bigcap\limits_{\beta<\alpha}\{ \delta:\delta=\sup(\delta \cap A^{\beta}_C(\bar{\cD}))     \}$ is a club. Suppose by contradiction that the set $A=\{\delta: \forall \beta < \alpha,~\delta \cap A^{\beta}_C(\bar{\cD}) \in \cD_\delta    \}$
is non-stationary. Thus for some club $D \subseteq C$ and for all $\delta \in D$, there exists some $\beta_\delta < \alpha$
such that $\delta \cap A^{\beta_\delta}_C(\bar{\cD}) \notin \cD_\delta$. As $\alpha < \omega_1,$ it follows from F\"{o}dor's lemma that there are a
stationary set $S \subseteq D$ and some  fixed
$\beta_*<\alpha$ such that
$$\delta \in S \implies \delta \cap A^{\beta_*}_C(\bar{\cD}) \notin \cD_\delta$$
On the other hand, by the $\diamondsuit_{\aleph_1}$-assumption, the set
\[
T=\{  \delta \in S:      \delta \cap A^{\beta_*}_C(\bar{\cD})=S_\delta                   \}
\]
is stationary. We may assume without loss of generality that for all $\delta \in T,  \delta \cap A^{\beta_*}_C(\bar{\cD})$ has order type $\omega.$ But then
$$\delta \in T \implies \delta \cap A^{\beta_*}_C(\bar{\cD})=S_\delta =E_\delta \in \cD_\delta,$$
which is a contradiction.

(2) Define $\gp$ of length $\omega_1$ such that for all $\alpha < \omega_1$,
\begin{itemize}
\item $\chi^{\gp}_\alpha = \chi_\alpha$,
\item $R^{\gp}_\alpha = [\cH((\bigcup_{\beta < \alpha}\chi_\beta)^+)]^{\leq \aleph_0}$,

\item ${\cE}^{\gp}_\alpha = {\cE}_\alpha$,

\item for $N \in {\cE}_\alpha,$ $D^{\gp}_{\alpha}(N)$ is defined as in Definition \ref{c27}(2).
\end{itemize}
Then $\gp$ is as required.
\end{proof}

\begin{lemma}
\label{c33}
\begin{enumerate}
\item If ${\cD}$ is
a $P$-filter on $\omega$ (or $P$-ultrafilter on $\omega$) and ${\gp}$ is a reasonable parameter obeying ${\cD}$,
then
for some $\bar{\cD},\bar{\cD}$ is an $\omega_1$-filter-sequence
(or $\omega_1$-ultrafilter-sequence) and ${\gp}$ obeys $\bar{\cD}$.

\item If ${\gp}$ is a $P$-point filter (or ultrafilter),
then ${\gp}$ is a non-$\Game$-loser.
\end{enumerate}
\end{lemma}

\begin{proof}
(1) For each limit ordinal $\delta< \omega_1$ fix a bijection $f_\delta: \omega \leftrightarrow \delta$
and set $\cD_\delta=\{f''[X]: X \in \cD   \}$. Then $\bar{\cD}$
is as required.

Proof of (2) is essentially similar to the proof of Lemma \ref{a18}.
\end{proof}
We now consider the case where the order type of the club sets $C_\delta$is higher  than $\omega$.

\begin{lemma}
\label{c39}
\begin{enumerate}
\item Assume

\begin{enumerate}
\item[$(a)$] $\kappa \leq \omega$ and  $\bar C = \langle C_{\delta,\ell}:\ell < k_\delta,\delta <
\omega_1 \text{ limit}\rangle$, where $1 + \kappa \le k_\delta \le \omega,
C_{\delta,\ell}$ is a closed unbounded subset of $\delta$ and
$\ell < m < k_\delta \Rightarrow C_{\delta,\ell} \cap C_{\delta,m}
= \emptyset$,

\item[$(b)$]   $\bbQ = \bbQ_{\bar C} = \{C:C$ is a closed bounded subset of
$\omega_1$ such that for every limit $\delta < \sup(C)$,  and for every $\ell
< k_\delta$ except $< 1 + \kappa$ many, $\delta \le \sup
(C \cap C_{\delta,\ell})\}$,\footnote{i.e., $\{\ell
< k_\delta:  \delta > \sup
(C \cap C_{\delta,\ell}) \}$ has size $< 1 + \kappa$.}

\item[$(c)$]   ${\gp}$ is a reasonable parameter,
obeying the $P$-ultrafilter ${\cD}$.
\end{enumerate}
Then  $\bbQ$ is a ${\gp}-$NNR$^0_\kappa$ forcing notion.

\item In part (1), if we add:
\begin{enumerate}
\item[] $N \in {\cE}^{\gp}_\alpha$ and $D^{\gp}_\alpha(N) =
\{Y:\{n:M_n \in Y\} \in {\cD}_N\},{\cD}_N$ a $P$-ultrafilter
 and $\ell < k_\delta \Rightarrow \{n < \omega:M_n \cap \omega_1
\in C_{\delta,\ell}\} = \emptyset \mod {\cD}_N$,
\end{enumerate}
then  we can
allow $C_{\delta,0} = C_{\delta,1}$.

\item Assume

\begin{enumerate}
\item[$(a)$]   $D_\delta$ is a family of subsets of $\dom(D_\delta)$, the
intersection $Y$ of any $<1 + \kappa$ of them satisfies,

\item[${{}}$]   $(*) \qquad\qquad \quad \exists n(\exists y_1,\dotsc,y_n \in Y)
[\delta > \sup (\bigcap\limits^n_{\ell =1} C_{\delta,y_\ell})]$,

\item[$(b)$]   $\bar C = \langle C_{\delta,x}:x \in \dom(D_\delta)$
and $\delta$ is a limit ordinal $< \omega_1 \rangle$,

\item[$(c)$]  $\langle C_{\delta,x}:x \in \dom(D_\delta) \rangle$
is a sequence of pairwise disjoint subsets of $\delta$,

\item[$(d)$]  $\bar X \in \prod\limits_{\delta < \omega_1} \dom(D_\delta)$,

\item[$(e)$]   $\bbQ_{\bar C,\bar X,\bar D} = \{C:C$ is a closed bounded
subset of $\omega_1$ such that for every limit $\delta \le \sup(C)$ we have
$(\exists x \in X_\delta)(\delta > \sup(C \cap C_{\delta,x}))\}$
ordered by  end extension,

\item[$(f)$]   $\bar{\cD}$ is a $P$-ultrafilter sequence,

\item[$(g)$]   ${\gp}$ is a reasonable parameter which obeys $\bar{\cD}$.
\end{enumerate}
Then $\bbQ$ is a ${\gp}-NNR^0_\kappa$ forcing notion.
\end{enumerate}
\end{lemma}

\begin{proof}
We  prove clause (1), as other items can be proved in a similar way. So let $\bold V_0$ be some transitive class with $\gp \in \bold V_0$  and suppose that $\bar{\bbQ} \in \bold V_0$ is an NNR$^0_\kappa$-iteration with $\bbP_\alpha=\lim(\bar{\bbQ})$ such that $\bold V = \bold V^{\bbP_\alpha}_0$
and $\Vdash^{\bold V_0}_{\bbP_\alpha} ``\name{\bbQ}_\alpha = \bbQ_{\name{\bar C}}$ is as
above". Set $\bbP_{\alpha+1}=\bbP_\alpha \ast \name{\bbQ}_\alpha$.  We have to show that items (a)-(d) of Definition \ref{a24} (for $\kappa=\aleph_0$) or \ref{a39} (for $\kappa < \aleph_0$) are satisfied. We only check clause (d), as other items are easier to prove. Suppose that in $\bold V_0,$
\begin{itemize}
\item $ N_0 \in \bigcup\limits_{\beta < \alpha}{\cE}^{\gp}_\beta,$
\item $N_0 \in N_1 \in
{\cE}^{\gp}_\alpha$,
\item $N_1 = \{M_n:n < \omega\}$,
\item $\cD_{N_0}$ is a $P$-ultrfiler as in Definition
\ref{c39}(2),

\item $ \bold G_m \subseteq N_1
\cap \bbP_\alpha$ is generic over $N_1$ for $m < k < 1 + \kappa,$

\item $\bigwedge\limits_{m< k}\bold [ \bold{G}_m \cap N_0=\bold G^* ]$,

\item  $p \in \bbP_{\alpha +1} \cap N_0$ is such that $p \restriction \alpha \in \bold G^*$,

\end{itemize}
 Clearly, \wilog \, $\langle M_n:n <
\omega \rangle, ~\cD_{N_0} \in N_1$.
For $\iota=0,1$ set $\delta_{N_\iota}=N_\iota \cap \omega_1$.
So for each $\ell < k_{\delta_{N_0}}$, we have $\name C_{\delta_{N_0},\ell}
[\bold G_m]$ is a closed subset of $\delta$ and for $\ell_1 < \ell_2 <
k_{\delta_{N_0}}$ we have $\name C_{N_0 \cap \omega_1,\ell_1}[\bold G_m]
\cap \name C_{N_0 \cap \omega_1,\ell_2}[\bold G_m] = \emptyset$.
So for some $\ell(m) \in \{0,1,\dotsc,k_{\delta_{N_0}}-1\}$ we have
$$\ell \ne \ell(m) \Rightarrow \name C_{N_0 \cap \omega_1,\ell}
[\bold G_m] = \emptyset \mod \cD_{N_0}.$$
Now let

$\hspace{1.5cm}$$B = \{n:$ if $\ell < k_{\delta_{N_0}},\ell \notin \{\ell(m):
m < k\}$ and $\ell < n$, then

$\hspace{2.8cm}$ $M_n \cap \omega_1 \notin
\name C_{\delta_{N_0},\ell}[\bold G_m]$ for $m < k$ and
$p \in M_n\}$.

Then $B$ belongs to $N_1 \cap \cD_{N_0}$.   Let $B = \{n_i:i <
\omega\}$ be an  increasing enumeration of $B$ and let $\langle {\cI}_n:n < \omega
\rangle$ list the dense open subsets of $\bbP_{\alpha +1}$ which belong
to $N_0$.  We  choose $p_i$, by induction on $i < \omega$, such that:

\begin{enumerate}
\item[(a)]  $p_i \in N_{1} \cap P_{\alpha +1}$,

\item[(b)]  $p_i \restriction \alpha \in
\bigcap\limits_{m < k} \bold G_m$,

\item[(c)]  $p_i \in \bigcap\{{\cI}_n:n < n_i,{\cI}_n \in N_{1}
\text{ and } i > 0\}$,

\item[(d)]  $p \le p_i$,

\item[(e)]  $p_i \le p_{i+1}$,

\item[(f)]  $p_{i+1} \backslash p_i \text{ is disjoint to } \bigcup
\{\name C_{\delta_{N_0},\ell}[\bold G_m]:\ell < k_{\delta_{N_0}},\ell <
n_i \text{ and } m < k \Rightarrow \ell \ne \ell(m)\}$.
\end{enumerate}

This is possible as, for each $i<\omega$,
  $$\bigcup
\{\name C_{\delta_{N_0},\ell}[\bold G_m] \cap M_{n_i} \cap \omega_1:\ell < k_{\delta_{N_0}},\ell <
n_i \text{ and } m < k \Rightarrow \ell \ne \ell(m)\}$$
 is a bounded subset of $M_{n_i}\cap \omega_1$.

Set
\[
\bold G^{**}=\{q \in \bbP_{\alpha+1} \cap N_0:  \exists i<\omega (p_i \leq q)   \}.
\]
Then $\bold G^{**} \subseteq  \bbP_{\alpha+1} \cap N_0$ is generic over $N_0, p \in \bold G^{**}$ and if we set
 $q=\bigcup\limits_{i<\omega}p_i,$ then $q$ witnesses that $\bold G^{**}$ has an upper bound in $\bbP_{\alpha+1}/ \name{G}_{\bbP_{\alpha}}$.
\end{proof}
The proof of the next lemma is similar to the above proofs.

\begin{lemma}
\label{c42}
\begin{enumerate}
\item The forcing notion $\bbQ_{\bar{C}}$ from Lemma \ref{c39}(1) is
$NNR^0_{\aleph_0}$-forcing notion for every ${\gp}$,
non-$\Game_{\id}$-loser.

\item If $\bar{\cD}$ is  an
$\omega_1$-filter sequence and ${\gp}$ is a reasonable parameter obeying $\bar{\cD}$, then any $(< \omega_1)$-proper forcing notion is
${\gp}$-proper.
\end{enumerate}
\end{lemma}

\section {Second preservation of not adding reals} \label{second}
In this section we  present our second preservation theorem.
We shall concentrate on the simple case.
\begin{definition}
\label{d3}
Let ${\gp}$ be a reasonable parameter and let $\bar{\bbQ}=\langle \bbP_i, \name{\bbQ}_i: i<  \ell g(\bar{\bbQ}) \rangle$
be an iteration of forcing notions.
We say that $\bar{\bbQ}$ is a ${\gp}$-NNR$^1_\kappa$
iteration, where $2 \le \kappa \le \aleph_0$,\footnote{we omit
$\kappa = \aleph_1$ for convenience.} when  for some $f_i,g_i \in
\bold \cF^{\gp}_{\dc}$, for $i < \ell g(\bar{\bbQ})$ (see
\ref{b999}), we have:

\begin{enumerate}
\item[$(a)$]  $\cf(\ell g({\gp})) > \ell g(\bar{\bbQ})$,

\item[$(b)$]  $\bar{\bbQ}$ is a countable support iteration of proper forcing
notions such that for each $i < \ell g(\bar{\bbQ}), \bbP_i$
adds no reals,
\footnote{this follows from other parts (close (d)), even for $\bbP_{i+1},i <
\ell g(\bar{\bbQ})$.}

\item[$(c)$]   (long properness) for each $i < \ell g(\bar{\bbQ})$,
we have $\Vdash_{\bbP_i} ``\name {\bbQ}_i$ is
$({\gp}^{\bbP_i},f_i)$-proper'',

\item[$(d)$]   ($\kappa$-anti w.d.)  if $i
  < \ell g(\bar{\bbQ})$ and $\beta \in g_i(\alpha)$, then  $\name{\bbQ}_i$ has
  $(\kappa,\alpha,\beta)$-anti w.d. above $\bbP_i$.
\end{enumerate}
\end{definition}

\begin{remark}
\label{d6}
\begin{enumerate}

\item $\kappa$ is the amount of ``$\bbD$-completeness", in other
words what versions of weak diamond we kill by our iteration.  So
the case $\kappa = \aleph_0$ is easier, and we first deal with it in Theorem \ref{d9}.

\item Note that we ask for $f_i \in \bold V$
and not a $\bbP_i$-name $\name f_i$ of such a function.
The reason is that if for $i < \ell g(\bar \bbQ), \bbP_i$ satisfies the
$\cf(\ell g({\gp}))$-c.c., then we can find $f'_i \in {\cF}^{\gp},
f'_i \ge \name f_i$. As in practice we usually have
$\cf(\ell g({\gp})) > |\bbP_\alpha|,$ there is no point at present for
$f_i$ to be a $\bbP_i$-name.

\item In clause (d) we have implicitly used:
\begin{enumerate}
\item[$(*)$]   if $\alpha \le \beta' < \beta$, then clause (d) for
$(\alpha,\beta)$ and $\kappa$ implies clause (d) for $(\alpha,\beta')$ and $\kappa$.
\end{enumerate}
This holds by clause (i) of Definition \ref{a3}.

\item We could replace $f_i$ by a club $E_i$ of $\ell g({\gp})$, letting
$f_i(\alpha) = E_i \backslash \alpha$.


\item In clause (c), for a club $C$ of $\ell g({\gp})$ we
catch our tail, that is $f_i(\alpha) \cap C = C \backslash \alpha$
for a club of $\alpha < \ell g({\gp})$.

\item In clause (d), much of the freedom/variation will be due to the
decision how ``similar" are $\langle G^\ell:\ell < k \rangle$ such that
$\bold G^{**}$ exists.  Here we demand
\mn
\begin{enumerate}
\item[$(\alpha)$]   $Y \in D^{\gp}_\beta(N_1)$.
\end{enumerate}
In \cite[Ch.VIII]{Sh:f}, it is essentially required that

\begin{enumerate}
\item[$(\beta)$]   $\bold G^0 \times \bold G^1 \times
\cdots \times \bold G^{k-1} \subseteq (\bbP_i \times \cdots
\times \bbP_i) \cap N_1$ ($k$ times) is generic over $N_1$.
\end{enumerate}

In \cite[Ch.V]{Sh:f}, it is required that
\begin{enumerate}
\item[$(\gamma)$]   the common $Y$ is a pre-determined increasing
sequence of models.
\end{enumerate}
Clause $(\beta)$ makes demand (d) in
Definition \ref{d3} easier, but the parallel of (c) is harder compared to clause
$(\alpha)$.
\end{enumerate}
\end{remark}
We now state and prove the main results of this section. We first deal with the
${\gp}$-NNR$^1_{\aleph_0}$ iterations.
\begin{theorem}
\label{d9}
Assume $\bar{\bbQ}$ is a ${\gp}$-$\text{NNR}^1_{\aleph_0}$ iteration,
${\gp}$ is a reasonable parameter and ${\gp}$ is a $\Game_f$-winner
for some $f \in {\cF}^{\gp}_{\club}$ (or at least is $\Game'_f$-non
loser).
\begin{enumerate}
\item Forcing with $\bbP_{\ell g(\bar{\bbQ})} = \Lim(\bar{\bbQ})$
does not add reals (so consequently adds no $\omega$-sequences, as we
are assuming properness).

\item If $i \le j \le \ell g(\bar{\bbQ})$, then

\begin{enumerate}
\item[$(b)'$]   $\bbP_j/\bbP_i$ is proper,

\item[$(c)'$]   $\bbP_j/\bbP_i$ is $({\gp},f_{i,j})$-proper, where
$f_{i,j} \in {\cF}^{\gp}_{\club}$ is increasing continuous and
is computable from the $f_\varepsilon \in {\cF}^{\gp}$
for $\varepsilon \in [i,j)$,

\item[$(d)'$]   we have the parallel of clause (d) in the following sense: if $i < j < \ell g(\bar{\bbQ}),$ then for some function $g \in \cF^{\gp}_{\text{cd}}$ in $\bold V$ and for all $\alpha < \ell g(\gp)$ and $\beta \in g(\alpha)$ we have $\bbP_j / \bbP_i$ is $(\aleph_0, \alpha, \beta)$-anti-w.d above $\bbP_i$.
\end{enumerate}
\end{enumerate}
\end{theorem}

\begin{proof}
The proof is by induction on $\ell g(\bar{\bbQ})$.  For notational simplicity
we assume that:

\begin{enumerate}
\item[$\boxtimes:$]  all $f_i$'s are  also in $\cF^{\gp}_{\nd}$,
so we can consider them as increasing and continuous functions from $\ell g({\gp})$ to $\ell g({\gp})$.  We also demand that the $f_{i,j}$'s are also like that, are
increasing continuous and moreover $f_{i,j}(f_{i,j}(\alpha)) = f_{i,j}
(\alpha)$, and they are $\ge f^*$ where $f^* \in {\cF}^{\gp}_{\nd}$ is
increasing continuous and ${\gp}$ is $\Game_{f^*}$-winner (or at least
$\Game'_{f^*}$-non-loser).
\end{enumerate}

\noindent
\underline{Case 1: $\ell g(\bar{\bbQ}) = 0$.}
This is trivial.

\noindent
\underline{Case 2: $\ell g(\bar{\bbQ}) = i(*) + 1$ is a successor ordinal.}
We show that items (1) and (2) are satisfied.

\underline{Clause (1)}:  $\bbP_{i(*)}$ adds no reals by the
induction hypothesis and $\Vdash_{\bbP_{i(*)}}$``$\name{\bbQ}_{i(*)}$ adds no reals'', by clause (d) in
Definition \ref{d3}, hence $\bbP_{i(*)+1} = \bbP_{i(*)} *
\name{\bbQ}_{i(*)}$ adds no reals.

\underline{Clause (2)}: We have to show that items $(b)', (c)'$ and $(d)'$ are satisfied.

\begin{enumerate}
\item[] \underline{Clause $(b)'$}: By \cite[Ch. III]{Sh:f},  $\bbP_j/\bbP_i$ is proper.

\item[] \underline{Clause $(c)'$}:
Given $i \le j \le \ell g(\bar Q)$, if $j < i(*)+1$ the conclusion follows
by the induction hypothesis.  So assume
$j =i(*) +1$.  If $i=j$, the required demand is trivial, so assume $i < j$.
If $i=i(*)$, use clause (c) of Definition \ref{d3} for $i$ to get the conclusion.
So assume that $i < i(*)$.  Let
\begin{itemize}
\item $f_{i,j,0} = f_{i(*)},$
\item $f_{i,j,m+1}
= f_{i(*)} \circ f_{i,i(*)} \circ f_{i,j,m}$,
\item
$f_{i,j}(\alpha) = \underset{m < \omega}{\sup} \,
f_{i,j,m}(\alpha)$.
\end{itemize}
Then the $f_{i, j}$'s are as required in $\boxtimes$. To
prove ``$\bbP_j/\bbP_i$ is $({\gp},f_{i,j})$-proper'', assume
that
\begin{enumerate}
\item[$(*)_1$]  $(a) \quad N \prec ({\cH}(\chi),\in)$ is countable,

\item[${{}}$]  $(b) \quad \{\bar{\bbQ},i,j,\alpha,\beta,f_{i,i(*)},f_{i(*)},
f_{i,j}\} \in N$,

\item[${{}}$]   $(c) \quad \alpha \le f_{i,j}(\alpha) \le
\beta < \ell g({\gp})$,
\
\item[${{}}$]  $(d) \quad q \in \bbP_i \text{ is }
(N,\bbP_i) \text{-generic}$,

\item[${{}}$]  $(e) \quad p \in N \cap \bbP_j,p \restriction i \le q$,

\item[${{}}$]  $(f) \quad Y \in D^{\gp}_\beta(N)$,

\item[${{}}$]  $(g) \quad q \Vdash ``Y \subseteq {\cM}_{\bbP_i}
[\name{\bold G}_{\bbP_i},N]$''.
\end{enumerate}
\mn
First we deal with version 2, and assume that ${\gp}$ is simple.
Choose $y^* \in N$ which codes enough information. Clearly
$\beta' = f_{i(*)}(\alpha)$ belongs to $N$.  So
$f_{i,i(*)}(\beta') \le \beta$, hence
by the induction hypothesis there are $q',Y'$ such that:
\begin{itemize}
\item $q \le q' \in \bbP_{i(*)},$

\item $p \restriction i(*) \le q',$
\item $q'$ is
$(N,P_{i(*)})$-generic,
\item $Y' \subseteq Y, Y'\in D^{\gp}_{\beta'}(N)$,
\item $q' \Vdash ``Y' \subseteq {\cM}_{\bbP_{i(*)}}
[\name{\bold G}_{\bbP_{i(*)}},N,y^*]"$.
\end{itemize}
Next, we apply clause (c) in the Definition \ref{d3}
for $i(*)$, so there are $q'',Y''$ such that
\begin{itemize}
\item $q' \le q'' \in \bbP_{i(*)+1} = \bbP_j,$
\item $p \le q'',$
\item $q''$ is $(N,\bbP_j)$-generic,
\item $Y'' \subseteq Y', Y'' \in D^{\gp}_\alpha(N)$,
\item $q'' \Vdash ``Y'' \subseteq {\cM}
[\name{\bold G}_{\bbP_j},N,y^*]"$.
\end{itemize}
The result follows immediately.
The proof for version 1 is similar.

\item[] \underline{Clause $(d)'$}:
Recall that we have demanded
$f_{i,j}(f_{i,j}(\alpha)) = f_{i,j}(\alpha)$ (see $\boxtimes$ at the
beginning of the proof).

Let $N_0,N_1,\alpha,\beta,i,j,p,k,q_\ell$
(for $\ell < k$), $\bold G^\ell$
(for $\ell < k$) and $\bold G^*$ be as in the assumptions of Definition \ref{d3}(d)
(see Definition \ref{b37}).

Without loss of generality $i < i(*)< j = i(*) +1$, since the other cases are trivial as in the proof of
clause $(c)'$ .
First choose $\bold G^{**} \in N_1$ for $\bold G^*,\alpha,
\beta,i,i(*)$.
For each $\ell < k$, if for some $s_\ell \in \bold G^\ell$ we have

\begin{enumerate}
\item[$(*)_2$] $\hspace{1.cm}$ $s_\ell \Vdash_{\bbP_i}$ ``\text{ there is an upper bound for } $\bold
G^{**}$ \text{in } $\bbP_{i(*)}/
\name{\bold G}_{\bbP_i}$'',
\end{enumerate}
then, as $\bold G^\ell$ is generic over
$N_0$, by  increasing $s_\ell$ if necessary,  there are $s_\ell \in \bold G^\ell$ and $r_\ell \in P_{i(*)} \cap
N_1$ such that $s_\ell$ forces that $r_\ell$ is an upper bound
for $\bold G^{**}$, and without loss of generality
$r_\ell \restriction i \le s_\ell$.  Now \wilog \,
\begin{center}
$\bold G^{\ell_1} =
\bold G^{\ell_2} \Rightarrow s_{\ell_1}
= s_{\ell_2}$
\end{center}
 and
\begin{center}
$\bold G_{\ell_1} \ne
\bold G_{\ell_2} \Rightarrow s_{\ell_1},s_{\ell_2}$ are incompatible.\footnote{see Remark \ref{remarkb37}.}
\end{center}
Now choose
$r \in \bbP_{i(*)} \cap N_1$ with domain $\subseteq i(*) \backslash i$ as
follows:
\begin{itemize}
\item $\dom(r) = \bigcup\limits_{\ell < k} \dom(r_\ell)
\backslash i$,
\item  $r(\alpha)=r_\ell(\alpha)$ if $s_\ell
\in \name{\bold G}_{\bbP_i},\ell < k$,
\item $r(\alpha)=\emptyset_{\bbP_\alpha}$
if this occurs for no $\ell$.
\end{itemize}
Renaming
$r \in N_i \cap \bbP_{i(*)},\dom(r) \subseteq i(*) \backslash i$ and
$s_\ell \in \bold G^\ell,r_\ell = s_\ell \cup r$ is above $\bold G^{**}$ in
$\bbP_{i(*)}$.  Let $\beta_\ell = f_{i,j,1 +\ell}(\alpha)$ for $\ell
\le k$.

We choose, by induction on $\ell \le k$, the objects
$Y_\ell,q'_\ell,M_\ell$ such that:

\begin{enumerate}
\item[$(*)_3$]  $(a) \quad Y_0 = Y$,

\item[${{}}$]  $(b) \quad M_0 = N_1$,

\item[${{}}$]  $(c) \quad N_0 \in M_{\ell +1}$,

\item[${{}}$]  $(d) \quad M_{\ell +1} \in M_\ell \cap {\cE}^{\gp}_{\beta_{k-\ell}}$,

\item[${{}}$]  $(e) \quad Y_{\ell +1} \subseteq Y_\ell$,

\item[${{}}$]  $(f) \quad Y_\ell \in D^{\gp}_{\beta_\ell}(M_\ell)$,

\item[${{}}$]  $(g) \quad M_{\ell +1} \in Y_\ell$,

\item[${{}}$]  $(h) \quad q_\ell \le q'_\ell \in P_{i(*)}$,

\item[${{}}$]  $(i) \quad q'_\ell \text{ is } (M_{\ell +1},\bbP_{i(*)})
\text{-generic}$,

\item[${{}}$]  $(j) \quad q'_\ell \text{ forces a value for }
\name{\bold G}_{\bbP_{i(*)}} \cap M_{\ell +1}$,

\item[${{}}$]  $(k) \quad q'_\ell \text{ is }
(N_0,\bbP_{i(*)}) \text{-generic}$,

\item[${{}}$]  $(l) \quad q'_\ell \Vdash ``Y_{\ell +1} \subseteq {\cM}_{\bbP_{i(*)}}
[\name{\bold G}_{\bbP_{i(*)}},M_\ell]"$,

\item[${{}}$]  $(m) \quad q'_\ell \restriction i \in G^\ell$.
\end{enumerate}

Now apply clause (d) of the definition for $i(*),N_0,M_k,
\langle q'_\ell:\ell < k \rangle,Y_k$ and $\bold G^{**}$ and get
$\bold G^{***}$ as required.

\end{enumerate}
\noindent
\underline{Case 3:  $\delta = \ell g(\bar{\bbQ})$ is a limit ordinal}.
We show that items (1) and (2)
are satisfied in this case as well.

\underline{Clause (1)}:  This follows from clause $(2)(d)'$ proved below.

\underline{Clause (2)}: Again, we have to check items $(b)'$, $(c)'$ and $(d)'$.  Let $f_{i,j}$ be fast enough functions.

\begin{enumerate}
\item[] \underline{Clause $(b)'$}: This is obvious.

\item[] \underline{Clause $(d)'$}:  We first prove clause $(d)'$ and later prove clause $(c)'$. As before, \wilog \, $i < j = \delta$.
Let $N_0,N_1,p,\bold G^*,\alpha,\beta,k < \aleph_0$
and $\bold G^\ell,q_\ell$ for $\ell < k$ be as in the assumptions of
clause (d) of Definition \ref{d3}.

Choose $\gamma \in N_0,
\alpha < \gamma < \beta$ such that $\gamma$ is large enough,
in particular,
$$i \le i' < j' < j \Rightarrow f_{i',j'}(\gamma)
= \gamma.$$

Let $\langle i_m:m < \omega \rangle \in N_1$ be such that
\begin{enumerate}
\item[$(*)_4$]  $(a) \quad i_0 = i$,

\item[${{}}$]  $(b) \quad i_m < i_{m+1}$,

\item[${{}}$]  $(c) \quad \underset{m < \omega}{\sup} i_m =
\sup(j \cap N_0)$.
\end{enumerate}

Choose $y^* \in N_1 \cap {\cH}(\chi_\gamma)$ coding enough information.  We choose
$$M_0, M_1, M_2, M_3, M_4 \in N_1 \cap {\cE}^{\gp}_\gamma \cap Y$$
 such that
\begin{enumerate}
\item[$(*)_5$]  $(a) \quad N_0 \in M_0 \in M_1 \in M_2 \in M_3 \in
  M_4$,

\item[${{}}$]  $(b) \quad Y \cap M_m \in D^{\gp}_\gamma(M_m) \text{ for } m < 5$.
\end{enumerate}

Choose $q'_\ell \in \bold G^\ell \cap M_4$ above $\bold G^\ell
\cap M_3$ so that $q'_\ell$ is
$(M_t,\bbP_i)$-generic for $t<4$.
Let $\langle {\cI}_m:m < \omega\rangle \in M_0$ list the dense open
subsets of $\bbP_j$ from $N_0$.
Now we shall use the diagonal argument and choose
$\name{\bold G}_{\bbP} \cap N_0, p_m \in \bbP_{i_m} \cap N_0, r_m \in
\bbP_{i_m}$.
We fulfill the above  in $M_4$, so that at the end can find a
solution in $N_1$, by using a canonical construction.

But to carry this, we need to have finitely many candidates for
$\name{\bold G}_{\bbP_{i_m}} \cap M_0$ with a common $Y_m$.
To get this in the inductive step, we need in step $m-1$ that for $M_1$ we
just have finitely many candidates for
$\name{\bold G}_{\bbP_{i_m}} \cap M_1$, and in turn to get this in
the step $m-1$, we use that in step $m-2$ for $M_2$ and
from every maximal antichain we choose a finite subset.
To get this we use that for $M_3$ we just ask $M_3
[\name{\bold G}_{\bbP_{i_{m-3}}}] \cap \bold V = M_3$.  So along the
way $N_0,M_0,M_1,M_2,M_3$ our induction demands go down, but slowly, so that
in each step $m$, advancing for say $M_0$, we have to preserve less than
really knowing $\name{\bold G}_{\bbP_{i_m}} \cap M_0$, and are
helped by our demand on $M_1$, just like in \cite[Ch. XVIII]{Sh:f}. So compared to \cite[Ch. V]{Sh:f}, we have a finite tower.

Thus we choose by induction on
$m < \omega$ the objects $r_m,\bold G^*_m,p_m,n_m$,
$\langle \bold G^\ell_m:\ell < n_m \rangle$ and $Y_m$ such that:
\begin{enumerate}
\item[$(*)_6$]   $(a) \quad r_m \in \bbP_{i_m} \cap M_4$,

\item[${{}}$]  $(b) \quad \dom(r_m) \subseteq [i,i_m)$,

\item[${{}}$]   $(c) \quad r_{m+1} \restriction i_m = r_m$,


\item[${{}}$]  $(d) \quad q'_\ell \cup r_m \in \bbP_{i_m}$
is $(M_t,\bbP_{i_m})$-generic for $t<4$,

\item[${{}}$]  $(e) \quad$ if $\ell < k,{\cJ} \subseteq \bbP_{i_m}$
is dense open and ${\cJ} \in M_2$, then  for some

$\hspace{0.8cm}$ finite
 ${\cJ} \subseteq {\cI} \cap M_2,{\cJ}$ is predense above
$q'_\ell \cup r_m$,


\item[${{}}$]  $(f) \quad n_m < \omega$, and
for $\ell < n_m, \bold G^\ell_m$ is a subset of
$\bbP_{i_m} \cap M_0$ generic over

$\hspace{0.8cm}$ $M_0$
 and $\bold G^\ell_m \in M_1$,

\item[${{}}$]  $(g) \quad \bold G^\ell_{m+1} \cap \bbP_{i_m}
\in \{\bold G^\ell_m:\ell < n_m\}$,

\item[${{}}$]   $(h) \quad n_0 =k$ and $\bold G^\ell_0 = \bold G^\ell \cap M_1$,

\item[${{}}$]  $(i) \quad q_\ell \cup r_m \Vdash
``\name{\bold G}_{\bbP_{i,m}} \cap M_1 \in \{\bold G^\ell_m:\ell < n_m\}"$,

\item[${{}}$]  $(j) \quad \bold G^*_m$ is a subset of
$\bbP_{i_m} \cap N_0$ generic over $N_0$,

\item[${{}}$]  $(k) \quad \bold G^*_m \subseteq \bold G^\ell_m$,
so that $\bold G^*_m \subseteq \bold G^*_{m+1}$ and $\bold G^*_0 = \bold G^*$,

\item[${{}}$]  $(l) \quad p_m \in \bbP_j \cap N_1,$

\item[${{}}$]   $(m) \quad p_0 = p$,

\item[${{}}$]  $(n) \quad p_m \restriction i_m \in \bold G^*_m$,

\item[${{}}$]  $(o) \quad p_m \le p_{m+1} \in {\cI}_m$,

\item[${{}}$]   $(p) \quad Y_m \subseteq {\cM}_{\bbP_{i_m}}[\bold G^\ell_m,M_0,y^*]$,

\item[${{}}$]   $(q) \quad Y_m \in D^{\gp}_\gamma(M_0)$.
\end{enumerate}
Let us now explain the induction construction.
If $m=0$, this is trivial, so suppose that it holds for $m$
and we do it for
$m+1$. This is done in several stages.

\noindent
 \underline{Stage A}:  Choosing $p_{m+1}$
is trivial, the demands are: $p_{m+1} \ge p_m,p_{m+1} \rest i_m
\in \bold G^*_m$ and $p_{m + 1} \in {\cI}_m$.

\noindent
\underline{Stage B}:  To choose $\bold G^*_{m+1}$, apply the
induction hypothesis using clause $(d)'$ of what we are proving
with $i_m,i_{m+1},\gamma,f_{i_m,i_{m+1}}(\gamma),N_0,M_1$ here
standing for $i,j,\alpha,\beta,N_0,N_1$ there.

\noindent
\underline{Stage C}:  Let $\{H^\ell_m:\ell < n_{m+1}\}$ list
the possibilities of $\name{\bold G}_{\bbP_{i_m}} \cap M_1$
(by clause (e) this exists).
Without loss of generality $H^\ell_m \cap M_0 = G^{h(\ell)}_m$, for some
function $h = h_m:n_{m+1} \rightarrow n_m$.
We choose $s^\ell_m \in \bbP_{i_{m+1}} \cap M_1$ above
$\bold G^*_{m+1}$, such that $s^\ell_m \restriction i_m \in
\bold G^{h(\ell)}_m$.  Now we repeat the argument of
the successor stage of shrinking $Y$, so we can find $t^\ell_m $ such that
\begin{itemize}
\item  $t^\ell_m \in \bbP_{i_{m+1}} \cap M_1$, above $s^\ell_m,$
\item $t^\ell_m \restriction i_m \in H^\ell_m,$
\item $t^\ell_n \Vdash
``\name{\bold G}_{\bbP_{i_{m+1}}} \cap M_0 =: \bold G^\ell_{m+1}"$,
\end{itemize}
and such that
\begin{enumerate}
\item[$(*)_7$]   $Y_{m+1} =: \bigcap\limits_{i < n_{m+1}} {\cM}
[\bold G^\ell_{m+1},M_0,y^*] \in D^{\gp}_\gamma(M_0)$.
\end{enumerate}
The rest is as in the proof of Theorem \ref{a27}.

Now, without
loss of generality the construction belongs to $N_1$.
So
$$\bold G^{**} = \{s \in \bbP_j \cap N_0:\bigvee\limits_{n < \omega}
s \le p_m\}$$
 is as required, as $q'_\ell =: q_\ell \cup
\bigcup\limits_m r_m \in \bbP_j \cap N_1$, and is above
$\bold G^{**}$ and $p \le q'_\ell$.
This finishes proving clause $(d)'$ in
the case $\ell g(\bar{\bbQ})$ is a limit ordinal.

\item[] \underline{Clause $(c)'$}:  Again, \wilog \, $i < j = \delta$.  So assume
$f_{i,j}(\alpha) \le \beta, \{i,j,\alpha,\beta\} \in N^* \in
{\cE}^{\gp}_\beta, q$ is $(N^*,\bbP_i)$-generic, $Y^* \in D^{\gp}_\beta
(N^*), q \in \bbP_\alpha$ and $q \Vdash ``Y^* \subseteq {\cM}_{\bbP_i}
[\name{\bold G}_{\bbP_i},N^*,y^*]"$ are given.  We prove the desired
conclusion by induction on $\alpha$.  For each $\alpha$, we would like to
simulate a play of $\Game_{\alpha,\beta}(N^*)$, supplying the challenger
with a strategy.  For this we apply the proof of clause $(d)'$.
Choose $N_0, N_1, M_0, \dotsc, M_4, q_0, \bold G_0, \bold G^*$
(and $k=1$) as there so that for some $\alpha' < \gamma' < \beta'$ as
there,  $\beta < \alpha'$ and
$N^*,q,Y^* \in N_0$.

During the construction,
 we demand $p_m \in N^* \cap \bbP_j$, so a generic for
$N_0$ is not necessarily created.  But still $p_m \le p_{m+1},p_m
\restriction i_m \in \bold G^*_m$.  Now $p_m$ will be played by the
chooser.   Now $g(1 + \alpha)$ will be a fixed point of
$f_{i_m,i_m + 1}$.  So we can add the demand $N^*[\bold G^*_m \cap N^*]
\cap \bold V = N^*$, i.e. $\bold G^*_m$ is generic over $N^*$ and

\begin{enumerate}
\item[$(*)_8$]  ${\cM}[\bold G^*_m \cap N^*,N^*,y^*] \in
D^{\gp}_{g(1 + \alpha)}(N^*)$.
\end{enumerate}
We will define the game so that the following are satisfied:
\begin{itemize}
\item The challenger chooses

\begin{enumerate}
\item[$(*)_9$]  $X_{m+1} = {\cM}[\bold G^*_m \cap N^*,N^*,y^*]
\cap \bigcup\limits_{\xi < g (1+ \alpha)} {\cE}^{\gp}_\xi
\cap \{M:p_m \in M\} \in N_0$.
\end{enumerate}
\item Now the chooser chooses $\alpha_m,\beta'_m$ and then the challenger chooses
$\beta_m \ge \beta'_m,f_{i_m,i_{m+1}}(\alpha)$ in $N^*_0 \cap j$ and the
chooser chooses $M^*_{m+1}$ such that $p_m \in M^*_{m+1}$.

\item Now the chooser chooses $Y_m \in D^{\gp}_{\beta_m}(M_0),Y_m \subseteq X_m \cap M_0,
Y_m \in N_0$.
\item Now we play $Z_m$ for the challenger as follows:
there is $p_{m+1} \ge p_m$ which is  $(M^*_{m+1},\bbP_{i_{m+1}})$-generic
such that
$p_{m+1} \restriction i_m \in \bold G^*_m$, $p'_{m+1}$ decides $\name{\bold G}_{\bbP_{i_{m+1}}} \cap N^*$ and
forces $Z_m = Y_m \cap {\cM}_{\bbP_{i_{m+1}}}[\name{\bold G}_{\bbP_{i_{m+1}}},
M^*_m,y^*] \in D^{\gp}_{\alpha_m}(M^*_m)$.
\end{itemize}
Let us now give the details.
Choose  $\langle i'_m:m < \omega
\rangle \in N_0$ such that $i_m \in N^*, i_0 = i, i_m < i_{m+1}$
and
$\sup\{i_m:m < \omega\} = \sup(N^* \cap j)$ and let
$\langle {\cI}'_m:m < \omega \rangle$ list the dense open subsets of
$\bbP_j$ from $N^*$.  For $\bold m < \omega$ let ${\cT}_{\bold m}$ be the
set of finite sequences ${\gx}$ from $M_4$ coding
\begin{itemize}
\item $\langle
r_{{\gx},m}:m \le \bold m \rangle,$
\item $\langle \bold G_{{\gx},m}:m \le
\bold m \rangle,$
\item $\langle p_{{\gx},m}:m \le \bold m \rangle,$
\item $\langle
n_{{\gx},m}:m \le \bold m \rangle,$
\item $\langle \bold G^\ell_{{\gx},m}:\ell \le
n_{{\gx},m},m \le \bold m \rangle,$
\item $\langle Y_{{\gx},m}:m \le \bold m
\rangle$,
\item   $(X_{{\gx},m},\alpha_{{\gx},m},\beta'_{{\gx},m},
\beta_{{\gx},m},M_{{\gx},n},y'_{{\gx},m},M'_{{\gx},m},y_{{\gx},m})$
for $m \le \bold m$,
\item  $Z_{{\gx},m}$ for $m < \bold m$,
\end{itemize}
satisfying
items (a)-(k), (m), (n), (p) and (q) from the proof of $(d)'$ above and
\begin{enumerate}
\item[$(*)_{10}$]   $(l)' \quad p_m \in \bbP_j \cap N^*$,

\item[${{}}$]   $(o)' \quad p_m \le p_{m+1} \in {\cI}'_m$,

\item[${{}}$]  $(r)' \quad r_{{\gx},m}$ is $(N^*,\bold G_{i_m})$-generic for
$m \le \bold m$,

\item[${{}}$]   $(s)' \quad \langle (X_{{\gx},m},\alpha_{{\gx},m},
\beta'_{{\gx},m},\beta_{{\gx},m},M_{{\gx},m},Y_{{\gx},m},
M'_{{\gx},m},Z_{{\gx},m'}):m \le \bold m \rangle$

$\hspace{0.8cm}$ belongs
  to $N$ and
is an initial segment of a play of the game

$\hspace{0.8cm}$ $\Game'_{\alpha,\beta}(N^*,{\gp})$ or
 just $\Game'_{\alpha,\beta}(N^*,N,
{\gp})$,\footnote{ note that in the $\bold m$-th move the challenger has
  not yet chose
$Z_{{\gx},\bold m}$, (see clause (e) of Definition \ref{b6}(1).}

\item[${{}}$]  $(t)' \quad Z_{{\gx},m} \subseteq Y_{m+1}$,

\item[${{}}$]  $(u)' \quad y_{{\gx},m}$ codes $p_m, \langle i_m:m < \omega
\rangle$,

\item[${{}}$]  $(v)' \quad f_{i_m,i_{m+1}}(\alpha_m) \le
\beta'_m$ for $m \le \bold m$.
\end{enumerate}

We let ${\gx} \triangleleft {\gy}$ to have the natural meaning for ${\gx}
\in {\cT}_{{\bold m}_1},{\gy} \in {\cT}_{{\bold m}_2},
\bold m_1 < \bold m_2$.
Note
that
\begin{enumerate}
\item[$\boxtimes_1$]    ${\cT}_{\bold m} \subseteq N$ for $\bold m < \omega$,

\item[$\boxtimes_2$]   ${\cT_0} \ne \emptyset$,

\item[$\boxtimes_3$]   if $\bold x \in {\cT}_{\bold m}$, then ${\gx}$ is an initial
segment of a play
of the game $\Game_{\alpha,\beta}(N^*,{\gp})$ (see clause $(s)'$ above).
\end{enumerate}
Now we show that $(*)_{11} \Rightarrow (*)_{12},$ where
\begin{enumerate}
\item[$(*)_{11}$]  $(a)$ \quad  $M'_{\bold m} \in Y_{{\gx},\bold m} \cap
{\cE}^{\gp}_{\alpha_{\bold m}} \cap (M_{{\gx},\bold m} \cup
\{M_{{\gx},\bold m} \cap {\cH}(\chi^{\gp}_{\bold m})\}$
satisfies $y_{\bold m},y'_{\bold m} \in$

$\hspace{0.8cm}$ $M'_{\bold m}$,

\item[${{}}$]  $(b)$ \quad $Z_{\bold m}
\subseteq {\cD}_{\alpha_{{\gx},m}}(M'_m),$

\item[${{}}$]  $(c)$ \quad $Z_{\bold m} \subseteq
Y_{\bold m}$ (hence $Z_{\bold m} \subseteq X_{{\gx},\bold m})$,

\item[${{}}$]  $(d)$ \quad $X_{\bold m +1} \in D^{\gp}_\beta(N^*) \cap X_{{\gx},\bold m}$ is such that
$Z_{\bold m} \subseteq X_{\bold m +1}$,

\item[${{}}$]  $(e)$ \quad $\alpha_{\bold m +1} \in
\alpha \cap N^*$,

\item[${{}}$]  $(f)$ \quad  $\beta'_{\bold m +1} \in \beta \cap N^* \backslash
p_{i_{\bold m +1},i_{\bold m_r}}(\alpha_{\bold m +1}),$

\item[${{}}$]  $(g)$ \quad $y'_{\bold m +1} \in
N \cap {\cH}(\chi^{\gp}_{\alpha_{\bold m +1}})$ and $y'_{\bold m + 1} \in M_{\bold m + 1}$,

\item[${{}}$]  $(h)$ \quad $\beta_{\bold m} \in \beta \cap N \backslash
\beta'_n \backslash \alpha_n$ and $M_{\bold m +1} \in
X_{\bold m + 1} \cap {\cE}^{\gp}_{\beta_{\bold m + 1}},$

\item[${{}}$]  $(i)$ \quad $y_{\bold m + 1} \in M_{\bold m + 1} \cap {\cH}(\chi^{\gp}
_{\alpha_{\bold m}})$,

\item[${{}}$]  $(j)$ \quad $Y_{\bold m + 1} \in N \cap D^{\gp}_{\beta_{\bold m + 1}}(M_{\bold m + 1})$,

\item[${{}}$]  $(k)$ \quad any
$M'_{\bold m + 1} \in Y_{\bold m + 1} \cap {\cE}^{\gp}_{\alpha_m} \cap
(M_{\bold m + 1} \cup \{M_{\bold m + 1} \cap {\cH}(\chi^{\gp}_{\alpha_m})\}$

$\hspace{0.8cm}$ satisfies $y_{\bold m + 1},y'_{\bold m +1} \in M'_{\bold m + 1}$.
\end{enumerate}
and
\begin{enumerate}
\item[$(*)_{12}$]   there is ${\gy} \in {\cT}_{\bold m +1}$
such that
 ${\gx} \triangleleft  {\gy}$
and
$(Z_{{\gy},\bold m},
X_{{\gy},\bold m + 1}$, $\alpha_{{\gy},\bold m +1},
\beta'_{{\gy},m + 1},y'_{{\gy},\bold m + 1}$,
$\beta_{{\gy},\bold m + 1},y_{{\gy},\bold m +1},
M_{{\gy},\bold m + 1},Y_{{\gy},\bold m + 1},M'_{{\gy},\bold m+1})$
is equal to $(Z_{\bold m},X_{\bold m + 1},\alpha_{\bold m + 1},
\beta'_{\bold m +1}$, $y'_{\bold m + 1},\beta_{\bold m + 1},y_{\bold m + 1}$,
$M_{\bold m + 1},Y_{\bold m + 1},M'_{\bold m + 1})$.

\end{enumerate}

To see this, note that $f_{i_{\bold m},i_{\bold m+1}}(\alpha_{\bold m}) \le \beta_{\bold m}$, hence
 $\bbP_{i_{\bold m + 1}}/\bbP_{i_{\bold m}}$ is $({\gp},\alpha_{\bold m},
\beta_{\bold m})$-proper; thus let
$\bold G_{i_{\bold m}} \subseteq \bbP_{i_{\bold m}}$ be generic
over $\bold V,r_{\bold m} \in \bold G_{i_{\bold m}}$ to the model
$M'_{{\gx},\bold m}$ and the set $Y_{{\gx},\bold m}$.
Thus we can describe a strategy for the challenger in the game
$\Game_{\alpha,\beta}(N^*,{\gp})$ (or $\Game'_{\alpha,\beta}(N^*,N,
{\gp})$) delaying his choice of $M'_{\bold m},Z_{\bold m}$ to the $(\bold m+1)$-th move, he
just chose on the side ${\gx}_{\bold m} \in {\cT}_{\bold m}$ which ``codes'' what they played
so far and preserve ${\gx}_{\bold m} \triangleleft {\gx}_{\bold m+1}$.

By $\boxtimes_3$ this is possible, all possible choices of the chooser are
allowed, that is this gives a well defined strategy for the challenger.
Now take ${\gy} \in {\cT}_{\bold m +1}$ be such that for all $\bold m$, ${\gx}_{\bold m} \triangleleft {\gy}$

As the challenger does not have a winning strategy, there is a play where the chooser wins. This gives us
$$\bold G''=\bigcup\limits_{n < \omega} \bold G^*_n \cap N^*$$
 with a bound.
Also, there  is such a choice $\langle {\gx}_{\bold m}:\bold
m < \omega \rangle$ with $\bigcup\{(M'_{{\gx}_{\bold m+1},\bold m}\} \cup
Y_{{\gx}_{\bold m+1},\bold m}:\bold m < \omega\} \in D^{\gp}_\alpha(N)$ and
$q'=\bigcup\limits_{\bold m < \omega} r_{\bold m}$ is as required.
\end{enumerate}
The proof is complete.
\end{proof}
Now we deal with the case
of ${\gp}$-NNR$^1_{\kappa}$ iteration, where  $2 \leq \kappa < \aleph_2$. The adaptation for the proof of Theorem \ref{d9} when $2 \leq \kappa < \aleph_2$ should be clear.

\begin{theorem}
\label{d20}
Assume $\bar{\bbQ}$ is a ${\gp}$-NNR$^1_{\kappa}$ iteration where  $2 \leq \kappa < \aleph_2$,
${\gp}$ is a reasonable parameter and ${\gp}$ is a $\Game_f$-winner
for some $f \in {\cF}^{\gp}_{\club}$ (or at least is $\Game'_f$-non
loser).
\begin{enumerate}
\item Forcing with $\bbP_{\ell g(\bar{\bbQ})} = \Lim(\bar{\bbQ})$
does not add reals.

\item If $i \le j \le \ell g(\bar{\bbQ})$, then

\begin{enumerate}
\item[$(b)'$]   $\bbP_j/\bbP_i$ is proper,

\item[$(c)'$]   $\bbP_j/\bbP_i$ is $({\gp},f_{i,j})$-proper, where
$f_{i,j} \in {\cF}^{\gp}_{\club}$ is increasing continuous and
is computable from the $f_\varepsilon \in {\cF}^{\gp}$
for $\varepsilon \in [i,j)$,

\item[$(d)'$]   if $i < j < \ell g(\bar{\bbQ}),$ then for some function $g \in \cF^{\gp}_{\text{cd}}$ in $\bold V$ and for all $\alpha < \ell g(\gp)$ and $\beta \in g(\alpha)$ we have $\bbP_j / \bbP_i$ is $(\aleph_0, \alpha, \beta)$-anti-w.d above $\bbP_i$.
\end{enumerate}
\end{enumerate}
\end{theorem}

\begin{PROOF}{\ref{c12}}
Similar to the proof of Theorem \ref{d9}, with some changes
as in the proof of Theorem \ref{a42}.
\end{PROOF}

\begin{remark}
\label{d23}
We may be interested in non-proper forcing notions, say
semi-proper and UP ones (see \cite[Ch. X, XI, XV]{Sh:f}).
Here the change from reasonable parameter ${\gp} =
{\gp}^V$ to ${\gp}^{\bold V[\bold G]}$ is more serious as
$\{N \cap \chi_\alpha:N \in {\cE}^{{\gp}^{\bold V[\bold G]}}\}$
is in general not equal to $\{N \cap \chi_\alpha:N \in
{\cE}^{\frak p}_\alpha\}$.  This is treated in \cite{Sh:311}.

\end{remark}

\section {Forcing axioms compatible with CH} \label{spelling}

As is well known, iteration theorems give us consistency of axioms and in this section we  present a few of such examples.
We  consider $\kappa \in \{2,\aleph_0\}$, but could also have $\kappa =
\aleph_1$ at some points.

\begin{definition}
\label{b488}
Suppose $\gp$ is an o.b. parameter. Then $\Ax^\alpha_\lambda({\gp},\kappa,0)$ means: if $\bbQ$ is
$\aleph_2$-e.c.c. ($\aleph_2$-pic if $\lambda=\aleph_2$) and an $NNR^0_\kappa$-forcing notion
for ${\gp},$  ${\cI}_\beta$ is a dense open subset of $\bbQ$, for
$\beta < \beta^* < \lambda$, and $\name S_i$ is a $\bbQ$-name of
a stationary subset of $\omega_1$, for $i < i^* < \alpha$, then
for some directed $\bold G \subseteq \bbQ$ we have:
\begin{itemize}
\item $\beta < \beta^* \Rightarrow \bold G \cap {\cI}_\beta \ne \emptyset,$
\item $i< i^* \Rightarrow \name S_i[\bold G] = \{\gamma < \omega_1:(\exists r \in \bold G)
(r \Vdash_{\bbQ} ``\gamma \in \name S_i")\}$ is a stationary
subset of $\omega_1$.
\end{itemize}
We remove $\alpha$ when $\alpha=0.$
\end{definition}
Now we use our results on preservation of being
an $\text{NNR}^0_\kappa$-forcing notion for ${\gp}$ to get the  consistency of $\Ax^\alpha_\lambda({\gp},\kappa,0)$.

\begin{lemma}
\label{b51}
\begin{enumerate}
\item If ${\gp}$ is an o.b. parameter of length $\ell g({\gp}) = \omega_1$ which is non-$\Game'_{\id}$-loser,
 and if $\bar{\bbQ}$ is a countable support iteration
such that for $\alpha < \ell g(\bar{\bbQ}),~
\Vdash_{\bbP_\alpha} ``\name{\bbQ}_\alpha$ is an
$NNR^0_\kappa$-forcing notion for ${\gp}$'', then  $\bar{\bbQ}$ is an
$NNR^0_\kappa$-iteration for ${\gp}$

\item Assume $\CH + \mu = \mu^{< \mu} \ge \lambda$. If ${\gp}$ is a
non$-\Game'_{\id}$-loser o.b. parameter, $\chi^{\gp}_0 > 2^\lambda$,
then  for some $\aleph_2$-e.c.c. ($\aleph_2$-pic, if $\lambda=\aleph_2$)
$NNR^0_\kappa$-forcing notion $\bbP$ of size
$\mu$ we have $\Vdash_{\bbP} ``\Ax_\lambda({\gp},\kappa,0)$''.

\end{enumerate}
\end{lemma}

\begin{proof}
(1). Follows from Theorem \ref{a27},

(2).  It follows using a suitable countable support iteration of length $\mu,$ forcing all instances of the axiom $\Ax_\lambda({\gp},\kappa,0)$ at some stage of the iteration. Clause (1) and Theorem \ref{ecc and pic and chain condition} guarantee that the iteration  is as required..
\end{proof}

\begin{definition}
\label{f11}
Assume $\lambda = \lambda^{< \lambda} \gg \aleph_1
  \ge \kappa \ge 2$ and $\gp$ is a reasonable parameter such that $\lambda < \chi^{\gp}_0$ and
 $\lambda \le \cf(\ell g(\gp))$. Let also
 \begin{center}
  $\bbR = \bbR_{\lambda,\gp} = (\{\bar{\bbQ}:\bar{\bbQ} \in
  \cH(\lambda)$ is a $\gp$-NNR$^0_\kappa$ iteration$\} , \leq_{\bbR})$
  \end{center}
   where
   $$\bar{\bbQ}^1
  \le_{\bbR} \bar{\bbQ}^2 \Leftrightarrow \bar{\bbQ}^1 = \bar{\bbQ}^2
  \rest \ell g(\bar{\bbQ}^1).$$

\begin{enumerate}
\item  $\name{\bbQ}$ is absolutely
   $(\lambda,\gp,\bbR)$-NNR$^0_\kappa$ forcing above $\bar{\bbQ}$,
   when

\begin{enumerate}
\item[$(a)$]  $\bar{\bbQ} \in \bbR$,

\item[$(b)$]  $\name{\bbQ}$ is a $\Lim(\bar \bbQ)$-name of a forcing
  notion from $\cH(\lambda)^{\bold V[\bold{G}_{\Lim(\bar{\bbQ}}]}$,

\item[$(c)$]  if $\bar{\bbQ} \le_{\bbR} \bar{\bbQ}^1$, then
  $\name{\bbQ}$ is  $(\bar{\bbQ}^1,\bar{\bbQ},\gp)$-NNR$^0_\kappa$,
  which means that there is $\bar{\bbQ}^2 \in \bbR,\bar{\bbQ}^1
  \le_{\bbR} \bar{\bbQ}^2$ and $\bar{\bbQ}^2_{\ell g(\bar{\bbQ}^1)}
  = \name{\bbQ}$ so $\name{\bbQ}$ is a $\Lim(\bar{\bbQ}^1)$-name $\in
  \cH(\lambda)$.
\end{enumerate}

\item  $(\name{\bbQ},\name{\bar{\cI}})$ is an absolute
$(\lambda,\gp,\bbR)$-NNR$^0_\kappa$-problem above $\bar{\bbQ}$, when
\begin{enumerate}
\item[$(a)$]  $\name{\bbQ}$ is a $\Lim(\bar{\bbQ})$-name of a forcing notion
  from $\cH(\lambda)^{\bold V[\bold{G}_{\Lim(\bar{\bbQ})}]}$,

\item[$(b)$]  $\name{\bar \cI}$ is a $\Lim(\bar{\bbQ})$-name for a sequence of $< \lambda$ subsets
  of $\name{\bbQ}$,

\item[$(c)$]  if $\bar{\bbQ} \le_{\bbR} \bar{\bbQ}^1$, then
there is $\bar{\bbQ}^2$ such that $\bar{\bbQ}^1 \le_{\bbR}
\bar{\bbQ}^2$ and $\Vdash_{\Lim(\bar{\bbQ}^2)}
``(\name{\bbQ}, \name{\bar \cI})$ is solved'', which means there is a directed
$\bold G \subseteq \name{\bbQ}$ meeting $\name \cI_\varepsilon$ for
every $\varepsilon < \ell g(\name{\bar \cI})"$.
\end{enumerate}
\end{enumerate}
\end{definition}

\begin{lemma}
\label{f2}
Suppose that $\CH$ holds,  $\lambda = \lambda^{< \lambda} \gg \aleph_1
  \ge \kappa \ge 2$ and $\gp$ is a reasonable parameter such that $\lambda < \chi^{\gp}_0$ and
 $\lambda \le \cf(\ell g(\gp))$.
Then there is a proper $\lambda$-c.c. forcing
  notion  $\bbP_*$of cardinality $\lambda$, such that
\begin{enumerate}
\item[(a)] Forcing with $\bbP_*$  adds no reals,
\sn
\item[(b)]  $\bbP_* = \Lim(\bar{\bbQ}_*),$ where $\bar{\bbQ}_*$ is a countable
support
  iteration $\langle \bbP_\alpha,\name{\bbQ}_\alpha:\alpha <
  \lambda\rangle$ with $\Vdash_{\bbP_\alpha} ``|\bbQ_\alpha| < \lambda$'',
  such that $\bar{\bbQ}_* \rest \alpha
\in \cH(\lambda)$ for $\alpha < \lambda$.
  In particular, $\bbP_* = \bigcup\limits_{\alpha < \lambda} \bbP_\alpha$,

\item[(c)]  $\bar{\bbQ}_*$ is $ \gp$-NNR$^1_\kappa$-iteration,

\item[(d)] if $\bold I$ is a dense open subset of
  $$\bbR = \bbR_{\lambda,\gp} = (\{\bar{\bbQ}:\bar{\bbQ} \in
  \cH(\lambda) \text{ is a } \gp\text{-NNR}^0_\kappa \text{ iteration}\} , \leq_{\bbR}),$$
   where $\bar{\bbQ}^1
  \le_{\bbR} \bar{\bbQ}^2 \Leftrightarrow \bar{\bbQ}^1 = \bar{\bbQ}^2
  \rest \ell g(\bar{\bbQ}^1)$,  and if $\bold I$ is
definable in $(\cH(\lambda),\in)$ from a parameter, then
$\lambda =
  \sup\{\alpha < \lambda:\bar{\bbQ}_* \rest \alpha \in \bold I\}$,

\item[(e)]  if $\bbQ \in \bold V^{\bbP_*}$ is a forcing
  notion of cardinality $\aleph_1$, so without loss of
generality whose set of elements is a subset of
  $\omega_1$, and if  $\alpha < \lambda$ is such that $\bbQ \in \bold V^{\bbP_\alpha}$
 and $\bbQ$ is absolute $(\gp \rest
  \alpha,\lambda, \bbR)$-proper, then
 for some
$\beta \in (\alpha,\lambda),\bbQ_\beta = \bbQ$,
hence $\bold V^{\bbP_*}
  \models \Ax_\lambda(\bbQ)$.

\item[(f)]  similarly for $(\bbQ,\bar{\cI})$.
\end{enumerate}
\end{lemma}

\begin{proof}
First note that  $\bbR$ is non-empty and with no maximal member. Indeed,
the iteration of length zero belongs to $\bbR$ and if $\bar{\bbQ} =
\langle \bbP_\beta,\name{\bbQ}_\beta:\beta < \alpha\rangle \in \bbR$,  we can define
$\bar{\bbQ}' \in \bbR$ above it by letting $\bar{\bbQ}' =
\langle \bbP_\beta,\name{\bbQ}_\beta:\beta \leq \alpha\rangle,$
where $\bbP_\alpha=\lim(\bar{\bbQ})$ and $\Vdash_{\bbP_\alpha}$``$\name{\bbQ}_\alpha =
({}^{\omega_1>}2,\triangleleft)$''.

Now fix $\Phi: \lambda \to \cH(\lambda)$ such that for each $x \in \cH(\lambda),$ the set
$\Phi^{-1}\{x\}$ is stationary in $\lambda.$
Let $\bar{\bbQ}_*=\langle \bbP_\alpha,\name{\bbQ}_\alpha:\alpha < \lambda\rangle$
be a countable support iteration of forcing notions, where at stage $\alpha$,
if $\Phi(\alpha)$ is a $\bbP_\alpha$-name for a forcing notion as in $(e)$,
then we let $\Vdash_{\bbP_\alpha}$``$\name{\bbQ}_\alpha=\Phi(\alpha)$''. Otherwise,
$\name{\bbQ}_\alpha$ is forced to be the trivial forcing notion.

Using the preservation theorems we have proved earlier, we can easily show that $\bbP_*=\lim(\bar{\bbQ}_*)$
is as required.
\end{proof}
It is natural to restrict ourselves to the linear case, but this is
not a real difference when we allow to change the $f$ a little.
\begin{definition}
\label{f13}
Let $\gp$ be a reasonable parameter.

\begin{enumerate}
\item $\gp$ is linear,   if  whenever $\alpha < \ell g(\gp),N \in
   \cE^{\gp}_\alpha$ and $Y \in D^{\gp}_\alpha(N)$, then there is $Z$ such
   that:
\begin{enumerate}
\item[$(a)$]  $Z \in D^{\gp}_\alpha(N)$,

\item[$(b)$]  $Z \subseteq Y$,

\item[$(c)$]  if $a \in Z \cap \cE^{\gp}_\gamma$ then $N \rest a \prec N
  \rest \cH(\chi_\gamma)$,

\item[$(d)$]  $Z$ is linear, which means

\item[${{}}$]  $(\alpha) \quad Z$ is well ordered by $\in$ (and by
  $\subseteq$),

\item[${{}}$]  $(\beta) \quad$ if $a \in Z$ then $Z \cap a \in N$ and
  $\langle N \rest a:a \in Z\rangle$ is $\subseteq$-increasing
  continuous.
\end{enumerate}

\item ${\gp}$ is linearly standard when it is standard and linear.

\item Assume  $f \in \cF^{\gp}$ and $g
   \in \cF^{\gp}$ is defined as
   $$g(\alpha) = \bigcup\{f(\beta):\beta \in
   f(\alpha)\}.$$  Let $\gq = \gp^{[f]}$ be defined as in $\gp$, except that for each
   $\alpha <\ell g(\gp)$ and $N$,
   $D^{\gq}_\alpha(N)= \{Y \in D^{\gp}_\alpha(N)$: for some $\beta
\in N \cap f(\alpha)$ there is $Z \subseteq Y,Z$ linear and $Z \in N
   \cap \cH(\chi^{\gp}_\beta)\}$.
   \end{enumerate}
\end{definition}

\begin{lemma}
\label{f15}
If $\gp$ is a reasonable parameter and $f, g$ and $\gq$ are as in Definition \ref{f13}(3),
then
\begin{enumerate}
\item[$(a)$] $\gq$ is a reasonable parameter,

\item[$(b)$]  if $f$ is increasing continuous then so is $g$,

\item[$(c)$]  if $f(\alpha) = \alpha$ then $g(\alpha) = \alpha$ and
  $D^{\gq}_\alpha(N) = D^{\gp}_\alpha(N)$.
\end{enumerate}
\end{lemma}
\begin{proof}
Straightforward.
\end{proof}
\begin{lemma}
\label{f17}
Assume $\gp$ is a reasonable parameter, $f \in \cF^{\gp}, \lambda=\lambda^{<\lambda}$ is large enough regular and $\bbR=\bbR_{\lambda, \gp}$.
Let $\bar{\bbQ} \in \bbR$ and $\bbP = \Lim(\bar{\bbQ})$.
\begin{enumerate}
\item If $\gp$ is linear and $\name{\bbQ}$ is a $\bbP$-name of a
($<^+\omega_1)$-proper forcing notion from $\cH(\lambda)$ and
$f \in \cF^{\gp}_{\dc}$, then, $\Vdash_{\bbP} ``\name{\bbQ}$ is
  $(\gp^{\bbP}, f)$-proper''.

 \item   if $\name{\bbQ}$  satisfies the $\kappa$-completeness
  system $\bbD \in \bold V$ over $\bbP$ and $\bbP$ forces it is an NNR proper forcing notion, then for some
  $f \in \cF^{\gp}_{\dc}$, we have
$(\bbP,\name{\bbQ})$ is $(\kappa,f)$-anti
  w.d. (see Definition \ref{b37}),

\end{enumerate}
\end{lemma}
\begin{proof}
(1). The proof is essentially the same as  the proof of  Lemma \ref{a48}.

(2). Define the function $f$ such that for each $\alpha < \ell g({\gp}),$ $f(\alpha)=[\theta_\alpha, \ell g({\gp}))$,
 where $\theta_\alpha \geq \alpha$ is large enough so that $\cH(\chi_{\theta_\alpha})$ contains all the relevant information
  and the cardinal
  $\theta$ from Definition \ref{completeness system} witnessing $\bbD$ is a completeness system is below $\chi_{\theta_\alpha}$.

  Suppose $\alpha < \ell g({\gp})$, $\beta \in f(\alpha)$ and suppose that $N_0, N_1, n, \langle p_\ell: \ell < n \rangle, \langle \bold G^\ell: \ell < n \rangle, \bold G^*, Y$ and $\name{q}$ are as in Definition \ref{b37}(1)(A). Without loss of generality, we can assume that the $p_\ell$'s are pairwise incompatible.

Pick some $M \in Y,$ so that for each $\ell < n, (\bold G^\ell \cap M)$ is $(\bbP \cap M)$-generic over $M$. Fix some $\ell < n.$ Consider the pair $(M[\bold G^\ell], \name{q}[\bold G^*]) \in \dom(\bbD)$.
By the assumption, there are $p'_\ell \geq p^\ell$ and $\name{\bold H^\ell}$ such that
\[
p'_\ell \Vdash \text{``} \name{\bold H^\ell} \text{~is in~}\Gen^+(M[\bold G^\ell], \name{\bbQ}, \name{q}) \text{''}.
\]
By extending $p'_\ell$ if necessary, we can assume that for some $\name{q}_\ell$ we also have
\[
p'_\ell \Vdash \text{``} \name{q_\ell} \geq \name{q} \text{~is an upper bound for~}\name{\bold H^\ell} \text{''},
\]
Set $\bold K^\ell = \bold G^\ell \ast \name{\bold H^\ell}$. As the iteration $\bbP \ast \name{\bbQ}$ does not add any new
$\omega$-sequences of elements of $\bold V$, and by our choice of $\bold G^*$, by extending $(p'_\ell, \name{q}_\ell),$ we may assume that for some
fixed $\bold G^{**}$ and for all $\ell < n,$ we have
\[
(p'_\ell, \name{q}_\ell) \Vdash \text{``} \bold K^\ell \cap N_0=\bold G^{**}          \text{''}.
\]
It follows that $\bold G^{**} \in \Gen(N_0, \bbP \ast \name{\bbQ})$.
Let $\name{q}'$ be such that for all $\ell < n, p'_\ell \Vdash$``$\name{q}'=\name{q}_\ell$''.
It is clear that $\langle p'_\ell: \ell < n \rangle, \name{q}'$ and $\bold G^{**}$ are as required by Definition \ref{b37}(1)(B).
\end{proof}
For the rest of this section, assume that $\CH$ holds, $\lambda=\lambda^{<\lambda}\gg \aleph_1 \geq \kappa \geq 2$ and $\gp$ is a reasonable parameter
such that $\lambda < \chi^{\gp}_0$ and $\lambda \leq \cf(\ell g(\gp))$. Let $\bbP_*$ be as in Lemma \ref{f2}.
\begin{lemma}
\label{f22}
Under the above assumptions, if $\bbQ \in \bold{V}^{\bbP_*}$ is one  of the following forcing notions, then it satisfies clause (e) of Lemma \ref{f2},
in particular $\bold V^{\bbP_*}
  \models \Ax_\lambda(\bbQ)$.
\begin{enumerate}
\item[$(a)$]  $\bbQ=\bbQ_{\bar{C}, \bar{u}}$ is as in Definition \ref{c3}(1),
\item[$(b)$]  $\name{\bbQ}=\bbQ_{\bar C}$ is as in Definition \ref{c9}, for $\bar C = \langle
  C_\delta:\delta < \omega_1$ limit$\rangle,$ where for each limit ordinal $\delta, \otp(C_\delta) = \omega.$
\sn
\item[$(c)$]   $\name{\bbQ}=\bbQ_{\bar C}$ is as in Definition \ref{c9}, for $\bar C = \langle
C_\delta:\delta < \omega_1$ limit$\rangle$, where  for some countable ordinal
$\gamma(*), \delta < \omega_1 \Rightarrow \otp(C_\delta) \le
\omega^{\gamma(*)}$.
\end{enumerate}
\end{lemma}
\begin{proof}
The $(\gp,f)$-properness of $\bbQ_{\bar{C}, \bar{u}}$ follows from Lemmas  \ref{c6} and \ref{f17}.
The $(\gp,f)$-properness of $\bbQ_{\bar{C}}$ follows from Lemma \ref{c12} for $\bar{C}$
as in (b) and from Lemma \ref{c18} for $\bar{C}$ as in (c).
The $(\gp,g)$-anti w.d. is
straightforward.
\end{proof}

Given an Aronszajn tree $T$, let  $\bbQ_T $ be the forcing notion of \cite[Ch. V, Definition 6.5]{Sh:f}.
Let also  $\bar{\cI}_T = \langle \cI_{T,\alpha}:\alpha
  < \omega_1\rangle$ where
  $$\cI_{T, \alpha} = \{(f, C, \Psi) \in \bbQ_T:T_{\le \alpha}
  \subseteq \dom(f)\}.$$

\begin{lemma}
\label{f25}
Under the above assumption, we have the following:
\begin{enumerate}
\item  If $\name T$ is a $\bbP_*$-name of an Aronszajn tree, then  the pair
   $(\name{\bbQ}_{\name T},\name{\bar\cI_T})$ is an absolute $(\lambda, \gp, \bbR)$-NNR$^0_{\aleph_1}$ problem over $\bbP_*$.

\item  If $\lambda$ is strongly inaccessible\footnote{We may
  avoid this, if we use iterations as in \cite[Ch. VIII]{Sh:f},
  i.e. $\bar{\bbQ} \in \bbR$ is only a class of $(\cH(\lambda),\in)$,
  satisfying a strong version of $\lambda$-c.c., so $\lambda =
  \aleph_2$ is sufficient}, then  every Aronszajn tree is special.
   \end{enumerate}
\end{lemma}
\begin{proof}
(1) follows from \cite[Ch. V, Theorem 6.1]{Sh:f} and Lemma \ref{f17}.

(2) is clear, as for any Aronszajn tree $T$, by (1), the forcing notion $\bbQ_T$ satisfies
clause (e) of Lemma \ref{f2}.
\end{proof}

\section {On
Moore's question}
\label{moore}

In this section we answer a question of Justin Moore about the consistency of  strong failure of club guessing sequences with CH.
\begin{definition}
\label{f28}
 Let {\bf cd}:$\cH(\aleph_1) \rightarrow \omega_1$ be one-to-one.
We say $E$ solves {\bf cd}, when \ $E$ is a club of $\omega_1$ and for
every $\alpha \in E$ we have {\bf cd}$(E \cap (\alpha +1)) < \min(E
\backslash (\alpha +1))$.
\end{definition}
Justin Moore asked the following question.
\begin{question}
\label{f31}
Is the following consistent:
\begin{enumerate}
\item[$(a)$]  $\CH$,
\item[$(b)$]  for every one-to-one function
{\bf cd} from $\cH(\aleph_1)$ to $\omega_1$, there is some $E$ which solves it.
\item[$(c)$]  if $\bar C = \langle C_\delta:\delta < \omega_1$
  limit$\rangle$, where $C_\delta \subseteq \delta =
  \sup(C_\delta)$ and $\otp(C_\delta) = \omega$, then for some club $E$ of
  $\omega_1, (\forall \delta)(\delta > \sup(C_\delta \cap E))$.
\end{enumerate}
\end{question}
We give a positive answer to the above  question by proving the following theorem.
\begin{theorem}
\label{f33}
Suppose $\CH$ holds, $\lambda=\lambda^{<\lambda}\gg \aleph_1 \geq \kappa \geq 2$ and $\gp$ is a reasonable parameter
such that $\lambda < \chi^{\gp}_0$ and $\lambda \leq \cf(\ell g(\gp))$. Let $\bbP_*$ be as in Lemma \ref{f2} and set
$\bold V_1 = \bold V^{\bbP_*}$.  Then $\bold V_1$ satisfies the requirements of Question \ref{f31}.
\end{theorem}
\begin{proof}
In $\bold V_1, \CH$ holds by Lemma \ref{f2}. Clause (c) of
\ref{f28} holds by Lemma \ref{f22}(b).
To show that clause (b) of \ref{f31} is satisfied, let {\bf cd}:$\cH(\aleph_1) \rightarrow \omega_1$ be a one-to-one
function.
Define $\bbQ = \bbQ_{\text{\bf cd}}$ as follows:
\begin{enumerate}
\item[$(a)$]  $p \in \bbQ$ iff $p$ is a closed bounded subset of
  $\omega_1$ satisfying
  $$(\forall \alpha \in p)[\alpha \ne \max(p)
  \Rightarrow \text{\bf cd}(p \cap (\alpha +1)) < \min(p \backslash (\alpha
  +1))].$$

\item[$(b)$]  $p \le_{\bbQ} q \iff ~p,q \in \bbQ$ and $p$ is an
  initial segment of $q$.
\end{enumerate}
It is easily seen that
$\bbQ_{\text{\bf cd}}$ is $(<^+ \omega_1)$-proper and that it satisfies clause (e)
of Lemma \ref{f2}. The result follows immediately.
\end{proof}

\end{document}